\documentclass[leqno,11pt]{amsart}
\usepackage{amsmath,amscd,amsthm,amsxtra}
\usepackage{epsfig,graphics,color,colortbl}
\usepackage{amssymb,latexsym}
\usepackage{mathabx}
\usepackage{mathrsfs,eucal,upgreek}
\usepackage[poly,all]{xy}
\usepackage{hyperref}
\usepackage{tikz-cd}
\usepackage{arydshln}
\usepackage{ytableau}
\usepackage{adjustbox}
\usepackage{centernot}
\usepackage{mathtools}
\usepackage{stmaryrd}
\usepackage{youngtab}
\usepackage{caption}
\usepackage{dirtytalk}
\usepackage[draft]{todonotes}

\usepackage{subfiles}

\newtheorem{thm}{\bf Theorem}[section]
\newtheorem{df}[thm]{\bf Definition}
\newtheorem{prop}[thm]{\bf Proposition}
\newtheorem{cor}[thm]{\bf Corollary}
\newtheorem{lem}[thm]{\bf Lemma}
\newtheorem{rem}[thm]{\bf Remark}
\newtheorem{ex}[thm]{\bf Example}

\numberwithin{equation}{section}

\newcommand{\B}{\mathbf{B}}

\newcommand{\cP}{\mathscr{P}}
\newcommand{\pf}{\noindent{\bfseries Proof. }}
\newcommand{\ov}{\overline}

\newcommand{\N}{\mathbb{N}}

\newcommand{\Z}{\mathbb{Z}}
\newcommand{\C}{\mathbb{C}}

\newcommand{\te}{\widetilde{e}}
\newcommand{\tf}{\widetilde{f}}
\newcommand{\g}{\mathfrak{g}}
\newcommand{\td}{\widetilde}

\newcommand{\mc}{\mathcal}
\newcommand{\mf}{\mathfrak}
\newcommand{\La}{\Lambda}

\newcommand{\la}{\lambda}

\newcommand{\blue}[1]{{\color{blue}#1}}
\newcommand{\red}[1]{{\color{red}#1}}

\newcommand{\tl}[1]{\substack{\scalebox{0.75}{#1}}}

\makeatletter
\newcommand{\oset}[3][-0.4ex]{%
  \mathrel{\mathop{#3}\limits^{
    \vbox to#1{\kern-2\ex@
    \hbox{$\scriptscriptstyle#2$}\vss}}}}
\makeatother

\begin{document}
\title[Branching rule from ${\rm GL}_n$ to ${\rm O}_n$]
{Flagged Littlewood-Richardson tableaux and branching rule for classical groups}

\author{IL-SEUNG JANG}

\address{Department of Mathematical Sciences, Seoul National University, Seoul 08826, Korea}
\email{is\_jang@snu.ac.kr}

\author{JAE-HOON KWON}

\address{Department of Mathematical Sciences and RIM, Seoul National University, Seoul 08826, Korea}
\email{jaehoonkw@snu.ac.kr}


\thanks{This work is supported by the National Research Foundation of Korea(NRF) grant funded by the Korea government(MSIT) (No.\,2019R1A2C108483311 and 2020R1A5A1016126).}

\begin{abstract}
We give a new formula for the branching rule from ${\rm GL}_n$ to ${\rm O}_n$ generalizing the Littlewood's restriction formula.
The formula is given in terms of Littlewood-Richardson tableaux with certain flag conditions which vanish in a stable range. As an application, we give a combinatorial formula for the Lusztig $t$-weight multiplicity $K_{\mu 0}(t)$ of type $B_n$ and $D_n$ with highest weight $\mu$ and weight $0$.

\end{abstract}

\maketitle
\setcounter{tocdepth}{1}

\section{Introduction}
Let $V^{\la}_{{\rm GL}_n}$ denote a complex finite-dimensional irreducible representation of the complex general linear group ${\rm GL}_n$ parametrized by a partition $\la$ of length $\ell(\la)\leq n$. 
Suppose that ${\rm G}_n$ is a closed subgroup ${\rm Sp}_n$ or ${\rm O}_n$, where $n$ is even for ${\rm G}_n={\rm Sp}_n$. 
Let $V^{\mu}_{\rm G_n}$ be a finite-dimensional irreducible ${\rm G}_n$-module parametrized by a partition $\mu$ with $\ell(\mu)\leq n/2$ for ${\rm G}_n={\rm Sp}_n$, and by a partition $\mu$ with $\ell(\mu)\leq n$ and $\mu'_1+\mu'_2\leq n$ for ${\rm G}_n={\rm O}_n$. Here $\mu'=(\mu'_i)_{i\geq 1}$ is the conjugate partition of $\mu$.  

Let 
\begin{equation}\label{eq:branching mult}
\left[V^{\la}_{{\rm GL}_n} : V^{\mu}_{\rm G_n} \right] = \dim {\rm Hom}_{{\rm G}_n}\left(V^{\mu}_{\rm G_n},V^{\la}_{{\rm GL}_n}\right)
\end{equation}
denote the multiplicity of $V^{\mu}_{\rm G_n}$ in $V^{\la}_{{\rm GL}_n}$.
In \cite{Lw-1,Lw-2}, Littlewood showed that if $\ell(\la)\leq n/2$, then 
\begin{equation}\label{eq:Littlewood restriction}
\begin{split}
\left[V^{\la}_{{\rm GL}_n} : V^{\mu}_{{\rm Sp}_n} \right]&=\sum_{\delta\in \cP^{(2)}}c^\la_{\delta'\mu}, \quad \
\left[V^{\la}_{{\rm GL}_n} : V^{\mu}_{{\rm O}_n} \right]=\sum_{\delta\in \cP^{(2)}}c^\la_{\delta\mu},
\end{split}
\end{equation}
where $c^{\alpha}_{\beta\gamma}$ is the Littlewood-Richardson coefficient corresponding to partitions  $\alpha,\beta,\gamma$, and $\cP^{(2)}$ denotes the set of partition with even parts. There has been numerous works on extending the Littlewood's restriction rules \eqref{eq:Littlewood restriction} for arbitrary $\la$ with $\ell(\la)\leq n$ (see \cite{EW,HTW} and also the references therein), but most of which are obtained in an algebraic way and hence given not in a subtraction-free way.

In \cite{Su}, Sundaram gave a beautiful combinatorial formula for \eqref{eq:branching mult} when ${\rm G}_n={\rm Sp}_n$, as the sum of the numbers of Littlewood-Richardson (LR) tableaux of shape $\la/\delta'$ with content $\mu$ satisfying certain constraints on their entries, which vanish in a stable range $\ell(\la)\leq n/2$. Recently, based on a certain combinatorial model for classical crystals called spinor model introduced by one of the authors \cite{K15,K18-3}, Lecouvey and Lenart obtained another formula for \eqref{eq:branching mult} when ${\rm G}_n={\rm Sp}_n$ in terms of LR tableaux with some flag conditions on their companion tableaux \cite{LL}. A conjectural bijection between these two formulas is also suggested in \cite{LL}. On the other hand, no orthogonal analogue of these formula has been known so far.

The goal of this paper is to give a combinatorial formula for
\eqref{eq:branching mult} when ${\rm G}_n={\rm O}_n$ for arbitrary $\la$ and $\mu$ in terms of LR tableaux with certain flag conditions on their companion tableaux which vanish in a stable range $\ell(\la)\leq n/2$.

Let us state our result more precisely when $n-2\mu'_1\geq 0$ for simplicity since the result for $n-2\mu'_1< 0$ is similar. 
Let ${\texttt {LR}}^{\la}_{\delta\mu^\pi}$ be the set of LR tableaux of shape $\la/\delta$ with content $\mu^\pi$, where $\mu^\pi$ is the skew Young diagram obtained by $180^\circ$-rotation of $\mu$. 
For $U\in {\texttt {LR}}^{\la}_{\delta\mu^\pi}$, let $T$ be the companion tableau of $U$ which is of shape $\mu^{\pi}$  (see Section \ref{subsec:notations}).
For $1 \le i \le \mu_1'$ and $1 \le j \le \mu_2'$, let
\begin{equation*}
\begin{split}
	\sigma_i & = \textrm{the $i$-th entry in the rightmost column of $T$ from bottom}\,, \\
	\tau_j & = \textrm{the $j$-th entry in the second rightmost column of $T$ from bottom}\,, \\
	m_i & = \min \left\{ n-\sigma_i+1, 2i-1 \right\}, \\
	n_j & = \textrm{the $j$-th smallest number in $\{\,j+1,\dots,n\,\} \setminus \{\,m_{j+1},\dots,m_{\mu'_1}\,\}$}\,,
\end{split}
\end{equation*}
(see Example \ref{ex:example for LR undervar}).
\,The main result in this paper is as follows (Theorem \ref{thm:non-levi branching}):
\begin{thm} \label{thm:main intro_revision}
We have
\begin{equation*}\label{eq:main result}
\left[V^{\la}_{{\rm GL}_n} : V^{\mu}_{{\rm O}_n} \right]=\sum_{\delta\in \cP^{(2)}}\underline{c}^\la_{\delta\mu}\,\,,
\end{equation*}
where $\underline{c}^\la_{\delta\mu}$ is the number of $U\in {\texttt {\em LR}}^{\la}_{\delta\mu^\pi}$ such that
\begin{equation} \label{eq:bound condition in main result}
	\quad \quad \quad \quad \tau_j + n_j \le n+1 \quad \quad \textrm{for \,\,$1 \le j \le \mu_2'$}\,.
\end{equation}
\end{thm}
\vskip 1mm
The formula of \eqref{eq:branching mult} for ${\rm G}_n = {\rm O}_n$ in terms of LR tableaux in ${\texttt {LR}}^{\la}_{\delta\mu^\pi}$ (not their companion tableaux) is also given in Corollary \ref{cor:LR version}.

The branching multiplicity \eqref{eq:branching mult} is equal to the one from $D_\infty$ to $A_{+\infty}$ from a viewpoint of Howe duality on a Fock space \cite{Wa}. 
To describe this multiplicity, we use the spinor model for crystal graphs of integrable highest weight modules of type $D_\infty$ \cite{K16}.
Unlike the case of ${\rm Sp}_n$ \cite{LL}, we have to develop in addition a non-trivial combinatorial algorithm on the spinor model of type $D$ called {\em separation} in order to have a description of branching multiplicity in terms of LR tableaux satisfying the condition \eqref{eq:bound condition in main result} for $\underline{c}^{\la}_{\delta\mu}$. 
This is a key ingredient in the proof of Theorem \ref{thm:main intro_revision}. 
We can also recover the formula \eqref{eq:Littlewood restriction} in the stable range directly from the above formula (see Corollary \ref{cor:littlewood}).

As an interesting byproduct, we obtain a new combinatorial realization for the Lusztig $t$-weight multiplicity $K_{\mu 0}(t)$ of type $B_n$ and $D_n$ with highest weight $\mu$ and weight $0$ or generalized exponents (Theorem \ref{thm:Kostka-Foulkes for BD}). This gives an orthogonal analogue of the result for type $C_n$ in \cite{LL}.

The paper is organized as follows. In Section \ref{sec:spinor} we review the spinor model for crystals of integrable highest weight modules of type $D_\infty$. In Section \ref{sec:separation}, we develop a separation algorithm for the spinor model of type $D_\infty$. In Section \ref{sec:branching}, we use the separation algorithm to derive a branching formula in Theorem \ref{thm:main1}, and show that it is equivalent to Theorem \ref{thm:main intro_revision}
(see Theorem \ref{thm:main result-flagged}).  In Section \ref{sec:genexp}, we give a combinatorial formula for the generalized exponents of type $B_n$ and $D_n$ following the idea in \cite{LL} for type $C_n$. Finally, we give a proof of Theorem \ref{thm:main1} in Section \ref{sec:proof of main}.\vskip 2mm

{\bf Acknowledgement} The authors would like to thank C. Lenart and C. Lecouvey for their interest in this work and helpful discussion on generalized exponents.
Also they wish to thank the anonymous referees for careful reading and helpful comments on the manuscript.


\section{Spinor model} \label{sec:spinor}

\subsection{Crystals} \label{subsec:crystals}
Let us give a brief review on crystals (we refer the reader to \cite{HK,Kas91,Kas95} for more details). Let $\g$ be the Kac-Moody algebra associated to a symmetrizable generalized Cartan matrix $A =(a_{ij})_{i,j\in I}$ indexed by a set $I$.
Let $P^\vee$ be the dual weight lattice, $P = {\rm Hom}_\Z( P^\vee,\Z)$ the weight lattice, $\Pi^\vee=\{\,h_i\,|\,i\in I\,\}\subset P^\vee$ the set of simple coroots, and $\Pi=\{\,\alpha_i\,|\,i\in I\,\}\subset P$ the set of simple roots of $\g$ such that $\langle \alpha_j,h_i\rangle=a_{ij}$ for $i,j\in I$. Let $P_+$ be the set of integral dominant weights. 

A {\it $\g$-crystal} (or simply a {\it crystal} if there is no confusion on $\g$) is a set $B$ together with the maps ${\rm wt} : B \rightarrow P$, $\varepsilon_i, \varphi_i: B \rightarrow \mathbb{Z}\cup\{-\infty\}$ and $\te_i, \tf_i: B \rightarrow B\cup\{{\bf 0}\}$ for $i\in I$ satisfying certain axioms. 
For $\Lambda\in P_+$, we denote by $B(\Lambda)$ the crystal associated to an irreducible highest weight $U_q(\g)$-module with highest weight $\Lambda$. 

Let $B_1$ and $B_2$ be crystals. A {\em tensor product $B_1 \otimes B_2$} is defined to be $B_1 \times B_2$ as a set with elements denoted by $b_1 \otimes b_2$, where
\begin{equation*}
\begin{split}
	& {\rm wt}(b_1 \otimes b_2) = {\rm wt}(b_1) + {\rm wt}(b_2), \\
	& \varepsilon_{i}(b_{1} \otimes b_{2}) = \max\{ \varepsilon_{i}(b_{1}), \varepsilon_{i}(b_{2})-\langle {\rm wt}(b_{1}), h_i \rangle \}, \\
	& \varphi_{i}(b_{1} \otimes b_{2}) = \max\{ \varphi_{i}(b_{1})+\langle {\rm wt}(b_2), h_i \rangle, \varphi_{i}(b_{2}) \}, \\
\end{split}	
\end{equation*}
{\allowdisplaybreaks
\begin{equation} \label{eq:tensor_product_rule}
\begin{split}
	& \te_{i}(b_{1} \otimes b_{2}) = \left\{ \begin{array}{cc} \te_{i}b_{1} \otimes b_{2} & \textrm{if} \ \varphi_{i}(b_{1}) \ge \varepsilon_{i}(b_{2}), \\ b_{1} \otimes \te_{i}b_{2} & \textrm{if} \ \varphi_{i}(b_{1}) < \varepsilon_{i}(b_{2}), \end{array} \right. \\
	& \tf_{i}(b_{1} \otimes b_{2}) = \left\{ \begin{array}{cc} \tf_{i}b_{1} \otimes b_{2} & \textrm{if} \ \varphi_{i}(b_{1}) > \epsilon_{i}(b_{2}), \\ b_{1} \otimes \tf_{i}b_{2} & \textrm{if} \ \varphi_{i}(b_{1}) \le \epsilon_{i}(b_{2}),\end{array} \right.
\end{split}	
\end{equation}}

for $i \in I$. Here, we assume that $\textbf{0} \otimes b_{2} = b_{1} \otimes \textbf{0} = \textbf{0}$. Then $B_1 \otimes B_2$ is a crystal.
For $b_1\in B_1$ and $b_2\in B_2$, we say that {\em $b_1$ is equivalent to $b_2$} if there is an isomorphism of crystals $\psi : C(b_1) \longrightarrow C(b_2)$ such that $\psi(b_1)=b_2$, where $C(b_i)$ is the connected component of $b_i$ in $B_i$ for $i=1,2$, and write $b_1\equiv_{\mf g} b_2$ or simply $b_1\equiv b_2$ if there is no confusion.

\subsection{Littlewood-Richardson tableaux} \label{subsec:notations}

In this subsection, we review some combinatorics related to the Littlewood-Richardson tableaux (LR tableaux, for short) (see \cite{Ful}).

Let $\Z_+$ denote the set of non-negative integers.
Let $\cP$ be the set of partitions or Young diagrams. 
We let $\cP_{\ell}=\{\,\la\in\cP\,|\,\ell(\la)\leq \ell\,\}$ for $\ell\geq 1$, 
where $\ell(\la)$ is the length of $\la$,
let $\cP^{(2)}=\{\,\la\in\cP\,|\,\la=(\la_i)_{i \geq 1}, \la_i\in 2\Z_+ \ (i\geq 1)\,\}$, 
and let $\cP^{(1,1)}=\{\,\la'\,|\,\la\in \cP^{(2)}\,\}$, 
where $\la'$ is the conjugate of $\la$.
Put $\cP^{(2,2)}=\cP^{(1,1)}\cap \cP^{(2)}$.
For $\Diamond\in \{(1,1), (2), (2,2)\}$ and $\ell\geq 1$, 
we put $\cP^{\Diamond}_\ell=\cP^\Diamond\cap \cP_\ell$.

For a skew Young diagram $\lambda/\mu$, we define ${SST}(\lambda/\mu)$ to be the set of semistandard tableaux of shape $\lambda/\mu$ with entries in $\mathbb{N}$. For $T\in SST(\lambda/\mu)$, let 
$w(T)$ be the word given by reading the entries of $T$ column by column from right to left and from top to bottom in each column, and let ${\rm sh}(T)$ denote the shape of $T$.  

Let $\la\in \cP$ be given. For $T\in SST(\lambda)$ and $a\in \N$, we denote by $a \rightarrow T$ the tableau obtained by the column insertion of $a$ into $T$.
For a word $w=w_1\dots w_r$, we define $(w\rightarrow T)=(w_r\rightarrow (\dots \rightarrow(w_1\rightarrow T)))$. For a semistandard tableau $S$, we define $(S\rightarrow T)=(w(S)\rightarrow T)$. 

Let $\la^\pi$ denote the skew Young diagram obtained from $\la$ by $180^\circ$ rotation. Let $H_\la$ and  $H_{\la^\pi}$ be the tableaux in $SST(\la)$ and $SST(\la^\pi)$, respectively, where the $i$-th entry from the top in each column is filled with $i$ for $i\geq 1$.

For $\la,\mu,\nu\in\cP$, let $\texttt{LR}^\la_{\mu \nu}$ be the set of Littlewood-Richardson tableaux $S$ of shape $\la/\mu$ with content $\nu$. There is a natural bijection from $\texttt{LR}^\la_{\mu \nu}$ to the set of $T\in SST(\nu)$ such that $(T\rightarrow H_\mu)=H_\la$, where each $i$ in the $j$th row of $S\in \texttt{LR}^\la_{\mu \nu}$ corresponds to $j$ in the $i$th row of $T$. We call such $T$ a companion tableau of $S\in \texttt{LR}^\la_{\mu \nu}$.

We also need the following anti-version of LR tableaux which will be used frequently in this paper.

\begin{df} \label{df:LR tableaux anti-dominated ver}
{\em
We define $\texttt{LR}^\la_{\mu\nu^\pi}$ to be the set of $S\in SST(\la/\mu)$ with content $\nu^\pi$ such that $w(T)=w_1\dots w_r$ is an {\em anti-lattice word}, that is,
the number of $i$ in $w_k\dots w_r$ is greater than or equal to that of $i-1$ for each $k\geq 1$ and $1<i\leq \ell(\nu)$. 
}
\end{df}
Let us call $S \in \texttt{LR}^\la_{\mu\nu^\pi}$ a Littlewood-Richardson tableau of shape $\la/\mu$ with content $\nu^\pi$.
As in case of $\texttt{LR}^\la_{\mu\nu}$, the map from $S\in \texttt{LR}^\la_{\mu\nu^\pi}$ to its companion tableau gives a natural bijection from $\texttt{LR}^\la_{\mu\nu^\pi}$ to the set of $T\in SST(\nu^\pi)$ such that $(T\rightarrow H_\mu)=H_\la$.

From now on, all the LR tableaux are assumed to be the corresponding companion tableaux unless otherwise specified. 

Finally, let us recall a bijection 
\begin{equation}\label{conjugation of LR}
\xymatrixcolsep{3pc}\xymatrixrowsep{0pc}\xymatrix{
\psi : {\texttt {LR}}^{\la'}_{\mu'\nu'}  \ar@{->}[r] & {\texttt {LR}}^{\la}_{\mu\nu^\pi},}
\end{equation}
which may be viewed as an analogue of Halon-Sundaram's bijection \cite{HS} for an anti-dominant content (cf. \cite[Appendix A.3]{Ful}, \cite[Remark 6.3]{LL} and references therein).

Let $S\in {\texttt {LR}}^{\la'}_{\mu'\nu'}$ be given, that is, $(S\rightarrow H_{\mu'})=H_{\la'}$.  
Let $S^1,\dots,S^p$ denote the columns of $S$ enumerated from the right. 
For $1\le i\le p$, let $H^i=(S^i \rightarrow H^{i-1})$ with $H^0=H_{\mu'}$  
so that $H^{p}=H_{\la'}$. 
Define $Q(S\rightarrow H_{\mu'})\in SST(\la/\mu)$ to be the tableau such that the horizontal strip ${\rm sh}(H^i)'/{\rm sh}(H^{i-1})'$ is filled with $1\leq i\leq p$. 

On the other hand, let $U\in {\texttt {LR}}^{\la}_{\mu\nu^\pi}$ be given, that is,  
${\rm sh}(U\rightarrow H_{\mu})=H_{\la}$. Let $U_i$ denote the $i$-th row of $U$ from the top, and let $H_i = (U_i\rightarrow H_{i-1})$ with $H_0=H_{\mu}$ for $1\le i\le p$. 
Define $Q(U\rightarrow H_{\mu})$ to be tableau such that the horizontal strip ${\rm sh}(H_i)/{\rm sh}(H_{i-1})$ is filled with $1\le i\le p$. 

Then for each $S\in {\texttt {LR}}^{\la'}_{\mu'\nu'}$, there exists a unique 
$U\in SST(\nu^\pi)$ such that $(U\rightarrow H_{\mu})=H_{\la}$ and 
$Q(U\rightarrow H_{\mu})=Q(S\rightarrow H_{\mu'})$. We define $\psi(S)=U$.
Since the correspondence from $S$ to $U$ is reversible, $\psi$ is a bijection from ${\texttt {LR}}^{\la'}_{\mu'\nu'}$ to ${\texttt {LR}}^{\la}_{\mu\nu^\pi}$.

\begin{ex} \label{ex:companion}
{\rm 
Let $\lambda = (7, 6, 4, 3, 2)$, $\mu = (6, 4, 2, 2)$, and $\nu = (2, 2, 2, 1, 1)$.
Let $S\in {\texttt {LR}}^{\la'}_{\mu'\nu'}$ be given by
\begin{equation*}
S=\ytableausetup {mathmode, boxsize=1.0em} 
\begin{ytableau}
\tl{1} & \tl{3} & \tl{3} & \tl{5} & \tl{7}   \\
\tl{2} & \tl{4} & \tl{6} & \none & \none \\
\end{ytableau}\quad .
\end{equation*}
The recording tableau $Q(S \rightarrow H_{\mu'})$ is given by
\begin{equation*}
\hskip 4.2cm Q(S \rightarrow H_{\mu'}) = \hskip -5cm
\begin{split}
\ytableausetup {mathmode, boxsize=1.0em} 
\begin{ytableau}
\none &\none & \none & \none & \none & \none & \none & \tl{1}  \\
\none &\none & \none & \none & \none & \tl{2}  & \tl{3} & \none  \\
\none &\none & \none & \tl{3} & \tl{4} & \none & \none & \none    \\
\none &\none & \none & \tl{4} & \none & \none & \none & \none \\
\none &\tl{5} & \tl{5} & \none & \none & \none & \none & \none \\
\end{ytableau}
\end{split}\quad . \hskip 1cm
\end{equation*}
Then the corresponding $U=\psi(S)\in {\texttt {LR}}^{\la}_{\mu\nu^\pi}$  with 
$Q(U\rightarrow H_{\mu})= Q(S \rightarrow H_{\mu'})$ is given by
\begin{equation*}
U= 
\raisebox{3ex}{
\ytableausetup {mathmode, boxsize=1.0em} 
\begin{ytableau}
 \none & \tl{1}  \\
 \none & \tl{2}  \\
 \tl{2} & \tl{3}    \\
 \tl{3} & \tl{4} \\
\tl{5} & \tl{5} \\
\end{ytableau}}
\end{equation*}	
}
\end{ex}
\vskip 2mm

\subsection{Spinor models} \label{subsec:spinor}
Let us recall the notion of spinor model of type $D$, which is a combinatorial model for the crystal $\B(\Lambda)$ ($\Lambda \in P_+$) when $\mf g$ is of type $D$ \cite{K16} (cf. \cite{K15, K18-3}). We keep the notations used in \cite{K18-3}.

We often assume that a horizontal line $L$ on the plane is given such that any box in a tableau $T$ is either below or above $L$, and denote by $T^{\texttt{body}}$ and $T^{\texttt{tail}}$ the subtableaux of $T$ placed above and below $L$, respectively.
For example, 
\begin{equation*} \label{eq:body and tail}
\begin{split}
T=\raisebox{5ex}{
\ytableausetup {mathmode, boxsize=0.9em} 
\begin{ytableau}
 \none & \none & \none & \none \\
 \none & \none & \tl{1} & \none \\
 \none & \none & \tl{2} & \none \\
 \none & \tl{1} & \tl{3} & \none \\
 \none[\!\!\!\!\mathrel{\raisebox{-0.5ex}{$\scalebox{0.45}{\dots\dots\dots\dots}$ }}] 
& \tl{2} & \tl{4} & \none[\quad\quad \mathrel{\raisebox{-0.5ex}{$\scalebox{0.45}{\dots\dots\dots\dots}$\ ${}_{\scalebox{0.75}{$L$}}$}}] \\
 \none & \tl{3} & \none & \none \\
 \none & \tl{5} & \none & \none \\
  \none & \none & \none & \none \\
\end{ytableau}}\quad\quad\quad
T^{\texttt{body}}=\raisebox{5ex}{
\ytableausetup {mathmode, boxsize=0.9em} 
\begin{ytableau}
 \none & \none & \none \\
 \none & \tl{1} & \none \\
 \none & \tl{2} & \none \\
 \tl{1} & \tl{3} & \none \\
 \tl{2} & \tl{4} & \none \\
\end{ytableau}}\quad
T^{\texttt{tail}}=\raisebox{5ex}{
\ytableausetup {mathmode, boxsize=0.9em} 
\begin{ytableau}
 \none & \none & \none \\
 \none & \none & \none \\
 \none & \none & \none \\
 \none & \none & \none \\
 \none & \none & \none \\
 \tl{3} & \none & \none \\
 \tl{5} & \none & \none \\
\end{ytableau}}\quad\quad
\end{split}
\end{equation*}
where the dotted line denotes $L$.

For a tableau $U$ with the shape of a single column, let {$\textrm{ht}(U)$} denote the height of $U$ and we put $U(i)$ (resp. $U[i]$) to be $i$-th entry of $U$ from bottom (resp. top). 
We also write
\begin{equation*} %
	U=\left(U(\ell),\dots , U(1)\right)=\left(U[1], \dots, U[\ell] \right),
\end{equation*} 
where $\ell={\rm ht}(U)$.
To emphasize gluing and cutting tableaux which appear on the horizontal line L,
we also use the notations
\begin{equation*} %
\begin{split}
& U^{\texttt{body}}\boxplus U^{\texttt{tail}} = U,\quad U \boxminus U^{\texttt{tail}} = U^{\texttt{body}}.
\end{split}
\end{equation*}

For $a,b,c\in \Z_+$, let $\lambda(a, b, c) = (2^{b+c}, 1^a)/(1^b)$.
Let $T$ be a tableau of shape $\lambda(a, b, c)$. We denote the left and right columns of $T$ by $T^{\texttt{L}}$ and $T^{\texttt{R}}$ , respectively. 

For $T \in {SST}(\lambda(a, b, c))$ and $0 \le k \le \min\{ a, b \}$, we slide down $T^{\texttt{R}}$ by $k$ positions to have a tableau $T'$ of shape $\lambda(a-k, b-k, c+k)$. We define $\mathfrak{r}_T$ to be the maximal $k$ such that $T'$ is semistandard.

\begin{df}\label{def:jdt}
{\rm
For $T\in SST(\la(a,b,c))$ with ${\mf r}_T=0$, we define  
\begin{itemize}
\item[(1)] $\mc E T$ to be the tableau in $SST(\la(a-1,b+1,c))$ obtained from $T$ by applying jeu de taquin sliding to the position below the bottom of $T^{\texttt{R}}$, when $a>0$,

\item[(2)] $\mc F T$ to be the tableau in $SST(\la(a+1,b-1,c))$ obtained from $T$ by applying jeu de taquin sliding to the position above the top of $T^{\texttt{L}}$, when $b>0$.

\end{itemize}
Here we assume that $\mathcal{E}T = \mathbf{0}$ and $\mathcal{F}T = \mathbf{0}$ when $a=0$ and $b=0$, respectively, where ${\bf 0}$ is a formal symbol. In general, for $T \in {SST}(\lambda(a, b, c))$ with $\mathfrak{r}_T=k$, we define $\mc{E}T = \mc{E}T'$ and $\mc{F}T = \mc{F}T'$, where $T'$ is obtained from $T$ by sliding down $T^{\texttt{R}}$ by $k$ positions and hence $\mf{r}_{T'}=0$. }
\end{df}

Let
{\allowdisplaybreaks
\begin{equation*} %
\begin{split}
	\mathbf{T}(a)  & = \left\{\, T \,|\, T \in {SST}(\lambda(a,b,c)),\ b,c \in 2\Z_+,\ \mathfrak{r}_T \le 1 \,\right\}\quad (a\in\Z_+), \\
	\overline{\mathbf{T}}(0) &   = \bigsqcup_{b, c\, \in 2\mathbb{Z}_+} {SST}(\lambda(0, b, c+1)),  \quad \ \ \mathbf{T}^{\textrm{sp}}  = \bigsqcup_{a \in \mathbb{Z}_+} {SST}((1^a)), \\
	\mathbf{T}^{\textrm{sp}+} &  = \{\, T \, |\,  T \in \mathbf{T}^{\textrm{sp}}, \, \mathfrak{r}_T = 0\, \}, \ \ \mathbf{T}^{\textrm{sp}-} = \{ \,T \, |\, T \in \mathbf{T}^{\textrm{sp}}, \, \mathfrak{r}_T = 1\, \},
\end{split}
\end{equation*}
}
where we define $\mathfrak{r}_T$ of $T \in \mathbf{T}^{\textrm{sp}}$ to be the residue of $\textrm{ht}(T)$ modulo 2.

For $T \in \mathbf{T}(a)$, we define the pairs $(T^{{\texttt{L}}*}, T^{{\texttt{R}}*})$ and $({}^{\texttt{L}} T, {}^{\texttt{R}} T)$ by
\vskip 1mm
\begin{equation} \label{eq:def_pairs}
\begin{split}
	(T^{{\texttt{L}}*}, T^{{\texttt{R}}*})  & = \ ((\mc{F}T)^{\texttt{L}},(\mc{F}T)^{\texttt{R}}), \quad \textrm{when $\mathfrak{r}_T = 1$}, \\
 	({}^{\texttt{L}} T, {}^{\texttt{R}} T) &  = 
	\begin{cases}
	((\mc{E}^a T)^{\texttt{L}},(\mc{E}^a T)^{\texttt{R}}), & \text{if $\mathfrak{r}_T = 0$,} \\
	((\mc{E}^{a-1} T)^{\texttt{L}},(\mc{E}^{a-1} T)^{\texttt{R}}), & \textrm{if $\mathfrak{r}_T = 1$.}
	\end{cases}
\end{split}	
\end{equation}
\vskip 1mm

\begin{df}\label{def:admissibility}
{\rm
Let $a, a' \in \mathbb{Z}_+$ be given with $a \ge a'$. We say a pair $(T,S)$ is {\em admissible}, and write $T \prec S$ if it is one of the following cases:
\begin{itemize}
	\item[(1)] $(T,S) \in \mathbf{T}(a) \times \mathbf{T}(a')$ or $\mathbf{T}(a) \times \mathbf{T}^{\textrm{sp}}$ with 
		{\allowdisplaybreaks
		\begin{equation*}
		\begin{split}
			& (\text{i}) \ \ \ \  \textrm{ht}(T^{\texttt{R}}) \le \textrm{ht}(S^{\texttt{L}}) - a' + 2\mathfrak{r}_T \mathfrak{r}_S, \\
			& (\text{ii}) \ \ \left\{ \begin{array}{ll} 
				T^{\texttt{R}}(i) \le {}^{\texttt{L}} S(i), & \textrm{if} \ \mf{r}_T \mf{r}_S  = 0, \\
				T^{{\texttt{R}}*}(i) \le {}^{\texttt{L}} S(i), & \textrm{if} \ \mathfrak{r}_T   \mathfrak{r}_S = 1, 
			\end{array} \right. \\
			& (\text{iii}) \ \ \left\{ \begin{array}{ll} 
				{}^{\texttt{R}} T(i+a-a') \le S^{\texttt{L}}(i), & \textrm{if} \ \mathfrak{r}_T  \mathfrak{r}_S = 0, \\
				{}^{\texttt{R}} T(i+a-a'+\varepsilon) \le S^{{\texttt{L}}*}(i), & \textrm{if} \ \mathfrak{r}_T   \mathfrak{r}_S = 1, 
			\end{array} \right. \\
		\end{split}
		\end{equation*}
		}
		for $i \ge 1$. Here $\varepsilon = 1$ if $S \in \mathbf{T}^{\textrm{sp}-}$ and $0$ otherwise, and we assume that $a' = \mathfrak{r}_S$, $S = S^{\texttt{L}} = {}^{\texttt{L}} S = S^{{\texttt{L}}*}$ when $S \in \mathbf{T}^{\textrm{sp}}$.
		
	\item[(2)] $(T,S) \in \mathbf{T}(a) \times \overline{\mathbf{T}}(0)$ with $T \prec S^{\texttt{L}}$ in the sense of (1), regarding $S^{\texttt{L}} \in \mathbf{T}^{\textrm{sp}-}$.
	
	\item[(3)] $(T,S) \in \overline{\mathbf{T}}(0) \times \overline{\mathbf{T}}(0)$ or $\overline{\mathbf{T}}(0) \times \mathbf{T}^{\textrm{sp}-}$ with $(T^{\texttt{R}}, S^{\texttt{L}}) \in \overline{\mathbf{T}}(0)$.
\end{itemize}}
\end{df}
\vskip 2mm

\begin{rem}\label{rem:admissibilty for spin column}
{\rm 
\begin{itemize}
	\item[(1)] For $T\in {\bf T}(a)$, we assume that the subtableau of single column with height $a$ is below L and hence equal to $T^{\texttt{tail}}$.
\begin{equation*}
\begin{split}
T=\ \raisebox{6ex}{
\ytableausetup {mathmode, boxsize=0.9em} 
\begin{ytableau}
 \none & \none & \none & \none \\
 \none & \none & \tl{1} & \none \\
 \none & \none & \tl{2} & \none \\
 \none & \tl{1} & \tl{3} & \none \\
 \none[\ \mathrel{\!\!\!\!\raisebox{-0.5ex}{$\scalebox{0.45}{\dots\dots\dots\dots}$ }}] 
& \tl{2} & \tl{4} & \none[\quad\quad \mathrel{\raisebox{-0.5ex}{$\scalebox{0.45}{\dots\dots\dots\dots}$\ ${}_{\scalebox{0.7}{$L$}}$}}] \\
 \none & \tl{3} & \none & \none \\
 \none & \tl{5} & \none & \none \\
\end{ytableau}}\quad\quad  \in {\bf T}(2)
\end{split}
\end{equation*} 

	\item[(2)] Let $S \in \mathbf{T}^{\textrm{sp}}$ be given, and let $\varepsilon$ be the residue of ${\rm ht}(S)$ modulo 2. 
We may assume that $S= U^{\texttt{L}}$ for some $U\in {\bf T}(\varepsilon)$, where $U^{\texttt{R}}(i)$ ($i\ge 1$) are sufficiently large. 
Then we have  $S = U^{\texttt{L}} = {}^{\texttt{L}} U$. 
We assume that $S^{\texttt{tail}}=U^{\texttt{tail}}$, which is non-empty when $\varepsilon=1$.
If $\varepsilon=1$, then $U^{\texttt{L}*}$ is obtained from $S$ by adding $U^{\texttt{R}}(1)$ to the bottom of $S$ since $\mf{r}_U=1$. Note that $\mf{r}_S=\varepsilon=\mf{r}_U$.

This implies that $T\prec S$ if and only if $T\prec U$  for $T\in {\bf T}(a)$. So we may understand the admissibility conditions in Definition \ref{def:admissibility}(1) for $(T,S) \in \mathbf{T}(a) \times \mathbf{T}^{\textrm{sp}}$ as induced from the ones for $(T,U) \in \mathbf{T}(a) \times \mathbf{T}(\varepsilon)$ (cf. Example \ref{ex:spin minus}).

	\item[(3)] Let $T\in \ov{\bf T}(0)$ be given. We assume that $T^{\texttt{L}}, T^{\texttt{R}}\in {\bf T}^{\rm sp-}$ so that $T^{\texttt{tail}}$ is non-empty. 
This means that $\left( T^{\texttt{L}} \right)^{\texttt{tail}}$ and $\left( T^{\texttt{R}} \right)^{\texttt{tail}}$ are non-empty in the sense of (2). 
\end{itemize}
}
\end{rem}
\vskip 2mm

From now on we assume that $\mf g$ is the Kac-Moody Lie algebras of type $D_\infty$, whose Dynkin diagram, set of simple roots $\Pi=\{\,\alpha_i\,|\,i\in I\,\}$, and fundamental weight $\Lambda_i$ $(i\in I)$ are
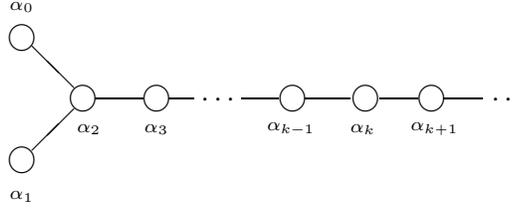
\begin{center}
\setlength{\unitlength}{0.16in} \hskip -3cm
\hskip 2cm \begin{picture}(24,5.8)
\put(6,0){\makebox(0,0)[c]{$\bigcirc$}}
\put(6,4){\makebox(0,0)[c]{$\bigcirc$}}
\put(8,2){\makebox(0,0)[c]{$\bigcirc$}}
\put(10.4,2){\makebox(0,0)[c]{$\bigcirc$}}
\put(14.85,2){\makebox(0,0)[c]{$\bigcirc$}}
\put(17.25,2){\makebox(0,0)[c]{$\bigcirc$}}
\put(19.4,2){\makebox(0,0)[c]{$\bigcirc$}}
\put(6.35,0.3){\line(1,1){1.35}} \put(6.35,3.7){\line(1,-1){1.35}}
\put(8.4,2){\line(1,0){1.55}} \put(10.82,2){\line(1,0){0.8}}
\put(13.2,2){\line(1,0){1.2}} \put(15.28,2){\line(1,0){1.45}}
\put(17.7,2){\line(1,0){1.25}} \put(19.8,2){\line(1,0){1.25}}
\put(12.5,1.95){\makebox(0,0)[c]{$\cdots$}}
\put(22,1.95){\makebox(0,0)[c]{$\cdots$}}
\put(6,5){\makebox(0,0)[c]{\tiny $\alpha_{0}$}}
\put(6,-1.2){\makebox(0,0)[c]{\tiny $\alpha_{1}$}}
\put(8.2,1){\makebox(0,0)[c]{\tiny $\alpha_{2}$}}
\put(10.4,1){\makebox(0,0)[c]{\tiny $\alpha_{3}$}}
\put(14.8,1){\makebox(0,0)[c]{\tiny $\alpha_{k-1}$}}
\put(17.15,1){\makebox(0,0)[c]{\tiny $\alpha_k$}}
\put(19.5,1){\makebox(0,0)[c]{\tiny $\alpha_{k+1}$}}
\end{picture}

\begin{equation*}
\begin{split}
&\Pi=\{\, \alpha_0=-\epsilon_1-\epsilon_2, \  \alpha_i=\epsilon_i-\epsilon_{i+1} \ (i\geq 1)\, \},\\
&\La_i =
\begin{cases}
\La_0+\epsilon_1, & \text{if $i=1$}, \\
2\La_0+\epsilon_1+\cdots+\epsilon_i, & \text{if $i>1$}.
\end{cases}
\end{split}
\end{equation*}
\end{center}
Here we assume that the index set for simple roots is $I=\Z_+$, and the weight lattice is $P=\Z\Lambda_0\oplus \left( \bigoplus_{i\geq 1}\Z\epsilon_i \right)$.
Let $\mathfrak{l}$ be the subalgebra of ${\mf g}$ associated to $\Pi\setminus\{\alpha_0\}$, which is of type $A_{+\infty}$.

Let $\B$ be one of $\mathbf{T}(a)$ $(a \in \mathbb{Z}_+)$, $\mathbf{T}^{\textrm{sp}}$, and $\overline{\mathbf{T}}(0)$. Let us describe the $\mf g$-crystal structure on $\B$.
Let $T \in \B$ given. 
Recall that $SST(\la)$ ($\la\in\cP$) has an $\mf l$-crystal structure \cite{KN}.
So we may regard $\B$ as a subcrystal of an $\mathfrak{l}$-crystal $\bigsqcup_{\lambda \in \cP} {SST}(\lambda)$ and hence define $\te_i T$ and $\tf_i T$ for $i \in {I} \setminus \{ 0 \}$.
Let ${\rm wt}_{\mf l}(T)=\sum_{i\geq 1}m_i\epsilon_i$ be the $\mf l$-weight of $T$, where $m_i$ is the number of occurrences of $i$ in $T$.
Next, we define $\te_0 T$ and $\tf_0 T$ as follows:

\begin{itemize}
	\item[(1)] When $\B = \mathbf{T}^{\textrm{sp}}$, we define $\te_0$ to be the tableau obtained from $T$ by removing a domino 
\raisebox{-.6ex}{{\tiny ${\def\lr#1{\multicolumn{1}{|@{\hspace{.6ex}}c@{\hspace{.6ex}}|}{\raisebox{-.3ex}{$#1$}}}\raisebox{-.6ex}
{$\begin{array}[b]{c}
\cline{1-1}
\lr{ 1 }\\
\cline{1-1}
\lr{ \!2\!}\\
\cline{1-1}
\end{array}$}}$}}
if $T$ has
\raisebox{-.6ex}{{\tiny ${\def\lr#1{\multicolumn{1}{|@{\hspace{.6ex}}c@{\hspace{.6ex}}|}{\raisebox{-.3ex}{$#1$}}}\raisebox{-.6ex}
{$\begin{array}[b]{c}
\cline{1-1}
\lr{ 1}\\
\cline{1-1}
\lr{ \!2\!}\\
\cline{1-1}
\end{array}$}}$}}
on its top, and \textbf{0} otherwise. We define $\tf_0 T$ in a similar way by adding
\raisebox{-.6ex}{{\tiny ${\def\lr#1{\multicolumn{1}{|@{\hspace{.6ex}}c@{\hspace{.6ex}}|}{\raisebox{-.3ex}{$#1$}}}\raisebox{-.6ex}
{$\begin{array}[b]{c}
\cline{1-1}
\lr{ 1 }\\
\cline{1-1}
\lr{ \!2\!}\\
\cline{1-1}
\end{array}$}}$}}.

	\item[(2)] When $\B = \mathbf{T}(a)$ or $\overline{\mathbf{T}}(0)$, we define $\te_0 T=\te_0 \left( T^{\texttt{R}} \otimes T^{\texttt{L}} \right)$ regarding $\B \subset \left( \mathbf{T}^{\textrm{sp}} \right)^{\otimes 2}$ by tentor product rule \eqref{eq:tensor_product_rule}. We define $\tf_0 T$ similarly.
\end{itemize}

Put
\begin{equation*}
\begin{split}
\textrm{wt}(T) &=
\begin{cases} 
2\La_0 + \textrm{wt}_{\mathfrak{l}}(T), & \textrm{if $T \in \mathbf{T}(a)$ or $\ov{\bf T}(0)$}, \\
\La_0 + \textrm{wt}_{\mathfrak{l}}(T), & \textrm{if $T\in {\bf T}^{\rm sp}$}.
\end{cases},\\
\varepsilon_i (T) & = \max \{\, k \, | \, \te_i^k T \neq \textbf{0}\, \} \quad 
\varphi_i (T) = \max \{\, k \, | \, \tf_i^k (T) \neq \textbf{0}\, \}.
\end{split}
\end{equation*}
Then $\B$ is a $\mathfrak{g}$-crystal with respect to $\te_i$ and $\tf_i$, $\varepsilon_i$, and $\varphi_i$ for $i \in I$.
By \cite[Proposition 4.2]{K16}, we have
\begin{equation*}
\begin{split}
	& \mathbf{T}(a) \cong \B(\La_{a}) \ \  (a \ge 2), \\ 
	& \mathbf{T}(0) \cong \B(2\La_{0}), \quad \overline{\mathbf{T}}(0) \cong \B(2\La_{1}), \quad \mathbf{T}(1) \cong \B(\La_{0}+\La_1), \\
	& \mathbf{T}^{\textrm{sp}-} \cong \B(\La_1), \quad \mathbf{T}^{\textrm{sp}+} \cong \B(\La_0).
\end{split}
\end{equation*}

For $n\geq 1$, let 
\begin{equation*}
{\mc P}({\rm O}_n)=\{\,\mu=(\mu_1,\cdots,\mu_n)\,|\,\mu_i\in\Z_+, \ \mu_1\geq\ldots \geq \mu_{n},\
\mu'_1+\mu'_2\leq n\,\},
\end{equation*}
where $\mu'=(\mu'_1, \mu'_2, \cdots )$ is the conjugate partition of $\mu$. 
Recall that $\mc{P}({\rm O}_n)$ parameterizes the complex finite-dimensional representations of the orthogonal group ${\rm O}_n$. 

We may also use $\mc{P}({\rm O}_n)$ to parametrize $P_+$ for $\mf g$.
More precisely, for $\mu\in \mc{P}({\rm O}_n)$, if we put 
\begin{equation*}
	\Lambda(\mu) = n\La_0+\mu'_1\epsilon_1+\mu'_2\epsilon_2 + \cdots,
\end{equation*}
then we have $P_+=\{\,\La(\mu)\,|\, \mu\in \bigsqcup_n \mc{P}({\rm O}_n)\,\}$ the set of dominant integral weights for $\mf g$. 

For $\mu \in \mathcal{P}({\rm O}_n)$, let $q_\pm$ and $r_\pm$ be non-negative integers such that 
\begin{equation*} 
\begin{cases}	
n-2\mu'_1 = 2q_+ + r_+, & \textrm{if} \ n-2\mu'_1 \ge 0, \\ 
2\mu'_1 - n = 2q_- + r_-, & \textrm{if} \ n-2\mu'_1 < 0, 
\end{cases}	
\end{equation*} 
where $r_\pm = 0, 1$. 
Let $\ov{\mu} = (\overline{\mu}_i) \in \cP$ be such that 
$\ov{\mu}'_1 = n-\mu'_1$ and $\ov{\mu}'_i = \mu'_i$ for $i \ge 2$ and let $M_+ = \mu'_1$ and $M_- = \ov{\mu}'_1$. 
Put
\begin{equation}\label{eq:hat{T}}
\widehat{\mathbf{T}}(\mu, n) = 
\begin{cases}
\mathbf{T}(\mu_1) \times \cdots \times \mathbf{T}(\mu_{M_+}) \times \mathbf{T}(0)^{\times q_+} \times (\mathbf{T}^{\textrm{sp+}})^{\times r_+}, & \textrm{if} \ n - 2\mu'_1 \ge 0, \\
\mathbf{T}(\ov{\mu}_1) \times \cdots \times \mathbf{T}(\ov{\mu}_{M_-}) \times \overline{\mathbf{T}}(0)^{\times q_-} \times (\mathbf{T}^{\textrm{sp}-})^{\times r_-}, & \textrm{if} \ n - 2\mu'_1 < 0.
\end{cases}
\end{equation}
Define
\begin{equation*} 
	\mathbf{T}(\mu, n) = 
	\{\, \mathbf{T} = (\dots,T_2,T_1) \in \widehat{\mathbf{T}}(\mu, n) \, | \, T_{i+1} \prec T_{i}   \ (i \ge 1) \,\}.
\end{equation*}

When considering $\mathbf{T} = (\dots,T_2,T_1) \in \widehat{\mathbf{T}}(\mu, n)$,
it is often helpful to imagine that $T_1, T_2, \dots$ are arranged from right to left on a plane, where the horizontal line $L$ separates $T_i^{\texttt{body}}$ and $T_i^{\texttt{tail}}$ simultaneously.

\begin{ex}
{\em
Let $n=8$ and $\mu = (4, 3, 3, 2) \in \mathcal{P}({\rm O}_8)$ and let ${\bf T} = (T_4,T_3, T_2, T_1)$ given by
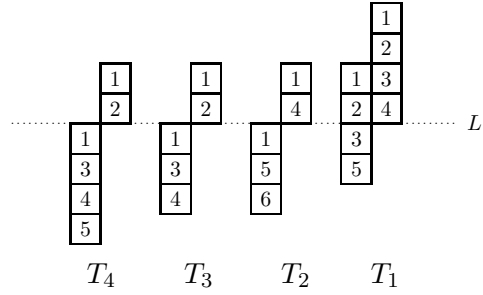
\begin{equation*}
\begin{split}
\ytableausetup {mathmode, boxsize=1.0em} 
& \begin{ytableau}
\none &\none & \none & \none & \none & \none & \none & \none & \none & \none & \none & \tl{1} & \none  \\
\none &\none & \none & \none & \none & \none & \none & \none & \none & \none & \none & \tl{2} & \none  \\
\none &\none & \tl{1} & \none & \none & \tl{1} & \none & \none & \tl{1} & \none & \tl{1} & \tl{3} & \none \\
\none[\!\!\!\!\mathrel{\raisebox{-0.7ex}{$\scalebox{0.45}{\dots\dots\dots\dots}$}}] &\none & \tl{2} & \none[\mathrel{\raisebox{-0.7ex}{$\scalebox{0.45}{\dots\dots}$}}] & \none & \tl{2} & \none[\mathrel{\raisebox{-0.7ex}{$\scalebox{0.45}{\dots\dots}$}}] & \none & \tl{4} & \none[\mathrel{\raisebox{-0.7ex}{$\scalebox{0.45}{\dots\dots}$}}]& \tl{2} & \tl{4} & \none[\ \mathrel{\raisebox{-0.7ex}{\quad $\scalebox{0.45}{\dots\dots\dots\dots}$\ ${}_{\scalebox{0.75}{$L$}}$}}] \\
\none &\tl{1} & \none & \none & \tl{1} & \none & \none & \tl{1} & \none & \none & \tl{3} & \none & \none  \\
\none &\tl{3} & \none & \none & \tl{3} & \none & \none & \tl{5} & \none & \none & \tl{5} & \none & \none  \\
\none &\tl{4} & \none & \none & \tl{4} & \none & \none & \tl{6} & \none & \none & \none & \none & \none  \\
\none &\tl{5} & \none & \none & \none & \none & \none & \none & \none & \none & \none & \none & \none  \\
\end{ytableau} \quad\quad \\
&\ \hskip 5mm T_4 \hskip 9mm T_3 \hskip 9mm T_2 \hskip 8mm T_1 \hskip 13mm 
\end{split}
\end{equation*} 
\vskip 3mm
Then we can check that $T_4 \prec T_3 \prec T_2 \prec T_1$. Thus, ${\bf T} \in {\bf T}(\mu, 8)$. 
where the dotted line denotes the common horizontal line L.
}
\end{ex}

We regard $\widehat{\bf T}(\mu,n)$ as a $\mf g$-crystal by identifying ${\bf T}=(\dots,T_2,T_1)\in \widehat{\bf T}(\mu,n)$ with $T_1\otimes T_2\otimes \dots$, and regard ${\bf T}(\mu,n)$ as its subcrystal. Then we have the following.

\begin{thm} \label{thm:fundamental_theorem} \cite[Theorem 4.3--4.4]{K16} 
For $\mu \in \mathcal{P}({\rm O}_n)$,
$\mathbf{T}(\mu, n)$ is a connected crystal with highest weight $\La(\mu)$. Furthermore,
we have $$\mathbf{T}(\mu, n) \cong \mathbf{B}(\Lambda(\mu)).$$  
\end{thm}
\qed

We call $\mathbf{T}(\mu, n)$ the {\em spinor model for $\mathbf{B}(\Lambda(\mu))$} (in type $D_\infty$).
From now on, we fix $\mu\in \mc{P}({\rm O}_n)$ throughout the paper.


\section{Separation algorithm} \label{sec:separation}

We introduce a combinatorial algorithm on ${\bf T}(\mu, n)$, so-called {\em separation}, which plays a crucial role in this paper.
Let us briefly explain its motivation and the result in this section.

First, we observe that the multiplicity \eqref{eq:branching mult} is equal to the number of $\mf{l}$-highest weight elements ${\bf T} \in {\bf T}(\mu, n)$ with highest weight $\la'$, that is, $\widetilde{e}_i {\bf T} = {\bf 0} \, (i \neq 0)$ and ${\bf T} \equiv_{\mf l} H_{\lambda'}$ from the see-saw pairs in the Howe duality of $({\rm O}_n, D_{\infty})$ \cite{Wa}.
To have a combinatorial description of ${\bf T}$, we then introduce an algorithm which gives an injective map 
\begin{equation}\label{eq:separation}
\xymatrixcolsep{3pc}\xymatrixrowsep{0pc}\xymatrix{
{\bf T}  \ \ar@{|->}[r] & \ (H_{(\delta')^{\pi}}\,,\, U)},
\end{equation}
where $\delta \in \cP^{(2)}_n$ and $U \in \texttt{LR}^{\lambda'}_{\delta' \mu'}$.
So this reduces our problem to characterizing the image of $\mf{l}$-highest weight elements under the map \eqref{eq:separation}. 

The idea of considering such a map, which we call {\em separation}, is basically from \cite{K18-3}, where the second author explains the stable branching rules \cite{HTW} in terms of crystals. 
The separation algorithm for types $B$ and $C$ is already present in \cite{K18-2}.
Roughly speaking, it is given by sliding horizontally the tails of ${\bf T}$ (using the jeu de taquin sliding) to the leftmost one as far as possible so that the resulting tableau gives $U$ and the remaining one in the body gives $H_{(\delta')^\pi}$.  
Also, the characterization of the image of \eqref{eq:separation} is given in \cite{LL} when ${\bf T}(\mu,n)$ is of type $C$.

However, the separation algorithm for types $B$ and $C$ does not work well on the spinor model of type $D$ due to more involved conditions for admissibility in ${\bf T}(\mu,n)$.

To overcome this difficulty, we introduce an operator {\em sliding} \eqref{eq:definition_S} which is given by a non-trivial sequence of jeu de taquin slidings, and also moves a tail in ${\bf T}$ by one position to the left horizontally (see Example \ref{ex:illustration of S}).

A key property is that our sliding is compatible with the type $A$ crystal structure on ${\bf T}(\mu,n)$ so that we obtain another $\mf{l}$-highest weight element $\td{\bf T}$ and ${\bf T}=\td{\bf T}\otimes U$ as an element in a crystal of type $A$, where $U$ is the leftmost column in ${\bf T}$ (see Lemmas \ref{lem:sliding one step} and \ref{lem:sliding one step-2}). 
Hence this enables us to define the map \eqref{eq:separation} by applying the sliding successively.
We refer the reader to Examples \ref{ex:3.13} and \ref{ex:separation for negative case} for its illustration.

\subsection{$\mathfrak{l}$-highest weight elements} \label{subsec:l-highest weight element}

\begin{df} \label{df:l-highest weight elements}
{\em 
Let
\begin{equation*}%
{\bf H}(\mu,n)=\{\,{\bf T}\,\,|\,\,{\bf T}\in \mathbf{T}(\mu, n),\ \te_i{\bf T}={\bf 0}\ (i\neq 0)\,\},
\end{equation*}
and call ${\bf T}\in {\bf H}(\mu,n)$ an {\em $\mf l$-highest weight element}.
}
\end{df}

The goal of this subsection is to give some necessary conditions for ${\bf T}\in {\bf T}(\mu,n)$ to be in ${\bf H}(\mu,n)$, 
which will be used when we define the separation algorithm.

Note that for ${\bf T}\in {\bf T}(\mu,n)$, we have ${\bf T}\in {\bf H}(\mu,n)$ if and only if ${\bf T}\equiv_{\mf l} H_\la$ for some $\la\in \cP$. Hence ${\bf H}(\mu,n)$ parametrizes the connected $\mf l$-crystals in ${\bf T}(\mu,n)$.

In this subsection, we assume that $n-2\mu'_1\geq 0$ (The case when $n-2\mu'_1 < 0$ will be considered in subsection \ref{subsec:negative case}).
Suppose that $n=2l+r$, where $l\geq 1$ and $r=0,1$.
Let ${\bf T} \in {\bf T}(\mu, n)$ be given and write 
\begin{equation}\label{eq:T notation}
{\bf T} = (T_l, \dots, T_1,T_0),
\end{equation}
where
$T_i\in {\bf T}(a_i)$ for some $a_i\in\Z_+$ ($1\leq i\leq l$), and $T_0\in {\bf T}^{\rm sp+}$ (resp. $T_0=\emptyset$) when $r=1$ (resp. $r=0$).
Let ${\rm sh}(T_i) = \lambda(a_i, b_i, c_i)$ and $\mathfrak{r}_{T_i} = \mathfrak{r}_i$ for $1\leq i\leq l$.

The lemma below follows directly from the tensor product rule \eqref{eq:tensor_product_rule}.

\begin{lem}\label{lem:criterion_highest_weight_elt}
Put $U_0=T_0$, $U_{2k-1} = T_k^{\texttt{R}}$ and $U_{2k} = T_k^{\texttt{L}}$ for $1\leq k\leq l$. 
Then ${\bf T}$ is an $\mathfrak{l}$-highest weight element if and only if $(U_i, \dots, U_0)$ is a $\mathfrak{l}$-highest weight element for $i\geq 0$, where we understand $(U_i, \dots, U_0)=U_0\otimes \dots \otimes U_i$ as an element of an $\mf l$-crystal. \qed
\end{lem}

\begin{df} \label{def:pseudo_H}
{\rm Let ${{\bf{H}}}^\circ(\mu, n)$ be a subset of ${\bf T}(\mu, n)$ consisting of ${\bf T} = (T_i)$ satisfying the following conditions: for each $i\geq 1$,
\begin{itemize}
 	\item[({\rm H0})] $T_0[k]=k$ for $k\geq 1$,
	
	\item[({\rm H1})] $T_i^{\texttt{L}}$ and $T_i^{\texttt{R}}$ are of the form 
		\begin{equation*}
		\begin{split}
			& T^{\texttt{R}}_i = \left(1, 2,  \dots ,  b_i + c_i - 1 , T_i^{\texttt{R}}(1)\right)\boxplus  \emptyset ,\\
			& T^{\texttt{L}}_i = \left(1, 2,  \dots ,  c_i - 1 ,  c_i   \right)\boxplus \left( T_i^{\texttt{L}}(a_i) ,   \dots ,T_i^{\texttt{L}}(1) \right),
		\end{split}
		\end{equation*}
		
	\item[({\rm H2})] the entries $T_i^{\texttt{R}}(1)$ and $T_i^{\texttt{L}}(a_i)$ satisfy 
		\begin{equation*}
		\begin{split}
			{\rm (i)} \ \ & \textrm{if} \ \mathfrak{r}_i = 0, \textrm{then} \ T^R_i(1) = b_i+c_i, \\
			{\rm (ii)} \ \ & \textrm{if $\mathfrak{r}_i = 1$, then}
				\left\{\begin{array}{l} T^{\texttt{R}}_i(1) =  b_i+c_i \ \textrm{or} \ T^{\texttt{R}}_i(1) \ge c_{i-1}+1+\mathfrak{r}_{i-1}, \\ T_i^{\texttt{L}}(a_i) = c_i + 1. \end{array}
				\right.
		\end{split}
		\end{equation*} 
\end{itemize} 
Here we assume that $c_0=\infty$ and ${\mf r}_0=0$.
}
\end{df}

\begin{lem} \label{cor:highest_weight_element}
For ${\bf T}\in {\bf H}^\circ(\mu,n)$, we have either $T^{\texttt{R}}_{i+1}(1) < T^{\texttt{L}}_i(a_i)$ or $T^{\texttt{R}}_{i+1}(1) > T^{\texttt{L}}_i(a_i)$ for each $i$. 
Furthermore, $T^{\texttt{R}}_{i+1}(1) > T^{\texttt{L}}_i(a_i)$ implies ${\mf r}_i{\mf r}_{i+1}=1$, and 
${\mf r}_i{\mf r}_{i+1}=0$ implies  $T^{\texttt{R}}_{i+1}(1) < T^{\texttt{L}}_i(a_i)$.
\end{lem}
\pf  Let $i=1,\dots,l-1$ be given. If $\mathfrak{r}_{i}\mathfrak{r}_{i+1} = 0$, then by Definition \ref{def:admissibility}(1)-(i) and (H1), we have $T^{\texttt{R}}_{i+1}(1) = b_{i+1}+c_{i+1} \le c_i < T^{\texttt{L}}_i(a_i)$.

Suppose that $\mathfrak{r}_{i}\mathfrak{r}_{i+1} = 1$. 
By (H2)(ii), we have $T^{\texttt{L}}_i(a_i)=c_i+1$, and $T^{\texttt{R}}_{i+1}(1)=b_{i+1}+c_{i+1}$ or $\geq c_i+2$.
If $T^{\texttt{R}}_{i+1}(1)\geq c_i+2$, then it is clear that $T^{\texttt{R}}_{i+1}(1) > T^{\texttt{L}}_i(a_i)$.
So we assume that $T^{\texttt{R}}_{i+1}(1) = b_{i+1}+c_{i+1}$. 
Note that $T^{\texttt{R}}_{i+1}(1) = b_{i+1}+c_{i+1} \le c_i+2$ by Definition \ref{def:admissibility}(1)-(i). 
If $b_{i+1} + c_{i+1} = c_i +2$, then $T^{\texttt{R}}_{i+1}(1) > T^{\texttt{L}}_i(a_i)$. 
If $b_{i+1} + c_{i+1} < c_i +2$, then $b_{i+1} + c_{i+1} \le c_i$ since both $b_{i+1}+c_{i+1}$ and $c_i$ are even. 
So $T^{\texttt{R}}_{i+1}(1) < T^{\texttt{L}}_i(a_i)$. 

Finally, suppose that $T^{\texttt{R}}_{i+1}(1) > T^{\texttt{L}}_i(a_i)$. 
If $\mathfrak{r}_{i}\mathfrak{r}_{i+1} =0$, then by Definition \ref{def:admissibility}(1)-(ii) we have $T^{\texttt{R}}_{i+1}(1) \leq T^{\texttt{L}}_i(a_i+1)<T^{\texttt{L}}_i(a_i)$ which is a contradiction. This proves the lemma.
\qed
\vskip 1mm

Now we verify that the ${\mf l}$-highest weight elements satisfy Definition \ref{def:pseudo_H}(H0)--(H2).
In particular, the admissibility in Definition \ref{def:admissibility} implies the condition (H2).

\begin{prop} \label{prop:highest_weight_vectors}
We have
${\bf{H}}(\mu, n) \subset {\bf{H}}^\circ(\mu, n).$
\end{prop}
\pf  
Suppose that ${\bf T} \in {\bf H}(\mu, n)$. 
By Lemma \ref{lem:criterion_highest_weight_elt}, it is clear that $T_0$ satisfies (H0).
Let $w(T_0)w(T_1)\dots w(T_l)=w_1w_2\dots w_\ell$, and 
\begin{equation} \label{prop:eq1}
	P_l=\left( w_\ell \rightarrow (\dots \rightarrow (w_3 \rightarrow (w_2 \rightarrow w_1))) \right),
\end{equation} 
By definition of ${\bf H}(\mu,n)$, we have $P_l=H_\nu$ for some $\nu\in \cP$. 
Let $h_l$ be the height of the rightmost column of $\nu$.

Let us use induction on $l$ to show that ${\bf T} \in {\bf H}^\circ(\mu, n)$.
We also claim that 
\begin{equation}\label{prop:eq0}
h_l=c_l + {\mf r}_l.
\end{equation}
Suppose that $l=1$. 
Since ${\bf T}$ is an $\mathfrak{l}$-highest weight element and ${\bf T} \equiv_{\mathfrak{l}} {\bf T}^{\texttt{R}} \otimes {\bf T}^{\texttt{L}}$ by Lemma \ref{lem:criterion_highest_weight_elt}, it is straightforward to check that ${\bf T}$ satisfies (H1) and (H2). It is clear that $c_1=h_1+{\mf r}_1$.

Suppose that $l > 1$. Let ${\bf T}_{l-1} = (T_{l-1}, \dots, T_1,T_0)$ and let $P_{l-1}$ be the tableau in \eqref{prop:eq1} corresponding to ${\bf T}_{l-1}$. 
By induction hypothesis, ${\bf T}_{l-1}$ satisfies (H1) and (H2). Put $${P}^{\flat}_{l} = (w(T^{\texttt{R}}_l) \rightarrow P_{l-1}).$$ 
Then ${P}^\flat_{l}=H_\eta$ for some $\eta\in\cP$ by Lemma \ref{lem:criterion_highest_weight_elt}, and  
${P}_{l} = (w(T^{\texttt{L}}_l) \rightarrow {P}^\flat_{l})$. 

{\em Case 1}. Suppose that $\mathfrak{r}_l = 0$. 
Note that by Definition \ref{def:admissibility}(1)-(i), we have $b_l + c_l \le c_{l-1}$. 
Also, by Definition \ref{def:admissibility}(1)-(ii), $T^{\texttt{R}}_l(i) \le {}^{\texttt{L}} T_{l-1}(i)$ for $1 \le i \le b_l+c_l$. 
By (H1) on ${\bf T}_{l-1}$, we have ${}^{\texttt{L}} T_{l-1}[k] = k$ for $1\le k\le \ov{c}_{l-1}$, where $\ov{c}_{l-1}=c_{l-1}+\mathfrak{r}_{l-1}$.  
Hence
\begin{equation} \label{prop:eq2}
	T^{\texttt{R}}_l(i) \le c_{l-1}-i+\mathfrak{r}_{l-1}+1,
\end{equation} 
for $1 \le i \le b_l+c_l$.  
Then \eqref{prop:eq2} implies that each letter of $w(T^{\texttt{R}}_l)$ is inserted to create a box to the right of the leftmost column of $P_{l-1}$ when we consider the insertion $(w(T^{\texttt{R}}_l) \rightarrow P_{l-1})$. 
Since ${P}^\flat_{l}=H_\eta$, we have $T^{\texttt{R}}_l[k]=k$ for $1\le k \le b_l+c_l$.

By semistandardness of $T^{\texttt{body}}_l$, we have
\begin{equation*} \label{prop:eq3}
	(T_l^{\texttt{L}})^{\texttt{body}}(i) \le T_l^{\texttt{R}}(i),
\end{equation*} 
for $i\geq 1$.
This implies that each letter of $w((T^{\texttt{L}}_l)^{\texttt{body}})$ is inserted to create a box to the right of the leftmost column of $P^\flat_{l}$ when we consider the insertion $w((T^{\texttt{L}}_l)^{\texttt{body}}) \rightarrow P^\flat_{l})$, and 
$T^L_l[k] = k$ for $1 \le k \le c_l$.
Hence ${\bf T}$ satisfies (H1), (H2), and \eqref{prop:eq0}.

{\em  Case 2.} Suppose that $\mathfrak{r}_l = 1$.
When $\mathfrak{r}_{l-1} = 0$, we see that $T_l$ satisfies the conditions (H1) and (H2) by the same argument in the previous case. In particular, \eqref{prop:eq2} implies that $T^{\texttt{R}}_l(1) = b_l + c_l$. Since $T^{\texttt{L}}_l(a_l) \le T^{\texttt{R}}_l(1)$ and $P^\flat_{l}=H_\eta$, we also have $T^{\texttt{L}}_l(a_l) = c_l + 1$ and \eqref{prop:eq0}. 

Now assume that $\mathfrak{r}_{l-1} = 1$. 
When $\mathfrak{r}_{l}\mathfrak{r}_{l-1} = 1$, we need to consider the $*$-pair $(T_l^{{\texttt{L}}*}, T_l^{{\texttt{R}}*})$ of $T_l$ in  \eqref{eq:def_pairs} (recall Definition \ref{def:admissibility}). Then,
by Definition \ref{def:admissibility}(1)-(ii) and the condition (H1) on ${\bf T}_{l-1}$, we have
\begin{equation*} %
	T^{{\texttt{R}}*}_l(i) \le {}^{\texttt{L}}T_{l-1}(i) = c_{l-1}-i+2.
\end{equation*} 
We claim that $T^{\texttt{R}}_l[k] = k$ for $ 1\le k \le b_l+c_l-1$.
Let $k$ be such that $T^{{\texttt{R}}*}_l(i) = T^{\texttt{R}}_l(i)$ for $1\leq i\leq k-1$, and $T^{{\texttt{R}}*}_l(i) = T^{\texttt{R}}_l(i+1)$ for $i\ge k$.  
Since $(\mc{F}T_l, T_{l-1}, \dots, T_1) \equiv_{\mathfrak{l}} (T_l, T_{l-1}, \dots, T_1)$, which is an $\mathfrak{l}$-highest weight element, 
we see from Lemma \ref{lem:criterion_highest_weight_elt} that
each letter of $w(T^{\texttt{R}*}_l)$ is inserted to create a box to the right of the leftmost column of $P_{l-1}$ when we consider the insertion $(w(T^{\texttt{R}*}_l) \rightarrow P_{l-1})$, and $T^{R*}_l[i] = i$ for $1\le i \le b_l+c_l-1$.
This implies that $T^R_l(k)$ is between $m$ and $m+1$ for some $m \in \mathbb{Z}_+$, and hence $k=b_l+c_l$ since $(w(T^{\texttt{R}*}_l) \rightarrow P_{l-1})$ is an ${\mf l}$-highest weight element. This proves the claim, and $T^{\texttt{R}}_l$ satisfies (H1).
Furthermore, the claim implies that $T^{\texttt{R}}_l(1)$ satisfies (H2)(ii) because $P^\flat_l$ is an $\mathfrak{l}$-highest weight element. 

We consider $T^{\texttt{L}}_l$. By the same argument as in {\em Case 1}, we have $T^{\texttt{L}}_l[j] = j$ for $1\le j \le c_l$, and $k=1$ (in the previous argument) implies that $T^{\texttt{L}}_l(a_l) \le b_l + c_l - 1$. 
Therefore, we have $T^{\texttt{L}}_l(a_l) = c_l+1$ since the tableau obtained by $(T^{\texttt{L}}_l(a_l) \rightarrow (w(T_l^{\texttt{body}}) \rightarrow P^\flat_{l}))$ is an $\mathfrak{l}$-highest weight element. Finally, we can check easily that \eqref{prop:eq0} holds.
\qed

\begin{ex}{\rm
Let $\mathbf{T} = (T_2, T_1)\in {\bf T}(2,2)$ with $\mathfrak{r}_{1} = \mathfrak{r}_{2} = 1$ given by
\begin{equation*}
\begin{split}
\ytableausetup {mathmode, boxsize=0.9em} 
\begin{ytableau}
\none & \none & \none & \none & \none & \tl{1} \\
\none & \none & \none & \none & \none & \tl{2} \\
\none & \tl{1} & \none & \none & \tl{1} & \tl{3} \\
\none & \tl{4} & \none & \none & \tl{2} & \tl{4} \\
\tl{1} & \none & \none & \none & \tl{3} & \none \\
\tl{5} & \none & \none & \none & \tl{5} & \none \\
\none & \none & \none & \none & \none & \none \\
\none[\tl{$T_2^{\texttt{L}}$} \ ] & \none[\ \tl{$T_2^{\texttt{R}}$}] & \none & \none & \none[\tl{$T_1^{\texttt{L}}$} \ ] & \none[\ \tl{$T_1^{\texttt{R}}$}]
\end{ytableau}
\end{split}
\end{equation*} 
We have $w(T_1)w(T_2)=(12341235)(1415)$ and the corresponding tableau \eqref{prop:eq1} is
\begin{equation*}
\begin{split}
\ytableausetup {mathmode, boxsize=0.9em} 
\begin{ytableau}
\none & \none & \none & \none & \none & \none\\
\tl{1} & \tl{1} & \tl{1} & \tl{1} & \none & \none\\
\tl{2} & \tl{2} & \none & \none & \none & \none\\
\tl{3} & \tl{3} & \none & \none & \none & \none\\
\tl{4} & \tl{4} & \none & \none & \none & \none\\
\tl{5} & \tl{5} & \none & \none & \none & \none\\
\none & \none & \none & \none & \none & \none\\
\end{ytableau}
\end{split}
\end{equation*} 
Thus ${\bf T}$ is an $\mathfrak{l}$-highest weight element.}
\end{ex}

\subsection{Sliding algorithm}\label{subsec:sliding}
In this subsection, we introduce a combinatorial algorithm which moves a tail in ${\bf T}\in {\bf H}^\circ(\mu,n)$ by one position to the left horizontally preserving the ${\mf l}$-crystal equivalence or Knuth equivalence. 
%
Then we obtain $\widetilde{\bf T} \in {\bf H}(\widetilde{\mu}, n-1)$ from ${\bf T}\in {\bf H}(\widetilde{\mu}, n-1)$ by applying the slidings repeatedly, where $\td{\mu}=(\mu_2,\mu_3,\dots)$. 
This is a key observation to define the separation inductively in the next subsection. 
\vskip 2mm

We need the following set given by
\begin{equation*}
	\mathbf{E}^n := \underset{{(u_n, \dots, u_1) \in \mathbb{Z}_+^n}}{\bigsqcup} 
	{SST}(1^{u_n}) \times \dots \times {SST}(1^{u_1}).
\end{equation*} 
Let $(U_n, \dots, U_1) \in \mathbf{E}^n$ given. 
For $1 \le j \le n-1$ and $\mathcal{X} = \mathcal{E}, \mathcal{F}$, we define
\begin{equation} \label{eq:E_j and F_j}
	\mathcal{X}_j (U_n, \dots, U_1) = 
	\begin{cases}
	(U_r, \dots, \mathcal{X}(U_{j+1}, U_j), \dots, U_1), 
	& \textrm{if $\mc{X}(U_{j+1}, U_j) \neq \mathbf{0}$,} \\
	\mathbf{0}, &  \ \textrm{if $\mc{X}(U_{j+1}, U_j) = \mathbf{0}$}.
	\end{cases}
\end{equation} 
where $\mathcal{X}(U_{j+1}, U_j)$ is understood to be $((\mc{X}U)^{\texttt{L}},(\mc{X}U)^{\texttt{R}})$ for $U\in SST(\la(a,b,c))$ with ${\mf r}_U=0$ and $(U^{\texttt{L}},U^{\texttt{R}})=(U_{j+1},U_j)$ (see Definition \ref{def:jdt}). 
Then ${\bf E}^n$ is a regular $\mf{sl}_n$-crystal with respect to $\mc{E}_j$ and $\mc{F}_j$ for $1\leq j\leq n-1$ \cite[Lemma 5.1]{K18-2}.
Furthermore, ${\bf E}^n$ can be viewed as an ${\mf l}$-crystal by identifying 
$(U_n,\dots, U_1)=U_1\otimes \dots \otimes U_n$ with respect to $\te_i$ and $\tf_i$ for $i\geq 1$. Then the operators ${\mc X}_j$ for $1\leq j\leq n-1$ commutes with $\te_i$ and $\tf_i$ so that ${\bf E}^{n}$ becomes a $({\mf l},\mf{sl}_n)$-bicrystal.

{Let us assume that $\mu\in \mc{P}({\rm O}_n)$ satisfies $n-2\mu'_1\geq 0$.} 
The case when $n-2\mu'_1<n$ will be discussed in subsection \ref{subsec:negative case}.

Consider the following embedding of sets given by 
\begin{equation}\label{eq:identification}
\xymatrixcolsep{3pc}\xymatrixrowsep{0pc}\xymatrix{
{\bf T}(\mu,n)  \ \ar@{->}[r] & \ \mathbf{E}^{n}  \\
\mathbf{T} = (T_l, \dots, T_1,T_0) \ar@{|->}[r] & (T_l^{\texttt{L}}, T_l^{\texttt{R}}, \dots, T_1^{\texttt{L}}, T_1^{\texttt{R}},T_0) }.
\end{equation}
We identify $\mathbf{T} = (T_l, \dots, T_1,T_0)\in {\bf T}(\mu,n)$ with its image ${\bf U}=(U_{2l},\dots,U_1,U_0)$ under \eqref{eq:identification} so that $T_0=U_0$ and $(T_{i+1}, T_i)$ is given by
\begin{equation*} 
(T_{i+1}, T_i)
=(U_{j+2}, U_{j+1}, U_j, U_{j-1})
=(T^{\texttt{L}}_{i+1},T^{\texttt{R}}_{i+1},T^{\texttt{L}}_i,T^{\texttt{R}}_i),
\end{equation*}
with $j=2i$ for $1\leq i\leq l-1$.

Let ${\bf T}\in {\bf H}^\circ(\mu,n)$ be given.
Our first step is to define an operator ${\mc S}_j$ on ${\bf T}$ for $j=2,4,\dots,2l-2$ by  
\begin{equation} \label{eq:definition_S}
\mathcal{S}_j = 
\begin{cases}
\mathcal{F}_j^{a_i}, & \text{if $U_{j+1}(1) < U_j(a_i)$},\\
\mathcal{E}_j\mathcal{E}_{j-1}\mathcal{F}_j^{a_i-1}\mathcal{F}_{j-1}, & \text{if $U_{j+1}(1) > U_j(a_i)$}.
\end{cases}
\end{equation}
Note that the operator $\mathcal{S}_j$ is well-defined by Corollary \ref{cor:highest_weight_element}.
We assume that $\mc{S}_j$ is the identity operator when $a_i=0$.
Let us describe $\mc{S}_j$ more explicitly.\vskip 2mm

\begin{lem}\label{lem:description of S}
Under the above hypothesis, we have
\begin{equation*}
{\mc S}_j{\bf U}=\left(\dots, U_{j+2},\td{U}_{j+1}, \td{U}_{j}, U_{j-1},\dots \right),
\end{equation*}
where
\begin{itemize}
\item[(i)] if $U_{j+1}(1) < U_j(a_i)$, then
\begin{equation*}
\td{U}_{j+1}=U_{j+1}\boxplus U_j^{\texttt{\em tail}},\quad
\td{U}_{j}=U_j\boxminus U_j^{\texttt{\em tail}},
\end{equation*}

\item[(ii)] if $U_{j+1}(1) > U_j(a_i)$, then
\begin{equation*}
\begin{split}
\td{U}_{j+1} &= \left(U_{j+1}(b_i+c_i), \dots, U_{j+1}(3) ) \boxplus ( U_{j+1}(2), \,U_j(a_i-1), \dots, U_{j}(1)\right), \\
\td{U}_{j} &=(U_j(a_i+c_i), \dots, U_j(a_i), U_{j+1}(1))\boxplus \emptyset.
\end{split}
\end{equation*}
\end{itemize}
\end{lem}
\pf Suppose that $U_{j+1}(1) < U_j(a_i)$. Then we have
\begin{equation*}
\mathcal{S}_j{\bf U}=\mathcal{F}_j^{a_i}{\bf U}
=(\dots,\mathcal{F}^{a_i}(U_{j+1}, U_i),\dots)
=(\dots,U_{j+1}\boxplus U_j^{\texttt{tail}}, U_j\boxminus U_j^{\texttt{tail}},\dots), 
\end{equation*}
which is given by cutting $U_j^{\texttt{tail}}$ and then putting it below $U_{j+1}$.\vskip 2mm

Next, suppose that $U_{j+1}(1) > U_j(a_i)$. 
By Lemma \ref{cor:highest_weight_element}, we have ${\mf r}_i{\mf r}_{i+1}=1$.
Then we have
$\mathcal{S}_j{\bf U}
=\mc{E}_j\mc{E}_{j-1}\mc{F}_j^{a_i-1}\mc{F}_{j-1}{\bf U}$.
Ignoring the components other than $(T_{i+1},T_i)$, we have
\begin{equation*} %
\begin{split}
\mc{E}_j\mc{E}_{j-1}\mc{F}_j^{a_i-1}\mc{F}_{j-1}&\left( T_{i+1}, T_i \right) \\
= &\mc{E}_j\mc{E}_{j-1}\mc{F}_j^{a_i-1}\left( U_{j+2}, U_{j+1},\mc{F}(U_j, U_{j-1}) \right) \\
= &\mc{E}_j\mc{E}_{j-1}\mc{F}_j^{a_i-1}\left( U_{j+2}, U_{j+1}, U^*_j, U^*_{j-1} \right) \\
= &\mc{E}_j\mc{E}_{j-1}\left( U_{j+2}, \mc{F}^{a_i-1}(U_{j+1}, U^*_j), U^*_{j-1} \right) \\
= &\mc{E}_j\mc{E}_{j-1}\left( U_{j+2}, U_{j+1} \boxplus {U}_j^{*\texttt{tail}}, U^*_j\boxminus{U}_j^{*\texttt{tail}}, U^*_{j-1} \right) \\
= &\mc{E}_j\left( U_{j+2}, U_{j+1}\boxplus {U}_j^{*\texttt{tail}}, \mc{E}(U^*_j\boxminus{U}_j^{*\texttt{tail}}, U^*_{j-1}) \right) \\
= &\mc{E}_j\left( U_{j+2}, U_{j+1}\boxplus {U}_j^{*\texttt{tail}}, U_j^{\uparrow}\boxminus{U}_j^{*\texttt{tail}}, U_{j-1} \right) \\
= & \left( U_{j+2}, \mc{E}(U_{j+1}\boxplus {U}_j^{*\texttt{tail}}, U_j^{\uparrow}\boxminus{U}_j^{*\texttt{tail}}), U_{j-1} \right) \\
= & \left( U_{j+2},\td{U}_{j+1},\td{U}_{j}, U_{j-1} \right),
\end{split}
\end{equation*} 
where 
\begin{equation*}\label{eq:calculation_operator_S2}
\begin{split}
& U^*_{j-1} = T^{\texttt{R}*}_i=(U_{j-1}(b_i+c_i), \dots, U_{j-1}(2)) \boxplus \emptyset, \\
& U^*_j = T^{\texttt{L}*}_i=(U_j(a_i+c_i), \dots, U_j(a_i), U_{j-1}(1)) \boxplus ( U_j(a_i-1), \dots, U_{j}(1)),\\
& U^{\uparrow}_j = (U_j(a_i+c_i), \dots, U_j(a_i) ) \boxplus ( U_j(a_i-1), \dots, U_{j}(1)).\\
\end{split}
\end{equation*}
This proves the lemma. 
\qed
\vskip 2mm

Let us call this algorithm to obtain $(\td{U}_{j+1},\td{U}_j)$ from $(U_{j+1},U_j)$ {\em sliding algorithm}.%

\begin{cor}\label{cor:spinor after S}
Under the above hypothesis, we have the following.
\begin{itemize}
\item[(1)] For $j=2$, there exists unique ${T}, S\in {\bf T}(0)$ such that 
$({T}^{\texttt{\em L}},{T}^{\texttt{\em R}})=(\td{U}_2, {U}_1)$ and
$({S}^{\texttt{\em L}},{S}^{\texttt{\em R}})=({U}_1, {U}_0)$ if $U_0$ is non-empty.

\item[(2)] For $j=2i$ with $1\leq i\leq l-1$, there exists a unique ${T}\in {\bf T}(a_i)$ such that 
$({T}^{\texttt{\em L}},{T}^{\texttt{\em R}})=(\td{U}_{j+1}, \td{U}_{j})$ and the residue of ${T}$ is 0  if $U_{j+1}(1) < U_j(a_i)$ and 1 if $U_{j+1}(1) > U_j(a_i)$.

\item[(3)] The pair $(U_{2l},\td{U}_{2l-1})$ forms a semistandard tableau when the columns are put together horizontally along L.
\end{itemize}

\end{cor}
\pf (1) and (3) follow directly from Definition \ref{def:admissibility} and the description of $(\td{U}_{j+1}, \td{U}_{j})$ in Lemma \ref{lem:description of S}.

By definition of $\mc{S}_j$, it is not difficult to see that 
$(\td{U}_{j+1}, \td{U}_{j})$ forms a semistandard tableau, say $T$ of shape $\la(a_i,b'_i,c'_i)$ for some $b'_i, c'_i\in\Z_+$ such that 
$(T^{\texttt{L}},T^{\texttt{R}})=(\td{U}_{j+1}, \td{U}_{j})$. 
The residue of $T$ follows immediately from the description of $(\td{U}_{j+1}, \td{U}_{j})$ in Lemma \ref{lem:description of S}. This proves (2).
\qed

\begin{cor}\label{cor:l-equivalence under S}
We have ${\mc S}_j{\mc S}_k={\mc S}_k{\mc S}_j$ for $j\neq k$, and ${\mc S}_2{\mc S}_4\dots{\mc S}_{2l-2}{\bf U}\equiv_{\mf l}{\bf U}$.
\end{cor}
\pf Since ${\mc S}_{j}$ changes only $U_j, U_{j+1}$ by Lemma \ref{lem:description of S}, it is clear that ${\mc S}_j{\mc S}_k={\mc S}_k{\mc S}_j$ for $j\neq k$. The ${\mf l}$-crystal equivalence follows from the fact that ${\bf E}^{n}$ is a $({\mf l},\mf{sl}_n)$-bicrystal.
\qed

\begin{ex} \label{ex:illustration of S}
\mbox{}

{\rm
(1) The following is an illustration of ${\mc S}_j$ when $U_{j+1}(1)<U_j(a_i)$.
\begin{equation*}
\hskip 0.7cm
\begin{split}
\ytableausetup {mathmode, boxsize=0.9em} 
\begin{ytableau}
	\none & \none & \none & \none & \none & \none & \none & \none\\
	\none & \none & \none & \none & \none  &  \tl{1} & \none\\
	\none & \none & \none & \none & \none &   \tl{2} & \none\\
	\none & \none & \tl{\blue 1} & \none &  \tl{\red 1} & \tl{3} & \none\\
	\none[\mathrel{\raisebox{-0.6ex}{$\scalebox{0.45}{\dots\dots}$}}] & \none & \tl{\blue {\bf 2}} & \none[\mathrel{\raisebox{-0.6ex}{$\scalebox{0.45}{\dots\dots}$}}] & \tl{\red 2} & \tl{4} & \none[\mathrel{\raisebox{-0.6ex}{$\scalebox{0.45}{\dots\dots}$}}] \\
	\none & \tl{1} & \none & \none &  \tl{\red {\bf 3}} & \none & \none\\
	\none & \tl{4} & \none & \none &  \tl{\red 5} & \none & \none\\
	\none & \none & \none & \none & \none & \none & \none & \none\\
	\none & \none & \none & \none & \none & \none & \none & \none\\
	\none & \none[] & \none[\quad \tl{$U_{j+1}$}] & \none & \none & \none[\tl{$U_j$} \quad ] & \none[] & \none
\end{ytableau}
\end{split} \quad {\xmapsto{\hspace*{1.5cm}}} \quad
\begin{split}
\ytableausetup {mathmode, boxsize=0.9em} 
\begin{ytableau}
	\none &\none & \none & \none & \none & \none & \none \\
	\none & \none & \none & \none & \none & \tl{1} & \none \\
	\none &\none & \none & \none & \none & \tl{2} & \none \\
	\none &\none & \tl{\blue 1} & \none & \tl{1} & \tl{3} & \none \\
	\none[\mathrel{\raisebox{-0.6ex}{$\scalebox{0.45}{\dots\dots}$}}] &\none & \tl{\blue 2} & \none[\mathrel{\raisebox{-0.6ex}{$\scalebox{0.45}{\dots\dots}$}}] & \tl{2} & \tl{4} & \none[\mathrel{\raisebox{-0.6ex}{$\scalebox{0.45}{\dots\dots}$}}] \\
	\none &\tl{1}  & \tl{\red 3}  & \none &  \none & \none & \none \\
	\none &\tl{4}  & \tl{\red 5}  & \none &  \none & \none & \none \\
	\none &\none & \none & \none & \none & \none & \none \\
	\none &\none & \none & \none & \none & \none & \none \\
	\none &\none[\ \ \ \ ] & \none[\quad \tl{$\td{U}_{j+1}$}] & \none & \none & \none[\tl{$\td{U}_j$}\quad ] & \none[] & \none
\end{ytableau}
\end{split}
\end{equation*}
	
(2) The following is an illustration of ${\mc S}_j$ when $U_{j+1}(1)>U_j(a_i)$.

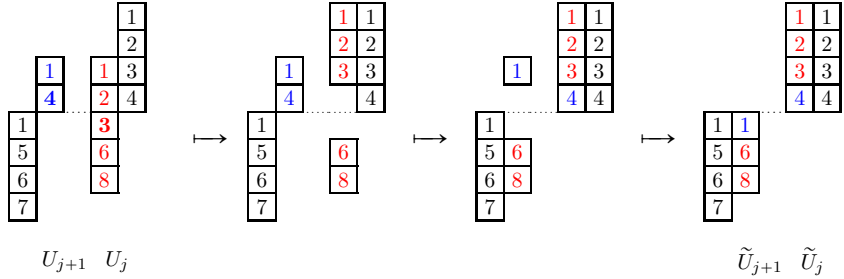
\begin{equation*}
\hskip 5mm
\begin{split}
\ytableausetup {mathmode, boxsize=0.9em} 
\begin{ytableau}
\none & \none & \none & \none & \tl{1} \\
\none & \none & \none & \none & \tl{2} \\
\none & \tl{\blue 1} & \none & \tl{\red 1} & \tl{3} \\
\none & \tl{\blue {\bf 4}} & \none[\mathrel{\raisebox{-0.6ex}{$\scalebox{0.45}{\dots\dots}$}}] & \tl{\red 2} & \tl{4} \\
\tl{1} & \none & \none &  \tl{\red {\bf 3}} & \none \\
\tl{5} & \none & \none &  \tl{\red 6} & \none \\
\tl{6} & \none & \none &  \tl{\red 8} & \none \\
\tl{7} & \none & \none &  \none & \none \\
\none & \none & \none & \none & \none & \none \\
\none[ ] & \none[\quad \tl{$U_{j+1}$}] & \none & \none & \none[\tl{$U_j$}\quad ] & \none[ ] 
\end{ytableau}
\end{split}\ \  {\xmapsto{\hspace*{0.25cm}}}\ \
\begin{split}
\ytableausetup {mathmode, boxsize=0.9em} 
\begin{ytableau}
\none & \none & \none & \tl{\red 1} & \tl{1} \\
\none & \none & \none & \tl{\red 2} & \tl{2} \\
\none & \tl{\blue 1} & \none & \tl{\red 3} & \tl{3} \\
\none & \tl{\blue 4} & \none[\mathrel{\raisebox{-0.6ex}{$\scalebox{0.4}{\dots\dots}$}}] & \none[\mathrel{\raisebox{-0.6ex}{$\scalebox{0.4}{\dots\dots}$}}] & \tl{4} \\
\tl{1} & \none & \none & \none & \none \\
\tl{5} & \none & \none & \tl{\red 6} & \none \\
\tl{6} & \none & \none & \tl{\red 8} & \none \\
\tl{7} & \none & \none & \none & \none \\
\none & \none & \none & \none & \none \\
\none & \none & \none & \none & \none \\
\end{ytableau}
\end{split}\ \ \xmapsto{\hspace*{0.25cm}}\ \
\begin{split}
\ytableausetup {mathmode, boxsize=0.9em} 
\begin{ytableau}
\none & \none & \none &  \tl{\red 1} & \tl{1} \\
\none & \none & \none &  \tl{\red 2} & \tl{2}  \\
\none & \tl{\blue 1} & \none & \tl{\red 3} & \tl{3} \\
\none & \none[\mathrel{\raisebox{-0.6ex}{$\scalebox{0.4}{\dots\dots}$}}] & \none[\mathrel{\raisebox{-0.6ex}{$\scalebox{0.4}{\dots\dots}$}}] & \tl{\blue 4} & \tl{4} \\
\tl{1} & \none & \none & \none & \none \\
\tl{5} & \tl{\red 6} & \none & \none & \none \\
\tl{6} & \tl{\red 8} & \none & \none & \none \\
\tl{7} & \none & \none & \none & \none \\
\none & \none & \none &  \none & \none \\
\none & \none & \none & \none & \none \\
\end{ytableau}
\end{split}\ \ \xmapsto{\hspace*{0.25cm}}\ \
\begin{split}
\ytableausetup {mathmode, boxsize=0.9em} 
\begin{ytableau}
\none & \none & \none &  \tl{\red 1} & \tl{1} \\
\none & \none & \none &  \tl{\red 2} & \tl{2} \\
\none & \none & \none & \tl{\red 3} & \tl{3} \\
\none & \none & \none[\mathrel{\raisebox{-0.6ex}{$\scalebox{0.45}{\dots\dots}$}}] & \tl{\blue 4} & \tl{4} \\
\tl{1} & \tl{\blue 1} & \none & \none & \none \\
\tl{5} & \tl{\red 6} & \none  & \none & \none \\
\tl{6} & \tl{\red 8} & \none  & \none & \none \\
\tl{7} & \none & \none & \none  & \none \\
\none & \none & \none & \none & \none & \none \\
\none[\ \ \ \ ] & \none[\quad \tl{$\td{U}_{j+1}$}] & \none & \none & \none[\tl{$\td{U}_j$}\quad ] & \none[]
\end{ytableau}
\end{split}
\end{equation*} 
	}
\end{ex}

Now consider
\begin{equation}\label{eq:S U}
\begin{split}
{\mc S}_2{\mc S}_4\dots{\mc S}_{2l-2}{\bf U}=(U_{2l},\td{U}_{2l-1},\td{U}_{2l-2},\dots,\td{U}_{3},\td{U}_{2},U_1,U_0)\in {\bf E}^{n},
\end{split}
\end{equation}
where $\td{U}_j$ is given in Lemma \ref{lem:description of S}. 
Let 
\begin{equation}\label{eq:tilde U}
\begin{split}
\td{\bf U}&= (\td{U}_{2l-1},\td{U}_{2l-2},\dots,\td{U}_{3},\td{U}_{2},U_1,U_0)\in {\bf E}^{n-1}.
\end{split}
\end{equation}

\begin{lem}\label{lem:sliding one step}
There exists a unique $\td{\bf T}\in {\bf T}(\td{\mu},n-1)$ 
corresponding to $\td{\bf U}$ under \eqref{eq:identification}, where $\td{\mu}=(\mu_2,\mu_3,\dots)$.
\end{lem}
\pf Let us define
\begin{equation}\label{eq:def of tilde T}
\td{\bf T}=
\begin{cases}
(\td{T}_{l-1},\dots,\td{T}_1,\td{T}_0), & \text{if $n=2l$},\\
(\td{T}_{l},\dots,\td{T}_1,\td{T}_0), & \text{if $n=2l+1$},\\
\end{cases}
\end{equation} 
as follows: 
\begin{itemize}
\item[(1)] if $n=2l$, then let $\td{T}_0=U_1$ and let $\td{T}_i\in {\bf T}(a_i)$ for $1\leq i\leq l-1$ such that 
$(\td{T}_i^{\texttt{L}},\td{T}_i^{\texttt{R}})=(\td{U}_{2i+1}, \td{U}_{2i}),$
given in Corollary \ref{cor:spinor after S}(2),
\item[(2)] if $n=2l+1$, then let $\td{T}_0=\emptyset$, $\td{T}_1\in {\bf T}(0)$ and $\td{T}_{i+1}\in {\bf T}(a_i)$ for $1\leq i\leq l-1$ such that 
$(\td{T}_1^{\texttt{L}},\td{T}_1^{\texttt{R}})=({U}_{1}, {U}_0)$,  
$(\td{T}_i^{\texttt{L}},\td{T}_i^{\texttt{R}})=(\td{U}_{2i+1}, \td{U}_{2i})$, 
given in Corollary \ref{cor:spinor after S}(1) and (2), respectively.
\end{itemize}

We have $\td{\bf T}\in \widehat{\bf T}(\td{\mu},n-1)$.
Let us show that $\td{\bf T}\in {\bf T}(\td{\mu},n-1)$.
For simplicity, let us assume that $n=2l$ since the proof for $n=2l+1$ is almost identical.

By Corollary \ref{cor:spinor after S}(1), we have $\td{T}_1\prec \td{T}_0$. 
So it suffices to show that $\td{T}_{i}\prec \td{T}_{i-1}$ for $2\leq i\leq l-1$.
This can be checked in a straightforward way 
using the fact that ${\bf T}\in {\bf H}^\circ(\mu,n)$ and Lemma \ref{lem:description of S} as follows.

Consider a triple $(T_{i+1}, T_i, T_{i-1})$ in ${\bf T}$. 
Recall that each $T_i$ satisfies (H1) and (H2).
Without loss of generality, let us consider $(T_3, T_2, T_1)$, which can be identified
with $(U_6, U_5, U_4, U_3, U_2, U_1)$ under the map \eqref{eq:identification}. 
Put
$$\mc{S}_{2}\mc{S}_{4}(T_{3}, T_2, T_{1}) = (\td{U}_6, \td{U}_5, \td{U}_4, \td{U}_3, \td{U}_2, \td{U}_1).$$  Note that $\td{U}_1 = U_1$ and $\td{U}_6=U_6$. Let $\lambda(a_j, b_j, c_j)$ be the shape of $T_j$ for $j=1, 2, 3$. Let $\td{T}_j$ be the tableau corresponding to $(\td{U}_{j+2}, \td{U}_{j+1})$ for $j=1, 2, 3$ in \eqref{eq:def of tilde T}.

We consider the following four cases. The other cases can be checked in a similar manner. 

{\bf\em Case 1. $(\mathfrak{r}_{3}, \mathfrak{r}_2 , \mathfrak{r}_{1})=(0,0,0)$.}
In this case, the operators $\mathcal{S}_2$ and $\mathcal{S}_4$ are just sliding $T_2^{\texttt{tail}}$ and $T_{1}^{\texttt{tail}}$ to the left horizontally.
Note that 
$\td{T}_1 \in {SST}(\lambda(a_{1}, c_{1}-b_2-c_2, b_2+c_2))$ and 
$\td{T}_2 \in {SST}(\lambda(a_2, c_2-b_{3}-c_{3}, b_{3}+c_{3}))$. 
It is straightforward to check that $\td{T}_2\prec \td{T}_1$.\vskip 2mm

{\bf\em Case 2. $(\mathfrak{r}_{3}, \mathfrak{r}_2 , \mathfrak{r}_{1})=(1,1,0)$.}
If $U_5(1) < U_4(a_2)$, then the proof is the same as in {\em Case 1}. 
So we assume that $U_5(1) > U_4(a_2)$.

Note that $\td{T}_1 \in {SST}(\lambda(a_{1}, c_{1}-b_2-c_2, b_2+c_2))$ and $\td{T}_2 \in SST(\lambda(a_2, c_2-b_{3}-c_{3}+4, b_{3}+c_{3}-2))$. 
By Corollary \ref{cor:spinor after S}(2), $\td{T}_2$ and $\td{T}_1$ have residue $1$ and $0$ respectively. We see that Definition \ref{def:admissibility}(1)-(i) holds on $(\td{T}_2,\td{T}_1)$. 

By \cite[Lemma 3.4]{K18-3}, we have $U_5(1)=\td{U}_4(1) \le \td{U}_3(a_1+1)=U_3(1)$, which together with (H1) on $T_2$ implies Definition \ref{def:admissibility}(1)-(ii) on $(\td{T}_2,\td{T}_1)$. 
Definition \ref{def:admissibility}(1)-(iii) on $(\td{T}_2,\td{T}_1)$ follows from 
the one on $(T_2, T_{1})$, (H1) and (H2) on $T_2$. Thus $\td{T}_2\prec\td{T}_1$. \vskip 2mm

{\bf\em Case 3. $(\mathfrak{r}_{3}, \mathfrak{r}_2 , \mathfrak{r}_{1})=(0,1,1)$.}
If $U_3(1) < U_2(a_{1})$, then the proof is the same as in {\em Case 1}. 
So we assume that $U_3(1) > U_2(a_1)$. 

Note that
$\td{T}_2 \in {SST}(\lambda(a_{1}, c_{1}-b_2-c_2+4, b_2+c_2-2))$ and 
$\td{T}_1 \in {SST}(\lambda(a_2, c_2-b_{3}-c_{3}, b_{3}+c_{3}))$. 
By Corollary \ref{cor:spinor after S}(2), $\td{T}_2$ and $\td{T}_1$ have residue $0$ and $1$ respectively. 
We see that Definition \ref{def:admissibility}(1)-(i) holds on $(\td{T}_2,\td{T}_1)$.
Definition \ref{def:admissibility}(1)-(ii) on $(\td{T}_2,\td{T}_1)$ follows from
(H1) on $T_2$.
Also, Definition \ref{def:admissibility}(1)-(iii) on $(\td{T}_2,\td{T}_1)$ 
follows from the one on $(T_2, T_1)$. Thus $\td{T}_2\prec\td{T}_1$. \vskip 2mm

{\bf\em Case 4.  $(\mathfrak{r}_{3}, \mathfrak{r}_2 , \mathfrak{r}_{1})=(1,1,1)$.}
If $U_3(1) < U_2(a_{1})$ or $U_5(1) < U_4(a_2)$, then the proof is the same as the one of {\em Case 1}-{\em Case 3}. 
So we assume that $U_3(1) > U_2(a_1)$ and $U_5(1) > U_4(a_2)$.

Note that 
$\td{T}_2 \in {SST}(\lambda(a_{1}, c_{1}-b_2-c_2+4, b_2+c_2-2))$ and 
$\td{T}_1 \in {SST}(\lambda(a_2, c_2-b_{3}-c_{3}+4, b_{3}+c_{3}-2))$ and both have residue $1$.
Since $\mathfrak{r}_2 = 1$, we have $b_2 \ge 2$ and $c_2+2 \le b_2+c_2$, 
which implies Definition \ref{def:admissibility}(1)-(i) on $(\td{T}_2,\td{T}_1)$.

Definition \ref{def:admissibility}(1)-(ii) and (iii) on $(\td{T}_2,\td{T}_1)$ follow from the same argument as in {\em Case 2}. Thus $\td{T}_2\prec\td{T}_1$.  

This completes the proof.
\qed

\begin{cor}\label{cor:sliding one step}
If ${\bf T}\in {\bf H}(\mu,n)$, then we have $\td{\bf T}\in {\bf H}(\td{\mu},n-1)$.
\end{cor}
\pf Since $\td{\bf T}$ corresponds to $\td{\bf U}$, we have $\td{\bf T}\in {\bf H}(\td{\mu},n-1)$ by Lemma \ref{lem:criterion_highest_weight_elt}.
\qed
\vskip 2mm

\subsection{Separation algorithm}\label{subsec:separation}

{Let us assume that $n-2\mu'_1\geq 0$.} 
Suppose that ${\bf T}\in {\bf H}(\mu,n)$ is given, which corresponds to ${\bf U}\in {\bf E}^n$ under \eqref{eq:identification}.

Let us define a tableau $\ov{\bf T}$ satisfying the following:
\begin{itemize}
\item[(S1)] $\ov{\bf T}$ is Knuth equivalent to ${\bf T}$, that is, $\ov{\bf T}\equiv_{\mf l}{\bf T}$,

\item[(S2)] $\ov{\bf T}^{\texttt{tail}}\in SST(\mu')$ and  $\ov{\bf T}^{\texttt{body}}\in SST((\delta')^\pi)$ for some $\delta\in \cP^{(2)}$.
\end{itemize}

We define $\ov{\bf T}$ inductively on $n$. If $n\leq 3$, then let $\ov{\bf T}$ is given by putting together the columns in ${\bf U}$ horizontally along L. 

Suppose that $n\geq 4$. 
By Corollary \ref{cor:sliding one step}, there exists a unique ${\bf S}\in {\bf H}(\td{\mu},n-1)$ corresponding to $\td{\bf U}$ in \eqref{eq:tilde U}. By induction hypothesis, there exists $\ov{\bf S}$ satisfying (S1) and (S2). Then we define $\ov{\bf T}$ to be the tableau obtained by putting together the leftmost column in ${\bf T}$ and $\ov{\bf S}$ (horizontally along L).

\begin{lem}\label{lem:ov{T} is well-defined}
The tableau $\ov{\bf T}$ satisfies {\em (S1)} and {\em (S2)}.
\end{lem}
\pf
By \eqref{eq:S U} and \eqref{eq:tilde U}, the leftmost columns in ${\bf T}$ and $\ov{\bf S}$ are $U_{2l}$ and $\td{U}_{2l-1}$, respectively. By Corollary \ref{cor:spinor after S}(3), we conclude that $\ov{\bf T}$ is semistandard. 

The condition (S1) holds since 
$$\ov{\bf T}\equiv_{\mf l}\ov{\bf S}\otimes U_{2l}\equiv_{\mf l} {\bf S}\otimes U_{2l}
\equiv_{\mf l}\td{\bf U}\otimes U_{2l}\equiv_{\mf l}{\bf U}\equiv_{\mf l}{\bf T}.$$
The condition (S2) follows directly from the definition of ${\mc S}_j$.
\qed \vskip 2mm

Roughly speaking, we may obtain $\ov{\bf T}$ from ${\bf T}$ by applying ${\mc S}_j$'s as far as possible so that no subtableau below L is movable to the left.

\begin{equation}\label{eq:shape after separation}
\raisebox{.9ex}{$\eta=$} 
\resizebox{.38\hsize}{!}{
\def\lr#1{\multicolumn{1}{|@{\hspace{.75ex}}c@{\hspace{.75ex}}|}{\raisebox{-.04ex}{$#1$}}}
\def\l#1{\multicolumn{1}{|@{\hspace{.75ex}}c@{\hspace{.75ex}}}{\raisebox{-.04ex}{$#1$}}}
\def\r#1{\multicolumn{1}{@{\hspace{.75ex}}c@{\hspace{.75ex}}|}{\raisebox{-.04ex}{$#1$}}}
\raisebox{-.6ex}
{$\begin{array}{cccccccccc}
\cline{8-9}
& & & & & & &\l{\ \ } & \r{ } & \\
\cline{6-7}
 & & & & & \l{\ \ } & & & \r{ }\\
\cline{5-5} 
& & & & \l{\ \ } &  &  & & \r{\!\!\!\!\!\!\!\!\!\!\!\!\!\!\!\!\!\!\!\!\!\!\!\!\!\!{}_{(\delta')^\pi} } \\
\cline{3-4}
& & \l{\ \ } & & & & & & \r{ } &  \\
\cline{2-2}\cdashline{1-10}[0.5pt/1pt]\cline{6-9}
& \l{\ \ } &  & & \r{\ \ }  &  & &  \\ 
\cline{4-5}
& \l{\ \ } & \r{\!\!\!\!\!\!\!\!\! {}_{\mu'}} & & &  & &  \\
& \l{\ \ } & \r{\ \ } & & & & &  \\  
\cline{3-3} 
& \lr{\ \ } & & & &  & &  \\  
\cline{2-2} 
\end{array}$}}\ \ \raisebox{-.7ex}{$L$}
\end{equation}

\begin{ex} \label{ex:3.13}
{\rm
Let $n=8$ and $\mu = (4, 3, 3, 2) \in \mathcal{P}({\rm O}_8)$. 
Let $\mathbf{T} = (T_4, T_3, T_2, T_1) \in {\bf H}(\mu,8)$ and ${\bf U}=(U_8,\dots,U_1)$ be given by
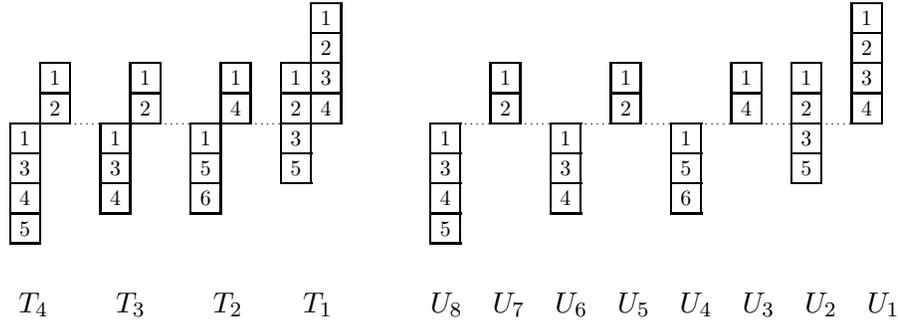
\begin{equation*}
\begin{split}
\ytableausetup {mathmode, boxsize=1.0em} 
&\begin{ytableau}
\none & \none & \none & \none & \none & \none & \none & \none & \none & \none & \tl{1} & \none  \\
\none & \none & \none & \none & \none & \none & \none & \none & \none & \none & \tl{2} & \none  \\
\none & \tl{1} & \none & \none & \tl{1} & \none & \none & \tl{1} & \none & \tl{1} & \tl{3} & \none \\
\none & \tl{2} & \none[\mathrel{\raisebox{-0.7ex}{$\scalebox{0.45}{\dots\dots}$}}] & \none & \tl{2} & \none[\mathrel{\raisebox{-0.7ex}{$\scalebox{0.45}{\dots\dots}$}}] & \none & \tl{4} & \none[\mathrel{\raisebox{-0.7ex}{$\scalebox{0.45}{\dots\dots}$}}]& \tl{2} & \tl{4} & \none \\
\tl{1} & \none & \none & \tl{1} & \none & \none & \tl{1} & \none & \none & \tl{3} & \none & \none  \\
\tl{3} & \none & \none & \tl{3} & \none & \none & \tl{5} & \none & \none & \tl{5} & \none & \none  \\
\tl{4} & \none & \none & \tl{4} & \none & \none & \tl{6} & \none & \none & \none & \none & \none  \\
\tl{5} & \none & \none & \none & \none & \none & \none & \none & \none & \none & \none & \none  \\
\none & \none & \none & \none & \none & \none & \none & \none & \none & \none & \none & \none \\
\end{ytableau} \quad\quad 
\begin{ytableau}
\none &\none &\none &\none &\none & \none & \none & \none & \none & \none & \none & \none & \none & \none & \tl{1} & \none  \\
\none &\none &\none &\none & \none & \none & \none & \none & \none & \none & \none & \none & \none & \none & \tl{2} & \none &\none  \\
\none & \none & \tl{1} & \none & \none &\none & \tl{1} & \none & \none &\none & \tl{1} & \none & \tl{1} &\none & \tl{3} & \none \\
\none &\none[\mathrel{\raisebox{-0.7ex}{$\scalebox{0.45}{\dots\dots}$}}] & \tl{2} & \none[\mathrel{\raisebox{-0.7ex}{$\scalebox{0.45}{\dots\dots}$}}] & \none & \none[\mathrel{\raisebox{-0.7ex}{$\scalebox{0.45}{\dots\dots}$}}] & \tl{2} & \none[\mathrel{\raisebox{-0.7ex}{$\scalebox{0.45}{\dots\dots}$}}] & \none &\none[\mathrel{\raisebox{-0.7ex}{$\scalebox{0.45}{\dots\dots}$}}] & \tl{4} & \none[\mathrel{\raisebox{-0.7ex}{$\scalebox{0.45}{\dots\dots}$}}] & \tl{2} & \none[\mathrel{\raisebox{-0.7ex}{$\scalebox{0.45}{\dots\dots}$}}] & \tl{4} & \none \\
\tl{1} & \none & \none &\none & \tl{1} & \none & \none &\none & \tl{1} & \none & \none &\none & \tl{3} & \none & \none  \\
\tl{3} & \none & \none &\none & \tl{3} & \none & \none &\none & \tl{5} & \none & \none &\none & \tl{5} & \none & \none  \\
\tl{4} & \none & \none &\none & \tl{4} & \none & \none &\none & \tl{6} & \none & \none & \none & \none & \none  \\
\tl{5} & \none & \none & \none & \none & \none & \none & \none & \none & \none & \none & \none  \\
\none & \none & \none & \none & \none & \none & \none & \none & \none & \none & \none & \none \\
\end{ytableau} \\
&\ T_4 \hskip 9mm T_3 \hskip 9mm T_2 \hskip 8mm T_1 \hskip 13mm 
U_8 \hskip 4mm U_7 \hskip 4mm U_6 \hskip 4mm U_5 \hskip 4mm U_4 \hskip 4mm U_3 \hskip 4mmU_2 \hskip 4mm U_1 \hskip 4mm
\end{split}
\end{equation*} 
  
Applying ${\mc S}_6{\mc S}_4{\mc S}_2$ to ${\bf U}$, we get
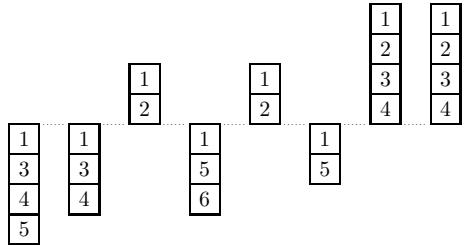
\begin{equation*}
\begin{split}
\ytableausetup {mathmode, boxsize=1.0em} 
\begin{ytableau}
\none &\none &\none &\none & \none & \none & \none & \none & \none & \none & \none & \none & \tl{1} &\none & \tl{1} & \none  \\
\none &\none &\none &\none & \none & \none & \none & \none & \none & \none & \none & \none & \tl{2} &\none & \tl{2} & \none  \\
\none & \none & \none &\none & \tl{1} & \none & \none &\none & \tl{1} & \none & \none &\none & \tl{3} &\none & \tl{3} & \none  \\
\none & \none[\mathrel{\raisebox{-0.7ex}{$\scalebox{0.3}{\dots\dots}$}}] & \none & \none[\mathrel{\raisebox{-0.7ex}{$\scalebox{0.3}{\dots\dots}$}}] & \tl{2} & \none[\mathrel{\raisebox{-0.7ex}{$\scalebox{0.3}{\dots\dots}$}}] & \none & \none[\mathrel{\raisebox{-0.7ex}{$\scalebox{0.3}{\dots\dots}$}}] & \tl{2} & \none[\mathrel{\raisebox{-0.7ex}{$\scalebox{0.3}{\dots\dots}$}}] & \none & \none[\mathrel{\raisebox{-0.7ex}{$\scalebox{0.3}{\dots\dots}$}}] & \tl{4} & \none[\mathrel{\raisebox{-0.7ex}{$\scalebox{0.3}{\dots\dots}$}}] & \tl{4} & \none \\
\tl{1} &\none & \tl{1} & \none & \none &\none & \tl{1} & \none & \none &\none & \tl{1} & \none & \none & \none & \none  \\
\tl{3} &\none & \tl{3} & \none & \none &\none & \tl{5} & \none & \none &\none & \tl{5} & \none & \none & \none & \none   \\
\tl{4} &\none & \tl{4} & \none & \none &\none & \tl{6} & \none & \none & \none & \none & \none & \none & \none   \\
\tl{5} & \none & \none & \none & \none & \none & \none & \none & \none & \none & \none & \none  \\
\none & \none & \none & \none & \none & \none & \none & \none & \none & \none & \none & \none  \\
\end{ytableau}
\end{split}
\end{equation*}  
The sequence of columns except the leftmost one (in gray) corresponds to ${\bf S}\in {\bf H}(\td{\mu},7)$ with $\td{\mu}=(3,3,2)$.
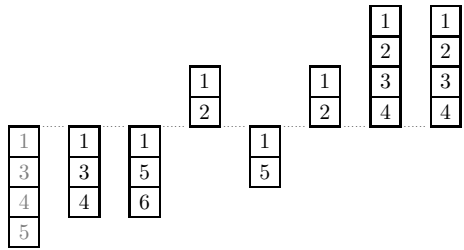
\begin{equation*}
\begin{split}
\ytableausetup {mathmode, boxsize=1.0em} 
\begin{ytableau}
\none &\none &\none &\none & \none & \none & \none & \none & \none & \none & \none & \none & \tl{1} &\none & \tl{1} & \none    \\
\none &\none &\none &\none & \none & \none & \none & \none & \none & \none & \none & \none & \tl{2} &\none & \tl{2} & \none   \\
\none &\none & \none & \none & \none &\none & \tl{1} & \none & \none &\none & \tl{1} & \none & \tl{3} &\none & \tl{3} & \none   \\
\none & \none[\mathrel{\raisebox{-0.7ex}{$\scalebox{0.3}{\dots\dots}$}}] & \none & \none[\mathrel{\raisebox{-0.7ex}{$\scalebox{0.3}{\dots\dots}$}}] &\none & \none[\mathrel{\raisebox{-0.7ex}{$\scalebox{0.3}{\dots\dots}$}}] &\tl{2} & \none[\mathrel{\raisebox{-0.7ex}{$\scalebox{0.3}{\dots\dots}$}}] & \none & \none[\mathrel{\raisebox{-0.7ex}{$\scalebox{0.3}{\dots\dots}$}}] & \tl{2} & \none[\mathrel{\raisebox{-0.7ex}{$\scalebox{0.3}{\dots\dots}$}}] & \tl{4} & \none[\mathrel{\raisebox{-0.7ex}{$\scalebox{0.3}{\dots\dots}$}}] & \tl{4} & \none  \\
\tl{\color{gray} 1} &\none & \tl{1} & \none & \tl{1} & \none & \none &\none & \tl{1} & \none & \none & \none & \none   \\
\tl{\color{gray} 3} &\none & \tl{3} & \none & \tl{5} & \none & \none &\none & \tl{5} & \none & \none & \none & \none & \none   \\
\tl{\color{gray} 4} &\none & \tl{4} & \none & \tl{6} & \none & \none & \none & \none & \none & \none & \none & \none   \\
\tl{\color{gray} 5} & \none & \none & \none & \none & \none & \none & \none & \none & \none & \none & \none   \\
\none & \none & \none & \none & \none & \none & \none & \none & \none & \none & \none & \none   \\
\end{ytableau} 
\end{split}
\end{equation*} 

Applying this process again to ${\bf S}$, we get
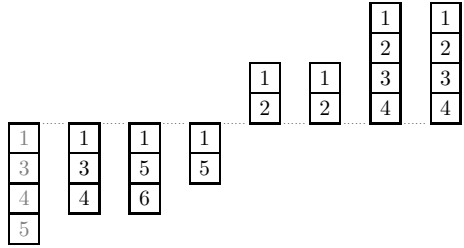
\begin{equation*}
\begin{split}
\ytableausetup {mathmode, boxsize=1.0em} 
\begin{ytableau}
\none &\none &\none &\none & \none & \none & \none & \none & \none & \none & \none & \none & \tl{1} &\none & \tl{1} & \none    \\
\none &\none &\none &\none & \none & \none & \none & \none & \none & \none & \none & \none & \tl{2} &\none & \tl{2} & \none   \\
\none &\none & \none & \none & \none &\none & \none & \none & \tl{1} &\none & \tl{1} & \none & \tl{3} &\none & \tl{3} & \none   \\
\none & \none[\mathrel{\raisebox{-0.7ex}{$\scalebox{0.3}{\dots\dots}$}}] & \none & \none[\mathrel{\raisebox{-0.7ex}{$\scalebox{0.3}{\dots\dots}$}}] &\none & \none[\mathrel{\raisebox{-0.7ex}{$\scalebox{0.3}{\dots\dots}$}}] &\none & \none[\mathrel{\raisebox{-0.7ex}{$\scalebox{0.3}{\dots\dots}$}}] & \tl{2} &\none[\mathrel{\raisebox{-0.7ex}{$\scalebox{0.3}{\dots\dots}$}}]  & \tl{2} & \none[\mathrel{\raisebox{-0.7ex}{$\scalebox{0.3}{\dots\dots}$}}]  & \tl{4} & \none[\mathrel{\raisebox{-0.7ex}{$\scalebox{0.3}{\dots\dots}$}}] & \tl{4} & \none  \\
\tl{\color{gray} 1} &\none & \tl{1} & \none & \tl{1} & \none & \tl{1} &\none & \none & \none & \none & \none & \none   \\
\tl{\color{gray} 3} &\none & \tl{3} & \none & \tl{5} & \none & \tl{5} &\none & \none & \none & \none & \none & \none & \none   \\
\tl{\color{gray} 4} &\none & \tl{4} & \none & \tl{6} & \none & \none & \none & \none & \none & \none & \none & \none   \\
\tl{\color{gray} 5} & \none & \none & \none & \none & \none & \none & \none & \none & \none & \none & \none   \\
\none & \none & \none & \none & \none & \none & \none & \none & \none & \none & \none & \none   \\
\end{ytableau} 
\end{split}
\end{equation*} 
Therefore, $\ov{\bf T}$ is given by
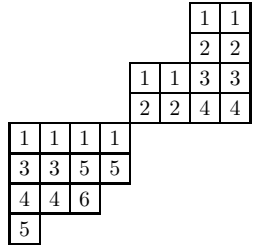
\begin{equation*}
\ytableausetup {mathmode, boxsize=1.0em} 
\begin{ytableau}
\none & \none & \none & \none & \none & \none & \tl{1} & \tl{1}   \\
\none & \none & \none & \none & \none & \none & \tl{2} & \tl{2}   \\
\none & \none &\none & \none & \tl{1} & \tl{1} & \tl{3} & \tl{3}  \\
\none &\none & \none &\none & \tl{2} & \tl{2} & \tl{4} & \tl{4}   \\
\tl{1} & \tl{1} & \tl{1} & \tl{1} &\none & \none & \none & \none   \\
\tl{3} & \tl{3} & \tl{5} & \tl{5} &\none & \none & \none & \none   \\
\tl{4} & \tl{4} & \tl{6} & \none & \none & \none & \none & \none  \\ 
\tl{5} & \none & \none & \none & \none & \none & \none & \none     \\
\none & \none & \none & \none & \none & \none & \none & \none      \\
\end{ytableau}
\end{equation*}

}
\end{ex}
\vskip 5mm
 
\begin{prop}\label{prop:body and tail}
Let ${\bf T}\in {\bf H}(\mu,n)$ be given with $n-2\mu'_1\geq 0$. 
Then
\begin{itemize}
\item[(1)] $\ov{\bf T}^{\texttt{\em body}}=H_{(\delta')^\pi}$ for some $\delta\in \cP^{(2)}$,

\item[(2)] $\ov{\bf T}^{\texttt{\em tail}}\in \texttt{\em LR}^{\la'}_{\delta' \mu'}$ if ${\bf T}\equiv_{\mf l} H_{\la'}$ for some $\la\in \cP$.
\end{itemize}
\end{prop}
\pf (1) We have ${\bf T}\equiv_{\mf l}\ov{\bf T}$, 
and $\ov{\bf T}\equiv_{\mf l} \ov{\bf T}^{\texttt{body}}\otimes \ov{\bf T}^\texttt{tail}$. Since ${\rm sh}(\ov{\bf T}^{\texttt{body}})=(\delta')^\pi$ for some $\delta\in\cP^{(2)}$, we should have $\ov{\bf T}^{\texttt{body}}=H_{(\delta')^\pi}$. 
(2) It follows from the fact that 
$\ov{\bf T}\equiv_{\mf l} H_{(\delta')^\pi}\otimes \ov{\bf T}^{\texttt{tail}}
\equiv_{\mf l}\left(\ov{\bf T}^{\texttt{tail}}\rightarrow H_{(\delta')^\pi}\right)=H_{\la'}$.
\qed 
\vskip 2mm

Let us call this combinatorial algorithm to obtain $(\ov{\bf T}^{\texttt{body}},\ov{\bf T}^{\texttt{tail}})$ from ${\bf T}$ {\em separation algorithm}.

\subsection{Separation algorithm when $n-2\mu'_1<0$}\label{subsec:negative case}
In this subsection, we consider the separation in the case $n-2\mu'_1<0$.
In this case, we need to deal with the tableaux with odd height in $\ov{\bf T}(0)$ and ${\bf T}^{\texttt{sp}-}$.
This is the reason why we consider the separation in this case separately. 
\vskip 2mm

Let us assume that $n-2\mu'_1< 0$.
Recall that $\overline{\mu} = (\ov{\mu}_i) \in \cP$ be such that $(\ov{\mu}')_1 = n-\mu'_1$ and $(\ov{\mu}')_i = \mu'_i$ for $i \ge 2$. 

Let ${\bf T}\in {\bf T}(\mu,n)$ be given. 
Suppose that $n=2l+r$, where $l\geq 1$ and $r=0,1$.
By \eqref{eq:hat{T}}, we have
\begin{equation}\label{eq:T notation-2}
{\bf T} = (T_l, \dots,T_{m+1},T_m,\dots,T_1,T_0),
\end{equation}
where $T_i\in {\bf T}(a_i)$ for some $a_i\in\Z_+$ ($m+1\leq i\leq l$), $T_i\in \ov{\bf T}(0)$ ($1\leq i\leq m$), and  $T_0\in {\bf T}^{\rm sp-}$ (resp. $T_0=\emptyset$) when $r=1$ (resp. $r=0$). Here $m=q_-$ in \eqref{eq:hat{T}}.
Under \eqref{eq:identification}, we identify ${\bf T}$ with
\begin{equation*}
{\bf U}=(U_{2l},\dots,U_{2m+1},U_{2m},\dots,U_1,U_0).
\end{equation*}
We may also assume that $U_i\in {\bf T}^{\textrm{sp}-}$ for $0\leq i\leq 2m$.
The following is an analogue of Definition \ref{def:pseudo_H} when $n-2\mu'_1<0$.
  
\begin{df}\label{def:pseudo_H-2}
{\rm
Let ${\bf H}^\circ(\mu,n)$ be the set of ${\bf T}\in {\bf T}(\mu,n)$ such that
\begin{itemize}
\item[(1)] $U_i[k]=k$ ($k\geq 1$) for $0\leq i\leq 2m$,

\item[(2)] $T_i$ satisfies (H1) and (H2) in Definition \ref{def:pseudo_H}  for $m+1\leq i\leq l$.
\end{itemize}
}
\end{df}  
  
\begin{prop}\label{prop:highest_weight_vectors-2}
We have
${\bf{H}}(\mu, n) \subset {\bf{H}}^\circ(\mu, n).$
\end{prop}
\pf Let ${\bf T}\in {\bf H}(\mu,n)$ be given.
 By Lemma \ref{lem:criterion_highest_weight_elt} and the admissibility of $T_{i+1}\prec T_i$ for $0\leq i\leq m-1$, we have $U_i[k]=k$ ($k\geq 1$) for $0\leq i\leq 2m$. Hence ${\bf T}$ satisfies (1).
The condition (2) can be verified by almost the same argument as in Proposition \ref{prop:highest_weight_vectors}.\qed

Now, let us define the tableau $\ov{\bf T}$ satisfying (S1) and (S2) in subsection \ref{subsec:separation}. We use induction on $n$.

Suppose that $n=3$, that is, ${\bf U}=(U_2,U_1,U_0)$. Take $U = H_{(1^a)}$ for some sufficiently large $a>0$. Then there exists a unique $T_0^\sharp\in {\bf T}(1)$ such that $((T_0^\sharp)^{\texttt{L}},(T_0^\sharp)^{\texttt{R}})=(U_0,U)$. 
Let ${\bf T}^\sharp=(T_1,T_0^\sharp)$. One can check that $T_1\prec T_0^\sharp$, and hence ${\bf T}^\sharp\in {\bf H}(\mu,4)$. 
Let $\ov{{\bf T}^\sharp}$ be the tableau defined in subsection \ref{subsec:separation}.
Then we define $\ov{\bf T}$ to be the one obtained from $\ov{{\bf T}^\sharp}$ by removing the rightmost column $U$. 
By Lemma \ref{lem:description of S}, $\ov{\bf T}$ does not depend on the choice of $U$ for all sufficiently large $a$, and hence is well-defined.

\begin{ex} \label{ex:spin minus}
{\rm Let $n=3$ and $\mu=(2,1)$.
\begin{equation*}
\quad {\bf T}=\hskip -8mm
\begin{split}
\ytableausetup{mathmode,boxsize=0.9em} 
\begin{ytableau}
\none & \none & \none & \none \\
\none & \none & \none & \none \\
\none & \none & \none & \none \\
\none & \tl{1} & \none & \tl{1} \\
\none & \tl{2} & \none & \tl{2} \\
\tl{1} & \tl{3} & \none & \tl{3} \\
\tl{2} & \tl{6} & \none[\mathrel{\raisebox{-0.6ex}{$\scalebox{0.3}{\dots\dots}$}}] & \tl{4} \\
\tl{3} & \none & \none &  \tl{5} \\
\tl{7} & \none & \none & \none \\
\none & \none & \none & \none \\
\end{ytableau}
\end{split}\ \xrightarrow{\hspace*{0.3cm}}\ {\bf T}^\sharp=\
\begin{split}
\ytableausetup {mathmode, boxsize=0.9em} 
\begin{ytableau}
\none & \none & \none &   \none &  \none[\color{gray} \tl{1}] \\
\none & \none & \none &   \none &  \none[\color{gray} \tl{2}] \\
\none & \none & \none &   \none &  \none[\color{gray} \tl{3}] \\
\none & \none & \none &   \none &  \none[\color{gray} \tl{4}]  \\
\none & \tl{1} & \none &   \tl{1} &   \none[\color{gray} \tl{5}]  \\
\none & \tl{2} & \none &  \tl{2} &  \none[\color{gray} \tl{6}]  \\
\tl{1}  & \tl{3} & \none &  \tl{3} &  \none[\color{gray} \tl{7}] \\
\tl{2}  & \tl{6} & \none[\mathrel{\raisebox{-0.6ex}{$\scalebox{0.3}{\dots\dots}$}}] &  \tl{4} &  \none[\color{gray} \tl{8}] \\
\tl{3}  & \none & \none &  \tl{5} &  \none \\
\tl{7}  & \none & \none &   \none &  \none  \\
\none & \none & \none &  \none &  \none  \\
\end{ytableau}
\end{split}\ \overset{\mathcal{S}}{\xrightarrow{\hspace*{0.3cm}}}\ \ov{{\bf T}^\sharp}=\
\begin{split}
\ytableausetup {mathmode, boxsize=0.9em} 
\begin{ytableau}
\none & \none &  \none &  \none[\color{gray} \tl{$1$}] \\
\none & \none &  \none &  \none[\color{gray} \tl{$2$}] \\
\none & \none &  \tl{1} &  \none[\color{gray} \tl{$3$}] \\
\none & \none &  \tl{2} &  \none[\color{gray} \tl{$4$}]  \\
\none & \none &  \tl{3} &  \none[\color{gray} \tl{$5$}]  \\
\none & \none &  \tl{4} &  \none[\color{gray} \tl{$6$}]  \\
\tl{1} & \tl{1} & \tl{5} &  \none[\color{gray} \tl{$7$}] \\
\tl{2} & \tl{2} & \tl{6} &  \none[\color{gray} \tl{$8$}] \\
\tl{3} & \tl{3} &   \none &  \none \\
\tl{7} & \none &   \none &  \none  \\
\none & \none &    \none &  \none  \\
\end{ytableau}
\end{split}\  \xrightarrow{\hspace*{0.3cm}}\  \ov{\bf T} = \
\begin{split}
\ytableausetup {mathmode, boxsize=0.9em} 
\begin{ytableau}
\none & \none & \none \\
\none & \none & \none \\
\none & \none & \tl{1} \\
\none & \none & \tl{2} \\
\none & \none & \tl{3}  \\
\none & \none & \tl{4} \\
\tl{1} & \tl{1} & \tl{5} \\
\tl{2} & \tl{2} & \tl{6} \\
\tl{3} & \tl{3} & \none \\
\tl{7} & \none & \none \\
\none & \none & \none \\
\end{ytableau} 
\end{split} 
\end{equation*}
}
\end{ex}

Suppose that $n\geq 4$. 
Let $${\bf V}=(U_{2l},\dots,U_{2m+1},U_{2m},U),$$ where $U=H_{(1^a)}$ for a sufficiently large $a$.
Note that there exists $T\in {\bf T}(1)$ such that $(T^{\texttt L},T^{\texttt R})=(U_{2m},U)$.
Since ${\bf U}\in {\bf H}^\circ(\mu,n)$ by Lemma \ref{prop:highest_weight_vectors-2}, we have ${\bf V}\in {\bf H}^\circ(\eta,N)$, where $\eta=\ov{\mu}\cup\{1\}=(\ov{\mu}_1,\dots,\ov{\mu}_{\ell},1)$ with $\ell=\ell(\ov{\mu})=l-m$ and $N=2\ell+2$.
Since $N-2\eta'_1=2\ell+2-2(\mu'_1+1)=0$, we may apply the sliding algorithm in subsection \ref{subsec:sliding} (see \eqref{eq:S U}) to have 
$$\widetilde{{\bf V}}=(\td{U}_{2l},\td{U}_{2l-1}\dots,\td{U}_{2m+1},\td{U}_{2m},U).$$ 
Let
$$\td{\bf U}=(\td{U}_{2l-1}\dots,\td{U}_{2m+1},\td{U}_{2m},U_{2m-1},\dots,U_1,U_0)\in {\bf E}^{n-1}.$$
We have the following analogue of Lemma \ref{lem:sliding one step}.

\begin{lem}\label{lem:sliding one step-2}
Under the above hypothesis, there exists a unique $\td{\bf T}\in {\bf H}(\td{\mu},n-1)$ corresponding to $\td{\bf U}$ under \eqref{eq:identification}, where $\td{\mu}=(\mu_2,\mu_3,\dots)$.
\end{lem}
\pf For $m\leq i\leq l-1$, there exists a unique $\td{T}_i\in {\bf T}(a_i)$ ($a_i\in\Z_+$) such that  
$(\td{T}_i^{\texttt{L}},\td{T}_i^{\texttt{R}})=(\td{U}_{2i+1}, \td{U}_{2i})$ by Lemma \ref{lem:sliding one step}. Also, for each $0\leq i\leq 2m-2$, there exists a unique $T\in \ov{\bf T}(0)$ such that $(\td{T}^{\texttt{L}},\td{T}^{\texttt{R}})=({U}_{i+1}, {U}_i)$ since $T_{j+1}\prec T_j$ for $0\leq j\leq 2m-1$ \eqref{eq:T notation-2}. 

Let us define
\begin{equation*} %
\td{\bf T}=(\td{T}_{l-1},\dots,\td{T}_m,\td{T}_{m-1},\dots,\td{T}_1,\td{T}_0)
\end{equation*} 
as follows:
\begin{itemize}
\item[(1)] let $\td{T}_i\in {\bf T}(a_i)$ such that $(\td{T}_i^{\texttt{L}},\td{T}_i^{\texttt{R}})=(\td{U}_{2i+1}, \td{U}_{2i})$ for $m\leq i\leq l-1$,

\item[(2)] if $n=2l$, then let $\td{T}_0=U_1$ and $\td{T}_i\in \ov{\bf T}(0)$ for $1\leq i\leq m-1$ such that 
$(\td{T}_i^{\texttt{L}},\td{T}_i^{\texttt{R}})=({U}_{2i+1},{U}_{2i})$,

\item[(3)] if $n=2l+1$, then let $\td{T}_0=\emptyset$ and $\td{T}_{i}\in \ov{\bf T}(0)$ for $1\leq i\leq m$ such that 
$(\td{T}_i^{\texttt{L}},\td{T}_i^{\texttt{R}})=({U}_{2i-1},{U}_{2i-2})$.
\end{itemize}

It is clear that $\td{T}_{i+1}\prec \td{T}_i$ for $0\leq i\leq m-2$.
By Lemma \ref{lem:sliding one step}, we have $\td{T}_{i+1}\prec \td{T}_i$ for $m\leq i\leq l-2$. Finally, by definition of $\td{\bf V}$, one can check without difficulty that $\td{T}_m\prec \td{T}_{m-1}$. Therefore, $\td{\bf T}\in {\bf T}(\td{
\mu},n-1)$. Since $\td{\bf T}$ is also an $\mf l$-highest weight element, we have $\td{\bf T}\in {\bf H}(\td{\mu},n-1)$ by Lemma \ref{lem:criterion_highest_weight_elt}.
\qed\vskip 2mm

Let ${\bf S}=\td{\bf T}\in {\bf H}(\td{\mu},n-1)$ in Lemma \ref{lem:sliding one step}. By induction hypothesis, there exists $\ov{\bf S}$ satisfying (S1) and (S2). Then we define $\ov{\bf T}$ to be the tableau obtained by putting together the leftmost column in ${\bf T}$ and $\ov{\bf S}$ (horizontally along L).

By the same arguments as in Lemma \ref{lem:ov{T} is well-defined}, we conclude that $\ov{\bf T}$ satisfies (S1) and (S2) in subsection \ref{subsec:separation}. 
Hence Proposition \ref{prop:body and tail} also holds in this case as follows by the same arguments.

\begin{prop}\label{prop:body and tail-2}
Let ${\bf T}\in {\bf H}(\mu,n)$ be given with $n-2\mu'_1<0$. 
Then
\begin{itemize}
\item[(1)] $\ov{\bf T}^{\texttt{\em body}}=H_{(\delta')^\pi}$ for some $\delta\in \cP^{(2)}$,

\item[(2)] $\ov{\bf T}^{\texttt{\em tail}}\in \texttt{\em LR}^{\la'}_{\mu'\delta'}$ if ${\bf T}\equiv_{\mf l} H_{\la'}$ for some $\la\in \cP$.
\end{itemize}
\end{prop}
\qed

\begin{ex} \label{ex:separation for negative case}
{\rm
Let $n=9$ and $\mu = (4, 3, 3, 2, 1) \in \mathcal{P}({\rm O}_9)$. 
We have $n-2\mu'_1 < 0$ and $\overline{\mu} = (4, 3, 3, 2)$. 
Let $\mathbf{T} \in {\bf T}(\mu, 9)$ be given by

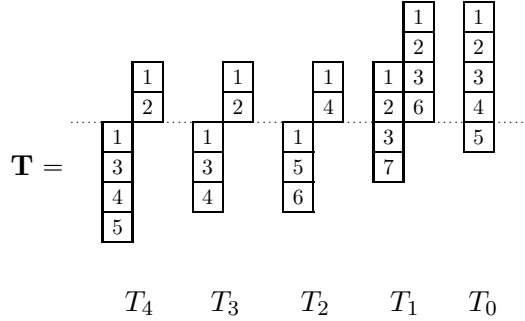
\begin{equation*}
\hskip 3cm{\bf T} = \hskip -4cm
\begin{split}
\ytableausetup {mathmode, boxsize=1.0em} 
\begin{ytableau}
\none & \none & \none & \none & \none & \none & \none & \none & \none & \none & \none & \tl{1} & \none &   \tl{1}  \\
\none & \none & \none & \none & \none & \none & \none & \none & \none & \none & \none & \tl{2} & \none &  \tl{2}  \\
\none & \none & \tl{1} & \none & \none & \tl{1} & \none & \none & \tl{1} & \none & \tl{1} & \tl{3} & \none &  \tl{3} \\
\none[\mathrel{\raisebox{-0.7ex}{$\scalebox{0.45}{\dots\dots}$}}] & \none & \tl{2} & \none[\mathrel{\raisebox{-0.7ex}{$\scalebox{0.45}{\dots\dots}$}}] & \none & \tl{2} & \none[\mathrel{\raisebox{-0.7ex}{$\scalebox{0.45}{\dots\dots}$}}] & \none & \tl{4} &\none[\mathrel{\raisebox{-0.7ex}{$\scalebox{0.45}{\dots\dots}$}}] & \tl{2} & \tl{6} & \none[\mathrel{\raisebox{-0.7ex}{$\scalebox{0.45}{\dots\dots}$}}] &  \tl{4} & \none[\mathrel{\raisebox{-0.7ex}{$\scalebox{0.45}{\dots\dots}$}}]\\
\none & \tl{1} & \none & \none & \tl{1} & \none & \none & \tl{1} & \none & \none & \tl{3} & \none & \none &  \tl{5} \\
\none & \tl{3} & \none & \none & \tl{3} & \none & \none & \tl{5} & \none & \none & \tl{7} & \none & \none &   \none \\
\none & \tl{4} & \none & \none & \tl{4} & \none & \none & \tl{6} & \none & \none & \none & \none & \none &  \none \\
\none & \tl{5} & \none & \none & \none & \none & \none & \none & \none & \none & \none & \none & \none &   \none \\
\none & \none & \none & \none & \none & \none & \none & \none & \none & \none & \none & \none & \none &   \none \\
\end{ytableau} \\
T_4 \hskip 7.5mm T_3 \hskip 8mm T_2 \hskip 8mm T_1 \hskip 6.5mm T_0 \hskip 3.5mm
\end{split}
\end{equation*}
It corresponds to 
\begin{equation*}
\hskip 2cm{\bf U} = \hskip -3cm
\begin{split}
\ytableausetup {mathmode, boxsize=1.0em} 
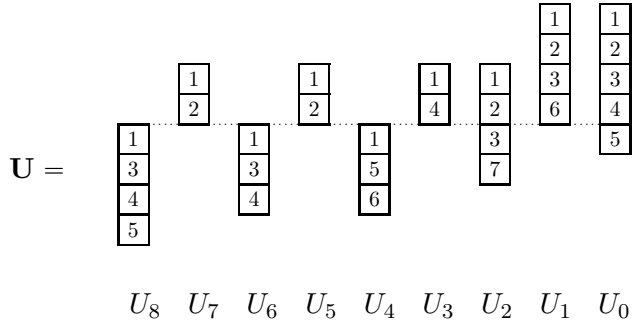
\begin{ytableau}
\none & \none & \none & \none &\none & \none & \none & \none & \none & \none & \none & \none & \none & \none & \tl{1} & \none &   \tl{1}  \\
\none & \none & \none & \none &\none & \none & \none & \none & \none & \none & \none & \none & \none & \none & \tl{2} & \none &  \tl{2}  \\
\none & \none & \tl{1} & \none & \none & \none &  \tl{1} & \none & \none & \none &  \tl{1} & \none & \tl{1} & \none &  \tl{3} & \none &  \tl{3} \\
\none & \none[\mathrel{\raisebox{-0.7ex}{$\scalebox{0.45}{\dots\dots}$}}] & \tl{2} & \none[\mathrel{\raisebox{-0.7ex}{$\scalebox{0.45}{\dots\dots}$}}] & \none & \none[\mathrel{\raisebox{-0.7ex}{$\scalebox{0.45}{\dots\dots}$}}] & \tl{2} & \none[\mathrel{\raisebox{-0.7ex}{$\scalebox{0.45}{\dots\dots}$}}] & \none & \none[\mathrel{\raisebox{-0.7ex}{$\scalebox{0.45}{\dots\dots}$}}] & \tl{4} & \none[\mathrel{\raisebox{-0.7ex}{$\scalebox{0.45}{\dots\dots}$}}] & \tl{2} & \none[\mathrel{\raisebox{-0.7ex}{$\scalebox{0.45}{\dots\dots}$}}] &  \tl{6} & \none[\mathrel{\raisebox{-0.7ex}{$\scalebox{0.45}{\dots\dots}$}}] &  \tl{4}\\
\tl{1} & \none & \none & \none & \tl{1} & \none & \none & \none &  \tl{1} & \none & \none & \none & \tl{3} & \none & \none & \none & \tl{5} \\
\tl{3} & \none & \none & \none &  \tl{3} & \none & \none & \none &  \tl{5} & \none & \none & \none & \tl{7} & \none & \none &   \none \\
\tl{4} & \none & \none & \none & \tl{4} & \none & \none & \none & \tl{6} & \none & \none & \none & \none & \none &  \none \\
\tl{5} & \none & \none & \none & \none & \none & \none & \none & \none & \none & \none & \none &   \none \\
\none & \none & \none & \none & \none & \none & \none & \none & \none & \none & \none & \none &   \none \\
\end{ytableau} \\
U_8 \hskip 3.5mm U_7 \hskip 3.5mm U_6 \hskip 3.5mm U_5 \hskip 3.5mm
U_4 \hskip 3.5mm U_3 \hskip 3.5mm U_2 \hskip 3.5mm U_1 \hskip 3.5mm U_0
\end{split}
\end{equation*} 
Putting $U=H_{(1^8)}$ at the rightmost column and applying the sliding algorithm, we get
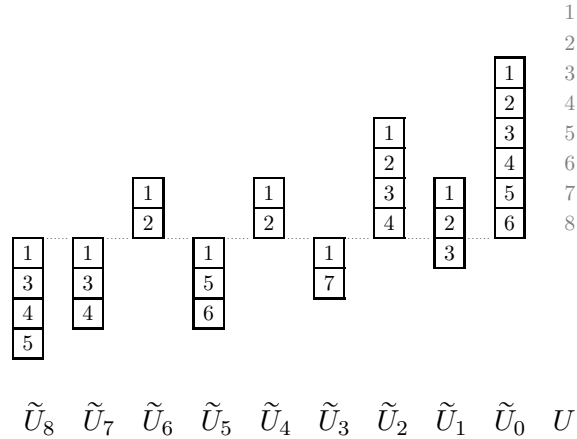
\begin{equation*}
\begin{split}
\ytableausetup {mathmode, boxsize=1.0em} 
\begin{ytableau}
\none & \none & \none & \none & \none & \none & \none & \none & \none & \none & \none & \none & \none & \none & \none & \none & \none & \none &   \none[\color{gray} \tl{1}]\\
\none & \none & \none & \none & \none & \none & \none & \none & \none & \none & \none & \none & \none & \none & \none & \none & \none & \none &   \none[\color{gray} \tl{2}] \\
\none & \none & \none & \none & \none & \none & \none & \none & \none & \none & \none & \none & \none & \none & \none & \none & \tl{1} & \none &   \none[\color{gray} \tl{3}]\\
\none & \none & \none & \none & \none & \none & \none & \none & \none & \none & \none & \none & \none & \none & \none & \none & \tl{2} & \none &   \none[\color{gray} \tl{4}]\\
\none & \none & \none & \none &\none & \none & \none & \none & \none & \none & \none & \none & \tl{1} & \none & \none & \none & \tl{3} & \none &   \none[\color{gray} \tl{5}]\\
\none & \none & \none & \none &\none & \none & \none & \none & \none & \none & \none & \none & \tl{2} & \none & \none & \none & \tl{4} & \none &   \none[\color{gray} \tl{6}]\\
\none & \none & \none & \none & \tl{1} & \none & \none & \none & \tl{1} & \none & \none & \none & \tl{3} & \none & \tl{1} & \none & \tl{5} & \none &   \none[\color{gray} \tl{7}]\\
\none & \none[\mathrel{\raisebox{-0.7ex}{$\scalebox{0.3}{\dots\dots}$}}] & \none & \none[\mathrel{\raisebox{-0.7ex}{$\scalebox{0.3}{\dots\dots}$}}] & \tl{2} & \none[\mathrel{\raisebox{-0.7ex}{$\scalebox{0.3}{\dots\dots}$}}] & \none & \none[\mathrel{\raisebox{-0.7ex}{$\scalebox{0.3}{\dots\dots}$}}] & \tl{2} & \none[\mathrel{\raisebox{-0.7ex}{$\scalebox{0.3}{\dots\dots}$}}] & \none & \none[\mathrel{\raisebox{-0.7ex}{$\scalebox{0.3}{\dots\dots}$}}] & \tl{4} & \none[\mathrel{\raisebox{-0.7ex}{$\scalebox{0.3}{\dots\dots}$}}] & \tl{2} & \none[\mathrel{\raisebox{-0.7ex}{$\scalebox{0.3}{\dots\dots}$}}] & \tl{6} & \none &   \none[\color{gray} \tl{8}]\\
\tl{1} & \none & \tl{1} & \none & \none & \none & \tl{1} & \none & \none & \none & \tl{1} & \none & \none & \none & \tl{3} & \none & \none  \\
\tl{3} & \none & \tl{3} & \none & \none & \none & \tl{5} & \none & \none & \none & \tl{7} & \none & \none & \none & \none &   \none \\
\tl{4} & \none & \tl{4} & \none & \none & \none & \tl{6} & \none & \none & \none & \none & \none & \none & \none &  \none \\
\tl{5} & \none & \none & \none & \none & \none & \none & \none & \none & \none & \none & \none &   \none \\
\none & \none & \none & \none & \none & \none & \none & \none & \none & \none & \none & \none &   \none \\
\end{ytableau} \\
\td{U}_8 \hskip 3.6mm \td{U}_7 \hskip 3.6mm \td{U}_6 \hskip 3.5mm \td{U}_5 \hskip 3.5mm
\td{U}_4 \hskip 3.5mm \td{U}_3 \hskip 3.5mm \td{U}_2 \hskip 3.5mm \td{U}_1 \hskip 3.5mm \td{U}_0 \hskip 3.5mm U\
\end{split}
\end{equation*} 
Then $\td{\bf U}=(\td{U}_7,\td{U}_6,\td{U}_5,\td{U}_4,\td{U}_3,\td{U}_2,\td{U}_1,\td{U}_0)$ corresponds to $\td{\bf T}\in {\bf H}(\td{\mu},8)$, with $\td{\mu}=(3,3,2,1)$. Repeating this process to $\td{\bf U}$ as in Example \ref{ex:3.13} (recall subsection \ref{subsec:separation}),
we get $\ov{\bf T}$
\begin{equation*}
\hskip 4cm\ov{\bf T}=\hskip -4cm
\begin{split}
\ytableausetup {mathmode, boxsize=1.0em} 
\begin{ytableau}
\none & \none & \none & \none & \none & \none & \none & \none &   \none \\
\none & \none & \none & \none & \none & \none & \none & \none & \tl{1}  \\
\none & \none & \none & \none & \none & \none & \none & \none & \tl{2}  \\
\none & \none & \none & \none & \none & \none & \none & \tl{1} & \tl{3}  \\
\none & \none & \none & \none & \none & \none & \none & \tl{2} & \tl{4}  \\
\none & \none & \none & \none & \tl{1} & \tl{1} & \tl{1} & \tl{3} & \tl{5} \\
\none & \none & \none & \none & \tl{2} & \tl{2} & \tl{2} & \tl{4} & \tl{6}\\
\tl{1} & \tl{1} & \tl{1} & \tl{1} & \tl{3} & \none & \none & \none & \none   \\
\tl{3} & \tl{3} & \tl{5} & \tl{7} & \none & \none & \none & \none & \none & \none   \\
\tl{4} & \tl{4} & \tl{6} & \none & \none & \none & \none & \none & \none & \none   \\
\tl{5} & \none & \none & \none & \none & \none & \none & \none & \none & \none & \none  \\
\end{ytableau}
\end{split}
\end{equation*}
Hence
\begin{equation*}
\hskip 2cm\ov{\bf T}^{\texttt{tail}}=\hskip-3cm
\begin{split}
\ytableausetup {mathmode, boxsize=1.0em} 
\begin{ytableau}
\none & \none & \none & \none & \none  \\
\none & \none & \none & \none & \none  \\
\none & \none & \none & \none & \none  \\
\none & \none & \none & \none & \none  \\
\none & \none & \none & \none & \none  \\
\none & \none & \none & \none & \none  \\
\tl{1} & \tl{1} & \tl{1} & \tl{1} & \tl{3}   \\
\tl{3} & \tl{3} & \tl{5} & \tl{7} & \none    \\
\tl{4} & \tl{4} & \tl{6} & \none & \none     \\
\tl{5} & \none & \none & \none & \none      \\
\end{ytableau}
\end{split}\quad\quad\quad\quad
\ov{\bf T}^{\texttt{body}}=\quad
\begin{split}
\ytableausetup {mathmode, boxsize=1.0em} 
\begin{ytableau}
  \none & \none & \none & \none & \tl{1}  \\
  \none & \none & \none & \none & \tl{2}  \\
  \none & \none & \none & \tl{1} & \tl{3}  \\
  \none & \none & \none & \tl{2} & \tl{4}  \\
 \tl{1} & \tl{1} & \tl{1} & \tl{3} & \tl{5} \\
 \tl{2} & \tl{2} & \tl{2} & \tl{4} & \tl{6}\\
 \none & \none & \none & \none & \none \\
 \none & \none & \none & \none & \none \\
 \none & \none & \none & \none & \none \\
 \none & \none & \none & \none & \none \\
\end{ytableau}
\end{split}
\end{equation*}
}
\end{ex}

%
\begin{rem} 
\mbox{}
{\rm 
\begin{itemize}
	\item[(1)] In \cite{JK20}, the authors define the separation algorithm (of type $D$) for an arbitrary element in ${\bf T}(\mu, n)$, generalizing the one (on the highest weight elements) in this paper and discuss its other application.
	\vskip 1mm
	
	\item[(2)] For types $B$ and $C$, the operator $\mathcal{S}_j$ is always $\mathcal{F}_j^{a_i}$ by Lemma \ref{cor:highest_weight_element}. Then, it is not difficult to check that the separation algorithm in this paper coincides with the one in \cite[Section 5.1]{K18-2} for types $B$ and $C$ (cf. Corollary \ref{cor:l-equivalence under S}).
\end{itemize}
}
\end{rem}

\section{Branching multiplicities}\label{sec:branching}

\subsection{Combinatorial description  of branching multiplicities} \label{subsec:combi_branching}

Suppose that $\la\in \cP_{n}$ and $\delta\in \cP^{(2)}_{n}$ are given. 
 Let 
\begin{equation}\label{eq:LR branching for d}
\texttt{LR}^{\mu}_{\lambda}(\mathfrak{d})
=\left\{\,{\bf T}\,\,\left|\,\,{\bf T}\in {\bf H}(\mu,n),\ {\bf T}\equiv_{\mf l} H_{\la'}\,\right.\right\}.
\end{equation}
Let $c^{\mu}_{\lambda}(\mathfrak{d}) = |\texttt{LR}^{\mu}_{\lambda}(\mathfrak{d})|$. 
Then $c^{\mu}_{\lambda}(\mathfrak{d})$ is equal to the multiplicity of irreducible highest weight $\mf l$-module with highest weight $\sum_{i \ge 1}\lambda_i^{'}\epsilon_i$ in the irreducible highest weight $\mf g$-module with highest weight $\La(\mu)$.

Let ${\bf \delta}^{\texttt{rev}} = (\delta_1^{\texttt{rev}},\dots,\delta_n^{\texttt{rev}})$ be the reverse sequence of $\delta=(\delta_1, \dots,\delta_n)$, that is, $\delta^{\texttt{rev}}_i=\delta_{n-i+1}$, for $1\leq i\leq n$. We put $p=\mu'_1$, $q=\mu'_2$, and $r=(\ov{\mu})'_1$ if $n-2\mu'_1<0$. 

\begin{df}\label{bounded orthogonal LR}
{\rm 
For $S \in \texttt{LR}^{\lambda'}_{\delta'\mu'}$, 
let $s_1\leq \dots \leq s_p$ denote the entries in the first row, and $t_1\leq \dots \leq t_q$ the entries in the second row of $S$.
Let $1 \le m_1 < \cdots < m_p < n$ be the sequence defined inductively from $p$ to $1$ as follows:
\begin{equation*}
m_i = \max\{\,k\,|\,\delta^{\texttt{rev}}_k\in X_i,\ \delta^{\texttt{rev}}_k<s_i\,\},
\end{equation*}
where 
\begin{equation*}
X_i =
\begin{cases}
\{\,\delta^{\texttt{rev}}_{i},\dots,\delta^{\texttt{rev}}_{2i-1}\,\}\setminus \{\delta^{\texttt{rev}}_{m_{i+1}},\dots,\delta^{\texttt{rev}}_{m_p}\}, & \text{if $1\leq i\leq r$}, \\
\{\,\delta^{\texttt{rev}}_{i},\dots,\delta^{\texttt{rev}}_{n-p+i}\,\}\setminus \{\delta^{\texttt{rev}}_{m_{i+1}},\dots,\delta^{\texttt{rev}}_{m_p}\}, & \text{if $r<i\leq p$},
\end{cases}
\end{equation*}
(we assume that $r=p$ when $n-2\mu'_1\geq 0$).
Let $n_1<\dots<n_q$ be the sequence such that $n_j$ is the $j$-th smallest integer in $\{j+1, \cdots, n\}\setminus \{ m_{j+1}, \cdots, m_p \}$ for $1 \le j \le q$.
Then we define $\ov{{\texttt {LR}}}^{\lambda'}_{\delta'\mu'}$ to be a subset of ${\texttt {LR}}^{\lambda'}_{\delta'\mu'}$ consisting of $S$ satisfying  
\begin{equation} \label{eq:condition_on_second_row}
t_j > \delta_{n_j}^{\texttt{rev}},
\end{equation}
for $1 \le j \le q$.
We put $\overline{c}^{\lambda}_{\delta\mu} = |\ov{\texttt{LR}}^{\lambda'}_{\delta'\mu'}|$.
}
\end{df}

\begin{rem} \label{rem:existence m_i}
{\rm
Let $S \in \texttt{LR}^{\lambda'}_{\delta'\mu'}$ be given. 
Let us briefly explain the well-definedness of the sequence $(m_i)_{1\le i\le p}$ in Definition \ref{bounded orthogonal LR}.
We may assume that $n - 2\mu'_1 \ge 0$ since the arguments for $n - 2\mu'_1 < 0$ are similar. 

It is enough to verify that $\delta_i^{\texttt{rev}} < s_i$ for $1 \le i \le p$.
Let 
$H'=(s_1 \rightarrow (s_2 \rightarrow \dots (s_p \rightarrow H_{\delta'})))$.
Then ${\rm sh}(H')/{\rm sh}(H_{\delta'})$ is a horizontal strip of length $p$.
If there exists $s_i$ such that $\delta_i^{\texttt{rev}} \ge s_i$, then we should have
$\ell(\la)>n$, which is a contradiction to $\la\in \cP_n$.
By definition of $m_i$, we also note that
\begin{equation*}
\begin{cases} 
			i \le m_i \le 2i-1 & \textrm{for $1 \le i \le r$}, \\
			i \le m_i \le n-p+ i & \textrm{for $r < i \le p$},
\end{cases}
\end{equation*}
where $r=p$ when $n-2\mu'_1\geq 0$.
}
\end{rem}

\begin{ex} \label{ex:example for LR var}
	{\rm
	Let $n=8$, $\mu = (2, 2, 2, 1, 1) \in \mathcal{P}({\rm O}_8)$, $\lambda = (5, 4, 4, 3, 2, 2) \in \cP_8$
	, and $\delta = (4,2,2,2,2) \in \cP_8^{(2)}$. Note that $n-2\mu_1' = -2 < 0$ and $r = (\ov{\mu})_1' = 3$.
	
	Let us consider the Littlewood-Richardson tableau $S \in \texttt{LR}^{\lambda'}_{\delta' \mu'}$ given by
	\begin{equation*}
		\hskip 5cm S=\hskip-5cm
		\begin{split}
		\ytableausetup {mathmode, boxsize=1.0em} 
			\begin{ytableau}
				\none & \none & \none & \none  & \none \\
				\tl{1} & \tl{3} & \tl{3} & \tl{3} & \tl{5}   \\
				\tl{2} & \tl{4} & \tl{4} & \none & \none    \\
				\none & \none & \none & \none  & \none \\
			\end{ytableau}
		\end{split} \quad . \quad\quad\quad
	\end{equation*}
	Then the sequences $(m_i)_{1 \le i \le 5}$ and $(n_j)_{1 \le j \le 3}$ are $(1, 3, 5, 7, 8)$ 
	and $(2, 4, 6)$, respectively, and $S$ satisfies the condition \eqref{eq:condition_on_second_row}:
	\begin{equation*}
		t_1 = 2 > 0 = \delta_{n_1}^{\texttt{rev}}, \quad t_2 = 4 > 2 = \delta_{n_2}^{\texttt{rev}}, \quad
		t_3 = 4 > 2 = \delta_{n_3}^{\texttt{rev}}.
	\end{equation*}
Hence $S \in \ov{{\texttt {LR}}}^{\lambda'}_{\delta'\mu'}$.
	}
\end{ex}

Now we are in a position to state the main result in this paper.
The proof is given in Section \ref{sec:proof of main}.

\begin{thm} \label{thm:main1} 
For $\mu \in \mathcal{P}({\rm O}_n)$ and $\la\in \cP_n$,
we have a bijection
\begin{equation*}
\xymatrixcolsep{3pc}\xymatrixrowsep{0pc}\xymatrix{
\texttt{{\em LR}}^{\mu}_{\lambda}(\mathfrak{d}) \ar@{->}[r] & \ \bigsqcup_{\delta \in \cP_{n}^{(2)}} \ov{\texttt{{\em LR}}}^{\lambda'}_{\delta'\mu'}  \\ 
\mathbf{T}   \ar@{|->}[r] & \ov{\bf T}^{\texttt{\em tail}} 
}.
\end{equation*}		
\end{thm}

\vskip 3mm

\begin{cor}\label{cor:main result}
Under the above hypothesis, we have
\begin{equation*}
c^{\mu}_{\lambda}(\mathfrak{d}) 
= \sum_{\delta \in \cP^{(2)}_n}\overline{c}^{\lambda}_{\delta\mu}.
\end{equation*}
\end{cor}

Let us give another description of $c^{\mu}_{\la}(\mf{d})$ which is simpler than $\ov{\texttt{LR}}^{\la'}_{\delta'\mu'}$, and also plays an important role in  Section \ref{sec:genexp}.

\begin{df}\label{flagged orthogonal LR+}{\rm 
For $U\in {\texttt {LR}}^{\la}_{\delta\mu^\pi}$ (see subsection \ref{subsec:notations}), let 
$\sigma_1>\dots>\sigma_p$ denote the entries in the rightmost column and $\tau_1>\dots>\tau_q$ the  second rightmost column of $U$, respectively.
Let $m_1<\dots<m_p$ be the sequence defined by
\begin{equation*}
m_i=
\begin{cases}
\min\{n-\sigma_i+1,2i-1\}, & \text{if $1\leq i\leq r$},\\
\min\{n-\sigma_i+1,n-p+i\}, & \text{if $r< i\leq p$}.
\end{cases}
\end{equation*} 
and let $n_1<\dots<n_q$ be the sequence such that $n_j$ is the $j$-th smallest number in $\{\,j+1,\dots,n\,\}\setminus\{\,m_{j+1},\dots,m_{p}\,\}$.
Then we define $\underline{\texttt {LR}}^{\la}_{\delta\mu}$ to be the subset of ${\texttt {LR}}^{\la}_{\delta\mu^\pi}$ consisting of $U$ such that
\begin{equation}\label{eq:flagged conditions}
\tau_j + n_j\le n+1,
\end{equation}
for $1\leq j\leq q$.
We put 
$\underline{c}^\la_{\delta\mu}=|\underline{\texttt {LR}}^{\la}_{\delta\mu}|$.
}
\end{df}

\begin{ex} \label{ex:example for LR undervar}
	{\rm
		We keep the assumption in Example \ref{ex:example for LR var} and consider the Littlewood-Richardson tableau $U \in \texttt{LR}^{\lambda}_{\delta \mu^{\pi}}$ given by
		\begin{equation*}
			\hskip 5cm U = \hskip -6cm
			\begin{split}
				\ytableausetup {mathmode, boxsize=1.0em} 
				\begin{ytableau}
					\none & \none \\
					 \none & \tl{1}  \\
					 \none & \tl{2}  \\
					 \tl{2} & \tl{3}  \\
					 \tl{3} & \tl{4} \\
					\tl{6} & \tl{6} \\
					\none & \none \\
				\end{ytableau} 
			\end{split}\quad .
		\end{equation*}
		The sequences $(m_i)_{1 \le i \le 5}$ and $(n_j)_{1 \le j \le 3}$ are $(1, 3, 5, 7, 8)$ 
	and $(2, 4, 6)$, respectively. Then $U$ satisfies the condition \eqref{eq:flagged conditions}:
		\begin{equation*}
			\begin{split}
				\tau_1 + n_1 = 6 + 2 = 8 \le 8+1 = n+1, \\
				\tau_2 + n_2 = 3 + 4 = 7 \le 8+1 = n+1, \\
				\tau_3 + n_3 = 2 + 6 = 8 \le 8+1 = n+1.
			\end{split}
		\end{equation*} 
		Hence $U \in \underline{\texttt {LR}}^{\la}_{\mu \delta}$.
	}
	\end{ex}

Now, one can show that $\ov{c}^{\lambda}_{\delta\mu}= \underline{c}^\la_{\delta\mu}$
by using the bijection $\psi$ \eqref{conjugation of LR}.

\begin{lem} \label{lem:equivalence}
The sequences $(m_i)_{1 \le i \le p}$ and $(n_j)_{1 \le j \le q}$ for $S$ in Definition \ref{bounded orthogonal LR} are equal to the ones for $U=\psi(S)$ in Definition \ref{flagged orthogonal LR+}.
\end{lem}
\pf
We assume that $n - 2\mu_1' \ge 0$. The proof for the case $n - 2\mu_1' < 0$ is similar.
Suppose that $S\in \texttt {LR}^{\la'}_{\delta'\mu'}$ is given. Let $s_1\leq \dots \leq s_p$ denote the entries in the first row of $S$.
Let $(m'_i)_{1 \le i \le p}$ and $(n'_j)_{1 \le j \le q}$ be the sequences for $S$ in Definition \ref{bounded orthogonal LR}.
Put $U=\psi(S)$. Let $\sigma_1>\dots>\sigma_p$ be the rightmost column of $U$ and let $(m_i)_{1 \le i \le p}$ and $(n_j)_{1 \le j \le q}$ be the sequences for $U$ as in Definition \ref{flagged orthogonal LR+}.

It is enough to show that $m_i' = m_i$ for $1 \le i \le p$, which clearly implies $n_j' = n_j$ for $1 \le j \le q$.
Let us enumerate the column of $\delta'$ by $n, n-1,\dots,1$ from left to right.
Consider the vertical strip $V^i:={\rm sh}(H^i)/{\rm sh}(H^{i-1})$ filled with $i$ for $1\leq i\leq p$ (recall \eqref{conjugation of LR} below). 
By definition of $\psi$ \eqref{conjugation of LR}, we see that the upper most box in $V^i$ is located in the $(n-\sigma_i+1)$-th column in $\delta'$.

Let $i\in \{1,\dots,p\}$ be given. 
First, we have $ m_i' \le n-\sigma_i + 1$ by definition of $m_i'$. Since $m_i' \le 2i-1$, we have $m_i' \le m_i=\min\{n-\sigma_i+1,2i-1\}$.
Next, we claim that $m_i\le m'_i$. If $n-\sigma_i+1 \le 2i-1$, then 
we have $\delta_{n-\sigma_i+1}^{\texttt{rev}} < s_i$, and hence $m_i\le n-\sigma_i+1 \le m_i'$ by definition of $m_i'$. 
If $n-\sigma_i+1 > 2i-1$, then we have $m_i = 2i-1 = m_i'$. This proves that $m_i = m'_i$.
\qed
\begin{thm}\label{thm:main result-flagged}
For $\mu\in {\mc P}({\rm O}_n)$, $\la\in \cP_{n}$ and $\delta\in \cP^{(2)}_{n}$, the bijection $\psi : {\texttt {\em LR}}^{\la'}_{\mu'\nu'}  \longrightarrow {\texttt {\em LR}}^{\la}_{\mu\nu^\pi}$ in \eqref{conjugation of LR} induces a bijection from $\overline{\texttt {\em LR}}^{\la'}_{\delta'\mu'}$ to $\underline{\texttt{\em LR}}^{\lambda}_{\delta\mu}$. 
\end{thm}
\pf Let $S \in {\texttt {LR}}^{\la'}_{\delta'\mu'}$ given and put $U=\psi(S)$. We keep the conventions in the proof of Lemma \ref{lem:equivalence}.
By definition of $\psi$, the second upper most box in \red{$V^j$} is located at the $(n-\tau_j+1)$-th column in $\delta'$. 
By Lemma \ref{lem:equivalence}, we see that $\delta_{n_j}^{\texttt{rev}} < t_j$ if and only if
$n-\tau_j + 1 \ge n_j$ or $\tau_j + n_j \le n+1$.
Therefore, $S\in \overline{\texttt {LR}}^{\la'}_{\delta'\mu'}$ if and only if $U\in \underline{\texttt{LR}}^{\lambda}_{\delta\mu}$.
\qed

\begin{ex} \label{ex:bijection}
{\rm 
Let $n=8$, $\mu = (2, 2, 2, 1, 1) \in \mathcal{P}({\rm O}_8)$, $\lambda = (5, 4, 4, 3, 2, 2) \in \cP_8$, and $\delta = (4,2,2,2,2) \in \cP_8^{(2)}$. Let $S$ be the Littlewood-Richardson tableau in Example \ref{ex:example for LR var}. 
We enumerate the columns of $S$ as follows: 
\begin{equation*}
	\begin{split}
	\ytableausetup {mathmode, boxsize=1.0em} 
		\begin{ytableau}
			\tl{\red{1}} & \none &\tl{\red{3}} &\none & \tl{\red{3}} &\none & \tl{\red{3}} &\none & \tl{\red{5}}   \\
			\tl{\blue{2}} & \none &\tl{\blue{4}} &\none & \tl{\blue{4}} & \none & \none &\none & \none    \\
			\none & \none &\none &\none & \none & \none &\none  &\none & \none \\
			\none[\scalebox{0.6}{$S^5$}] & \none &\none[\scalebox{0.6}{$S^4$}] & \none & \none[\scalebox{0.6}{$S^3$}] & \none &\none[\scalebox{0.6}{$S^2$}] &\none & \none[\scalebox{0.6}{$S^1$}] \\
		\end{ytableau}
	\end{split}
\end{equation*}
Then the insertion and recording tableaux are given by
\begin{equation*}
	\hskip 1cm S \rightarrow H_{\delta'} \ = \hskip -1.5cm 
	\begin{split}
	\ytableausetup {mathmode, boxsize=1.0em}  
	\begin{ytableau}
		\none & \none[\scalebox{0.5}{\color{gray}{$8$}}] & \none[\scalebox{0.5}{\color{gray}{$7$}}]  & \none[\scalebox{0.5}{\color{gray}{$6$}}] & \none[\scalebox{0.5}{\color{gray}{$5$}}] & \none[\scalebox{0.5}{\color{gray}{$4$}}] & \none[\scalebox{0.5}{\color{gray}{$3$}}] & \none[\scalebox{0.5}{\color{gray}{$2$}}] & \none[\scalebox{0.5}{\color{gray}{$1$}}] \\
		\none & \none[\scalebox{0.75}{\color{gray}{$\vdots$}}] & \none[\scalebox{0.75}{\color{gray}{$\vdots$}}]  & \none[\scalebox{0.75}{\color{gray}{$\vdots$}}] & \none[\scalebox{0.75}{\color{gray}{$\vdots$}}] & \none[\scalebox{0.75}{\color{gray}{$\vdots$}}] & \none[\scalebox{0.75}{\color{gray}{$\vdots$}}] & \none[\scalebox{0.75}{\color{gray}{$\vdots$}}] & \none[\scalebox{0.75}{\color{gray}{$\vdots$}}] \\
		\none & \none[\tl{1}] & \none[\tl{1}] & \none[\tl{1}] & \none[\tl{1}]  & \none[\tl{1}] & \red{\tl{1}} & \none \\
		\none & \none[\tl{2}] & \none[\tl{2}]  & \none[\tl{2}]  & \none[\tl{2}] & \none[\tl{2}] & \blue{\tl{2}} & \none \\
		\none & \none[\tl{3}] & \red{\tl{3}}  & \red{\tl{3}} & \red{\tl{3}} & \none & \none  & \none  \\
		\none & \none[\tl{4}] & \blue{\tl{4}}  & \blue{\tl{4}} & \none & \none & \none & \none\\
		\none & \red{\tl{5}} & \none  & \none & \none & \none & \none & \none\\
		\none & \none & \none  & \none & \none & \none \\
	\end{ytableau}
	\end{split} \quad , \ \
	\hskip 0.5cm Q(S \rightarrow H_{\delta'}) \ = 
	\begin{split}
		\ytableausetup {mathmode, boxsize=1.0em} 
		\begin{ytableau}
			\none & \none & \none & \none & \none  & \none  \\
			\none & \none & \none & \none & \none  & \none  \\
			\none & \none & \none & \none & \none  & \none  \\
			\none & \none & \none & \none & \none  & \tl{\red{1}}  \\
			\none & \none & \none  & \tl{\red{2}}  & \tl{\blue{3}} & \none  \\
			\none & \none & \none  & \tl{\red{3}} & \tl{\blue{4}} & \none    \\
			\none & \none & \none  & \tl{\red{4}} & \none & \none \\
			\none & \none & \none  & \none & \none & \none \\
			\none & \tl{\red{5}} & \tl{\blue{5}}  & \none & \none & \none \\
			\none & \none & \none  & \none & \none & \none \\
		\end{ytableau}
	\end{split}\quad . \hskip 1cm
\end{equation*}
Then $\psi(S)$ is obtained by
\begin{equation*}
	\hskip 2.5cm \left(
	\begin{split}
		\ytableausetup {mathmode, boxsize=1.0em} 
		\begin{ytableau}
			\none & \none & \none & \none & \none  & \none  \\
			\none[\scalebox{0.6}{\color{gray}{$\tl{8}\cdots$}}] & \none[\color{black}{\tl{1}}] & \none[\color{black}{\tl{1}}] & \none[\color{black}{\tl{1}}] & \none[\color{black}{\tl{1}}] & \tl{\red{1}}  \\
			\none[\scalebox{0.6}{\color{gray}{$\tl{7}\cdots$}}] & \none[\color{black}{\tl{2}}] & \none[\color{black}{\tl{2}}]  & \tl{\red{2}}  & \tl{\blue{2}} & \none  \\
			\none[\scalebox{0.6}{\color{gray}{$\tl{6}\cdots$}}] & \none[\color{black}{\tl{3}}] & \none[\color{black}{\tl{3}}]  & \tl{\red{3}} & \tl{\blue{3}} & \none    \\
			\none[\scalebox{0.6}{\color{gray}{$\tl{5}\cdots$}}] & \none[\color{black}{\tl{4}}] & \none[\color{black}{\tl{4}}]  & \tl{\red{4}} & \none & \none \\
			\none[\scalebox{0.6}{\color{gray}{$\tl{4}\cdots$}}] & \none[\color{black}{\tl{5}}] & \none[\color{black}{\tl{5}}] & \none  & \none & \none \\
			\none[\scalebox{0.6}{\color{gray}{$\tl{3}\cdots$}}] & \tl{\red{6}} & \tl{\blue{6}} & \none & \none  & \none \\
			\none[\scalebox{0.6}{\color{gray}{$\tl{2}\cdots$}}] & \none & \none & \none & \none  & \none  \\
			\none[\scalebox{0.6}{\color{gray}{$\tl{1}\cdots$}}] & \none & \none & \none & \none  & \none  \\
			\none & \none & \none[\quad \scalebox{0.7}{$\psi(S) \rightarrow H_{\delta}$}] & \none & \none  & \none  \\
		\end{ytableau}
	\end{split}   \quad , 
	\begin{split}
		\ytableausetup {mathmode, boxsize=1.0em} 
		\begin{ytableau}
			\none & \none & \none & \none & \none &  \none  \\
			\none & \none & \none & \none & \none &  \tl{\red{1}}  \\
			\none & \none & \none &  \tl{\red{2}}  & \tl{\blue{3}} & \none  \\
			\none & \none & \none &  \tl{\red{3}} & \tl{\blue{4}} & \none    \\
			\none & \none & \none &  \tl{\red{4}} & \none & \none \\
			\none & \none & \none &  \none & \none & \none \\
			\none & \tl{\red{5}} & \tl{\blue{5}} & \none  & \none & \none \\
			\none & \none & \none & \none & \none  & \none  \\
			\none & \none & \none & \none & \none  & \none  \\
			\none & \none & \none[\quad \scalebox{0.7}{$Q(\psi(S) \rightarrow H_{\delta})$}] & \none & \none  & \none  \\
		\end{ytableau} 
	\end{split} 
	\right)
	\ \xymatrixcolsep{3pc}\xymatrixrowsep{0pc}\xymatrix{
		\ar@{->}[r] & } \ 
	\begin{split}
		\ytableausetup {mathmode, boxsize=1.0em} 
			\begin{ytableau}
				\none & \none \\
				 \none & \tl{\red{1}}  \\
				 \none & \tl{\red{2}}  \\
				 \tl{\blue{2}} & \tl{\red{3}}  \\
				 \tl{\blue{3}} & \tl{\red{4}} \\
				\tl{\blue{6}} & \tl{\red{6}} \\
				\none & \none \\
			\end{ytableau} 
	\end{split} = \psi(S). \hskip 13.5cm
\end{equation*}
(Here the numbers in gray denote the enumeration of columns of $\delta'$.)
Thus we have
\begin{equation*}
	(\sigma_1, \sigma_2, \sigma_3, \sigma_4, \sigma_5) = (6, 4, 3, 2, 1), \quad \ \ (\tau_1, \tau_2, \tau_3) = (6, 3, 2). \quad \ \ 
\end{equation*}
Note that the enumeration of the rows of $U := \psi(S)$ is given by 
\begin{equation*}
	\begin{split}
	\ytableausetup {mathmode, boxsize=1.0em} 
	\begin{ytableau}
		\none[\scalebox{0.6}{$U_1$}] &\none &\none & \tl{\red{1}}&\none&\none  \\
		\none &\none &\none & \none&\none&\none \\
		\none[\scalebox{0.6}{$U_2$}] &\none &\none & \tl{\red{2}}&\none&\none  \\
		\none &\none &\none & \none&\none&\none \\
		\none[\scalebox{0.6}{$U_3$}] &\none &\tl{\blue{2}} & \tl{\red{3}}&\none&\none  \\
		\none &\none &\none & \none&\none&\none \\
		\none[\scalebox{0.6}{$U_4$}] &\none &\tl{\blue{3}} & \tl{\red{4}}&\none&\none \\
		\none &\none &\none & \none&\none&\none\\
		\none[\scalebox{0.6}{$U_5$}] &\none &\tl{\blue{6}} & \tl{\red{6}}&\none&\none \\
	\end{ytableau} 
	\end{split}.
\end{equation*}

Under the above correspondence, 
we observe that $\sigma_i \ (1 \le i \le 5)$ and $\tau_j \ (1 \le j \le 3)$ record
the positions of $s_i$ and $t_j$ in $\delta'$, respectively, and vice versa.
For example, the entry $\tau_2$ in $U^4$ is located at $(n-\tau_2+1)$-th row in $\delta$, 
which implies that $t_2$ is located at $(n-\tau_2+1)$-th column in $\delta'$.
This correspondence implies $\psi(S) \in \underline{\texttt{LR}}^{\lambda}_{\delta\mu}$ (cf. Example \ref{ex:example for LR undervar}).
} 
\end{ex}

\begin{cor}\label{cor:main result-2}
Under the above hypothesis, we have
\begin{equation*}
c^{\mu}_{\lambda}(\mathfrak{d}) 
= \sum_{\delta \in \cP^{(2)}_n}\underline{c}^{\lambda}_{\delta\mu}.
\end{equation*}
\end{cor}
\pf This follows from $\ov{c}^{\lambda}_{\delta\mu}= \underline{c}^\la_{\delta\mu}$.
\qed \vskip 2mm

We have another characterization of $\underline{c}^\la_{\delta\mu}$ in terms of usual LR tableaux (not companion tableaux) by considering the bijection between LR tableaux and their companion ones.

\begin{cor} \label{cor:LR version}
Let $U$ be an LR tableau of shape $\la/\delta$ with content $\mu^\pi$
and let $\sigma_i$ be the row index of the leftmost $\mu'_1-i+1$ in $U$ for $1 \le i \le \mu_1'$, and $\tau_j$ the row index of the second leftmost $\mu'_2-j+1$ in $U$ for $1 \le j \le \mu_2'$.
Let $m_1<\dots<m_{\mu'_1}$ be the sequence given by
$m_i=\min\{n-\sigma_i+1,2i-1\}$, 
and let $n_1\le\dots\le n_{\mu'_2}$ be the sequence such that 
$n_j$ is the $j$-th smallest number in 
$\{\,j+1,\dots,n\,\}\setminus\{\,m_{j+1},\dots,m_{\mu'_1}\,\}$.
Then,
$\underline{c}^{\la}_{\delta\mu}$ is equal to the number of $U$ such that
\begin{equation*}
\tau_j + n_j \le n+1,
\end{equation*}
for $1\leq j\leq \mu'_2$.
\end{cor}

We may recover the Littlewood's formula \eqref{eq:Littlewood restriction} from Corollary \ref{cor:main result-2}.
\begin{cor} \label{cor:littlewood}
Under the above hypothesis, if $\ell(\lambda) \le \frac{n}{2}$, then
\begin{equation*}
c^{\mu}_{\lambda}(\mathfrak{d}) 
= \sum_{\delta \in \cP^{(2)}_n}c^{\lambda}_{\delta\mu}.
\end{equation*}
\end{cor}
\pf We claim that $\texttt{LR}^{\lambda}_{\delta\mu^{\pi}}=\underline{\texttt{LR}}^{\lambda}_{\delta\mu}$.
Let $U \in \texttt{LR}^{\lambda}_{\delta\mu^{\pi}}$ be given.
Let $H'= (\sigma_1\rightarrow (\dots \rightarrow ( \sigma_p \rightarrow H_\delta)))$. 
Note that $\sigma_i+i-1\le \ell({\rm sh}(H'))=\ell(\la)\leq \frac{n}{2}$ for $1\le i\le p$.
So we have
\begin{equation}\label{eq:sigma bound}
n-\sigma_i+1 \ge 2i\quad (1 \le i \le p).
\end{equation}
Otherwise, we have
$n-i < \sigma_i + i -1 \le \frac{n}{2}$ and hence $\frac{n}{2} < i \le p=\mu'_1\le \ell(\la)$, which is a contradiction.
By definition $m_i$ and $n_j$, we have 
\begin{equation}\label{eq:m for stable}
	m_i = 2i-1, \quad n_j = 2j,
\end{equation} 
for $1 \le i \le p$ and $1 \le j \le q$. 
By \eqref{eq:sigma bound} and \eqref{eq:m for stable} we have
\begin{equation*}
\tau_j\leq \sigma_j\leq n-2j+1\quad (1\leq j\leq q),
\end{equation*}
which implies that $U$ satisfies \eqref{eq:flagged conditions}, that is, $U \in \underline{\texttt{LR}}^{\lambda}_{\delta\mu}$. This proves the claim.
By Theorem \ref{thm:main result-flagged}, we have ${c}^{\lambda}_{\delta\mu}= \underline{c}^\la_{\delta\mu}$.
\qed

\begin{rem}{\rm
For ${\bf T}=(T_l,\dots,T_0)\in \texttt{{LR}}^{\mu}_{\lambda}(\mathfrak{d})$, let ${\bf T}^{\texttt{tail}}=(T^{\texttt{tail}}_l,\dots ,T^{\texttt{tail}}_0)$.
We may regard ${\bf T}^{\texttt{tail}}$ as a column-semistandard tableau of shape $\mu'$ by putting together $T^{\texttt{tail}}_i$'s horizontally.
It is shown in \cite[Theorem 4.8]{K18-3} that if $\ell(\la)\leq n/2$, then the map sending ${\bf T}$ to ${\bf T}^{\texttt{tail}}$ gives a bijection
\begin{equation*}
\xymatrixcolsep{3pc}\xymatrixrowsep{0pc}\xymatrix{
\texttt{{LR}}^{\mu}_{\lambda}(\mathfrak{d}) \ar@{->}[r] & \ \bigsqcup_{\delta \in \cP_{n}^{(2)}} {\texttt{{LR}}}^{\lambda'}_{\delta'\mu'}
}.
\end{equation*}
By Lemma \ref{lem:equivalence} and \eqref{eq:m for stable}, we have $\ov{\bf T}^{\texttt{tail}}={\bf T}^{\texttt{tail}}$ if $\ell(\la)\leq n/2$, and hence Theorem \ref{thm:main1} recovers \cite[Theorem 4.8]{K18-3}.
}
\end{rem}

\begin{rem}{\rm
Let us briefly recall Sundaram's formula for \eqref{eq:branching mult} when ${\rm G}_n = {\rm Sp}_n$ \cite{Su}.
She constructs a bijection between the set of oscillating tableaux appearing in Berele's correspondence for ${\rm Sp}_n$ \cite{Be} and the set of pairs of the standard tableaux and LR tableaux with the symplectically fitting lattice word. Then it is shown that these LR tableaux count the branching multiplicity \eqref{eq:branching mult}. We remark that 
Lecouvey-Lenart provide a conjectural bijection between the Sundaram's LR tableaux and the flagged LR tableaux for type $C_n$ in \cite{LL}.

Recently, an orthogonal analogue of the above bijection \cite{Su} is given for ${\rm SO}_{2n+1}$ \cite{J}, where oscillating tableaux are replaced by vacillating tableaux, and the Sundaram's LR tableaux are replaced by so-called alternative orthogonal LR tableaux, which are in (highly non-trivial) bijection with $\texttt{{LR}}^{\mu}_{\lambda}(\mathfrak{d})$.

However, as far as we understand, alternative orthogonal LR tableaux do not seem to yield a formula similar to Sundaram's one, and recover Littlewood's restriction formula directly.
It would be very interesting to clarify this point.
}
\end{rem}

\begin{rem} \label{rem:for types B and C}
{\rm  
We may have an analogue of Theorem \ref{thm:main1} for types $B$ and $C$, that is, a multiplicity formula with respect to the branching from $B_\infty$ and $C_{\infty}$ to $A_{+\infty}$, respectively.

More precisely, let ${\bf T}^{\mf g}(\mu, n)$ be the spinor model for the integrable highest weight module over the Kac-Moody algebra of types $B_\infty$ and $C_\infty$ when $\mf g= \mf b$ and $\mf c$, respectively, corresponding to $\mu\in\cP({\rm G}_n)$ via Howe duality. 
Here $\cP({\rm G}_n)$ denotes the set of partitions parametrizing the finite-dimensional irreducible representations of an algebraic group ${\rm G}_n$ (see \cite[Section 2]{K18-3} for more details). 
 
For $\la\in \cP_n$, let $\texttt{LR}^{\mu}_{\lambda}(\mf g)$ be the set of ${\bf T}\in {\bf T}^{\mf g}(\mu, n)$ which is an $\mf l$-highest weight element with highest weight $\la'$ (cf. \eqref{eq:LR branching for d}). 
Let $\delta\in \cP_n^{\Diamond}$ be given, where $\Diamond = (1)$ for $\mf g=\mf b$ and $\Diamond=(1,1)$ for $\mf g=\mf c$ (here we understand $\cP^{(1)}=\cP$). 
Put
\begin{equation*}
\begin{split}
\ov{\texttt{LR}}^{\lambda'}_{\delta' \mu'} & 
= \left\{ S \in \texttt{LR}^{\lambda'}_{\delta' \mu'}  \,\Big|\, s_i > \delta^{\texttt{rev}}_{2i}\   (1 \le i \le \mu'_1) \right\},
\end{split}
\end{equation*} 
where $s_1\leq \dots \leq s_{\mu'_1}$ are the entries in the first row of $S$.

We may apply the same arguments in Section \ref{sec:separation} to ${\bf T}^{\mf g}(\mu, n)$. 
Then by Proposition \ref{prop:highest_weight_vectors} and Lemma \ref{lem:ov{T} is well-defined}, we have for ${\bf T} = (T_l, \dots, T_1) \in \texttt{LR}^{\mu}_{\lambda}(\mf g)$ that  
$${\bf T}^{\texttt{tail}}= (T_l^{\texttt{tail}}, \dots, T_0^{\texttt{tail}}) \in  \ov{\texttt{LR}}^{\lambda'}_{\delta'\mu'},$$ for some $\delta\in \cP_n^\Diamond$. 
Furthermore, the map
\begin{equation*}
\xymatrixcolsep{3pc}\xymatrixrowsep{0pc}\xymatrix{
\texttt{{LR}}^{\mu}_{\lambda}(\mf g) \ar@{->}[r] & \ \bigsqcup_{\delta \in \cP_{n}^\Diamond} \ov{\texttt{{LR}}}^{\lambda'}_{\delta'\mu'}  \\ 
\mathbf{T}   \ar@{|->}[r] & {\bf T}^{\texttt{tail}} 
}
\end{equation*}	 
is a bijection.  
The map $\psi$ in \eqref{conjugation of LR} induces a bijection from $ \ov{\texttt{LR}}^{\lambda'}_{\delta' \mu'}$ to $\underline{\texttt{LR}}^{\lambda}_{\delta \mu}$, where 
\begin{equation} \label{eq:flag condition for types B and C}
\begin{split}
\underline{\texttt{LR}}^{\lambda}_{\delta \mu} & 
= \left\{ U \in \texttt{LR}^{\lambda}_{\delta \mu^\pi}  \,\Big|\, \sigma_i+2i \leq n+1\ (1 \le i \le \mu'_1) \right\},
\end{split}
\end{equation} 
where $\sigma_1 > \dots > \sigma_{\mu'_1}$ are the entries in the rightmost column of $U$.
Therefore, 
\begin{equation*}
	c^{\mu}_{\lambda}(\mf g) = \sum_{\delta \in \cP_n^\Diamond}\ov{c}^{\lambda}_{\delta\mu}=\sum_{\delta \in \cP_n^\Diamond}\underline{c}^{\lambda}_{\delta\mu},
\end{equation*} 
where $c^{\mu}_{\lambda}(\mf g) = |\texttt{LR}^{\mu}_{\lambda}(\mf g)|$, 
$\ov{c}^{\lambda}_{\delta\mu} = |\ov{\texttt{LR}}^{\lambda'}_{\delta' \mu'}|$, and 
$\underline{c}^{\lambda}_{\delta\mu} = |\underline{\texttt{LR}}^{\lambda}_{\delta\mu}|$. 
This is a generalization of \cite[Theorem 4.8]{K18-3} for types $B$ and $C$, which also recovers \cite[Theorem 6.8]{LL} for type $C$.
}
\end{rem}

\subsection{Branching multiplicities of non-Levi type}  \label{sec:duality}
We assume that the base field is $\C$.
Let $V^\la_{{\rm GL}_n}$ denote the finite-dimensional irreducible ${\rm GL}_n$-module corresponding to $\la\in \cP_n$, and $V_{{\rm O}_n}^\mu$ the finite-dimensional irreducible module ${\rm O}_n$-module corresponding to $\mu\in \mc{P}({\rm O}_n)$. 

Then we have the following new combinatorial description of $\left[ V_{\rm{GL}_n}^{\la} : V_{{\rm O}_n}^{\mu} \right]$.

\begin{thm}\label{thm:non-levi branching}
For $\la\in \cP_{n}$ and $\mu\in {\mc P}({\rm O}_n)$, we have
\begin{equation*}
\left[ V_{\rm GL_n}^{\lambda} : V_{\rm O_n}^{\mu} \right] 
= \sum_{\delta \in \cP^{(2)}_n}\ov{c}^{\lambda}_{\delta\mu}
= \sum_{\delta \in \cP^{(2)}_n}\underline{c}^{\lambda}_{\delta\mu}.
\end{equation*}
\end{thm}
\pf It follows from the fact 
$\left[ V_{\rm GL_n}^{\lambda} : V_{\rm O_n}^{\mu} \right] = c_{\lambda}^{\mu}(\mathfrak{d})$ \cite[Theorem 5.3]{K18-3}, and  
Corollaries \ref{cor:main result} and \ref{cor:main result-2}.
\qed

\begin{ex} \label{ex:4.17}
{\rm %
Let us compare the formula in Theorem \ref{thm:non-levi branching} with the one by Enright and Willenbring in \cite[Theorem 4]{EW}. 

Let $\mu, \nu \in {\mc P}({\rm O}_n)$ be given by
\begin{equation*}
\begin{split}
	\mu & = (\LaTeXunderbrace{d, \LaTeXoverbrace{2, \dots, 2}^{a}, \LaTeXoverbrace{1, \dots, 1}^{b}, \LaTeXoverbrace{0, \dots, 0}^{c}}_{n}), \\
	\nu & = (\LaTeXunderbrace{d, \LaTeXoverbrace{2, \dots, 2}^{c}, \LaTeXoverbrace{1, \dots, 1}^{b}, \LaTeXoverbrace{0, \dots, 0}^{a}}_{n}).
\end{split}
\end{equation*} where $a, b, c, d$ are positive integers with $d \ge 2$. 
Then we have for $\la\in \cP_n$
\begin{equation*}
\left[ V_{\rm GL_n}^{\lambda} : V_{\rm O_n}^{\mu} \right]
= \sum_{\xi \in \cP_n^{(2)}} c^{\lambda'}_{\xi' \mu'}
- \sum_{\upsilon \in \cP_n^{(2)}} c^{\lambda'}_{\upsilon' \nu'},
\end{equation*}
(see \cite[Section 7 (7.11)]{EW}).
Suppose that $n=8$, $a = b = d= 2, c = 3$ and $\lambda = (5, 4, 4, 3, 2, 2, 0, 0) \in \cP_8$. Then it is straightforward to check that for $\xi, \upsilon \in \cP_8^{(2)}$
\begin{equation*}
\begin{split}
	c^{\lambda'}_{\xi' \mu'} & = \begin{cases}
						1, & \textrm{if $\xi = (4, 2, 2, 2, 2)$ or $(4, 4, 2, 2)$}, \\
						0, & \textrm{otherwise},
					\end{cases} \\
	c^{\lambda'}_{\upsilon' \nu'} & = \begin{cases}
						1, &  \textrm{if $\upsilon = (4, 2, 2, 2)$}, \\
						0, & \textrm{otherwise}.
					\end{cases}
\end{split}
\end{equation*} 
Hence we have
\begin{equation*}
\left[ V_{\rm GL_8}^{\lambda} : V_{\rm O_8}^{\mu} \right]
=
\sum_{\xi \in \cP_8^{(2)}} c^{\lambda'}_{\xi' \mu'}
- \sum_{\upsilon \in \cP_8^{(2)}} c^{\lambda'}_{\upsilon' \nu'}
= 2 - 1 = 1.
\end{equation*}
On the other hand, the following tableaux $S_{\alpha}$ and $S_{\beta}$ are the unique tableaux in $\texttt{LR}^{\lambda'}_{\alpha' \mu'}$ and $\texttt{LR}^{\lambda'}_{\beta' \mu'}$, respectively, where $\alpha = (4, 2, 2, 2, 2)$ and $\beta = (4, 4, 2, 2)$:
\begin{equation*}
\hskip 3cm S_{\alpha}=\hskip-3.5cm
\begin{split}
\ytableausetup {mathmode, boxsize=1.0em} 
\begin{ytableau}
\none & \none & \none & \none & \none  \\
\tl{1} & \tl{3} & \tl{3} & \tl{3} & \tl{5}   \\
\tl{2} & \tl{4} & \tl{4} & \none & \none    \\
\none & \none & \none & \none & \none  \\
\end{ytableau}
\end{split}\quad\quad\quad\quad
S_{\beta}= \ 
\begin{split}
\ytableausetup {mathmode, boxsize=1.0em} 
\begin{ytableau}
\none & \none & \none & \none & \none  \\
\tl{1} & \tl{1} & \tl{3} & \tl{3} & \tl{5}   \\
\tl{2} & \tl{2} & \tl{4} & \none & \none    \\
\none & \none & \none & \none & \none  \\
\end{ytableau}
\end{split}
\end{equation*}
We see that $S_{\alpha} \in \ov{\texttt{LR}}^{\lambda'}_{\alpha' \mu'}$ and $\psi(S_{\alpha}) \in \underline{\texttt{LR}}^{\lambda}_{\alpha \mu}$ 
(see Examples \ref{ex:example for LR var} and \ref{ex:bijection}).
On the other hand, for $S_{\beta}$, the sequence $(m_i)_{1 \le i \le 5}$ and $(n_j)_{1 \le j \le 3}$ are given by
$(1,3,5,6,8)$ and $(2,4,7)$, respectively.
Then $S_{\beta} \notin \ov{\texttt{LR}}^{\lambda'}_{\beta' \mu'}$ since $t_3 = 4 = \delta_{n_3}^{\texttt{rev}}$.
We can also check that $\psi(S_{\beta}) \notin \underline{\texttt{LR}}^{\lambda}_{\beta \mu}$ (cf. Example \ref{ex:bijection}).
By Theorem \ref{thm:non-levi branching}, we have 
\begin{equation*}
\left[ V_{\rm GL_8}^{\lambda} : V_{\rm O_8}^{\mu} \right] 
= \sum_{\delta \in \cP^{(2)}_8}\ov{c}^{\lambda}_{\delta\mu} 
= \sum_{\delta \in \cP^{(2)}_n}\underline{c}^{\lambda}_{\delta\mu}
= 1.
\end{equation*}
}
\end{ex}

\begin{rem} \label{rem:for symplectic groups}
{\em
\mbox{}

\begin{itemize}
	\item[(1)] For ${\rm G}_n = {\rm Sp}_n$, we also have
\begin{equation*}
	\left[ V_{\rm GL_n}^{\lambda} : V_{\rm Sp_n}^{\mu} \right]  = \sum_{\delta \in \cP_n^{(1,1)}}\ov{c}^{\lambda}_{\delta\mu}=\sum_{\delta \in \cP_n^{(1,1)}}\underline{c}^{\lambda}_{\delta\mu},
\end{equation*}
where $\ov{c}^{\lambda}_{\delta\mu}$ and $\underline{c}^{\lambda}_{\delta\mu}$ are given in Remark \ref{rem:for types B and C}. 
	\vskip 2mm
	
	\item[(2)] The flag condition in \eqref{eq:flag condition for types B and C} is different from the one in \cite[Section 6.3]{LL} because we use the bijection \eqref{conjugation of LR} whose image is the set of LR tableaux with anti-lattice word (cf. \cite[Theorem 6.2]{LL}).
\end{itemize}
}
\end{rem}
\vskip 2mm


\section{Generalized exponents}\label{sec:genexp}

\subsection{Generalized exponents} \label{subsec:generalized exponent}

Let $\mf g$ be a simple Lie algebra of rank $n$ over $\C$,
and $G$ the adjoint group of $\mf{g}$.
Let $S({\mf g})$ be the symmetric algebra generated by $\mf{g}$, 
and $S(\mf{g})^G$ the space of $G$-invariants with respect to the adjoint action.
Let $\mc{H}(\mf{g})$ be the space 
of polynomials annihilated by $G$-invariant differential operators with constant coefficients and no constant term. 
It is shown by Kostant \cite{Ko} that $S(\mf{g})$ is a free $S(\mf{g})^G$-module generated by $\mc{H}(\mf{g})$, that is,
\begin{equation*}
S(\mf{g})=S(\mf{g})^G\otimes \mc{H}(\mf{g}).
\end{equation*} 

Let $t$ be an indeterminate.
Let $\Phi^+$ denote the set of positive roots and $\Phi=\Phi^+\cup -\Phi^+$ the set of roots of $\mf g$. 
We define the graded character of $S(\mf{g})$ by
\begin{equation}\label{eq:char of S(g)}
{\rm ch}_t S(\mf{g}) =\dfrac{1}{(1-t)^n\prod_{\alpha\in \Phi}(1-te^\alpha)}.
\end{equation}
Then it is also shown in \cite{Ko} that the graded character of $\mc{H}(\mf{g})$ is determined by
\begin{equation}\label{eq:S(g) and H(g)}
{\rm ch}_t S(\mf{g}) = \frac{{\rm ch}_t\mc{H}(\mf{g})}{\prod_{i=1}^n(1-t^{d_i})},
\end{equation}
where $d_i=m_i+1$ for $i=1,\dots ,n$ and $m_i$ are the classical exponents of $\mf{g}$.

For $\mu\in P_+$, let $V_{\mf g}^{\mu}$ be the irreducible representation of $\mf{g}$ with highest weight $\mu$.
The generalized exponent associated to $\mu\in P_+$ is a graded multiplicity of $V^\mu_{\mf g}$ in $\mc{H}(\mf g)$, that is, 
\begin{equation*} 
E_t(V_{\mf g}^\mu) 
= \sum_{k\geq 0} \dim{\rm Hom}_{\mf{g}}(V_{\mf g}^{\mu},\mc{H}^k(\mf{g})) t^k,
\end{equation*}
where $\mc{H}^k(\mf{g})$ is the $k$-th homogeneous space of degree $k$.
It is shown in \cite{He} that 
\begin{equation*}%
E_t(V_{\mf g}^{\mu}) = K_{\mu 0}^{\mf g}(t),
\end{equation*}
where $K_{\mu 0}^{\mf g}(t)$ is the Lustig $t$-weight multiplicity for $V_{\mf g}^\mu$ at weight $0$. 
In other words, we have
\begin{equation*}
{\rm ch}_t\mc{H}(\mf{g}) = \sum_{\mu\in P_+} K_{\mu 0}^{\mf g}(t) {\rm ch}V_{\mf g}^{\mu}.
\end{equation*}

In \cite{LL}, a new combinatorial realization of $K_{\mu 0}^{\mf{sp}_{n}}(t)$ is given 
in terms of the LR tableaux which give a branching formula for \eqref{eq:Littlewood restriction} for ${\rm G}_n={\rm Sp}_n$. The goal of this section is to give combinatorial formulas for $K_{\mu 0}^{\mf{so}_{n}}(t)$ following the idea in \cite{LL} as a main application of Theorem \ref{thm:main1}.

\subsection{Combinatorial model of generalized exponents for $\mf{so}_n$} \label{subsec:combi model of generalized expoentns}
Suppose that $\mf{g}=\mf{so}_n$ for $n\geq 3$, that is, ${\mf g}=\mf{so}_{2m+1}, \mf{so}_{2m}$ for some $m$.
We assume that the weight lattice for $\mf{g}$ is $P=\bigoplus_{i=1}^m\Z\epsilon_i$ so that
$\Phi^+$ is 
$\{\,\epsilon_i\pm \epsilon_j, \epsilon_k  \,|\,1\leq i < j\leq m, 1\leq k\leq m\,\}$,  and
$\{\,\epsilon_i\pm \epsilon_j  \,|\,1\leq i < j\leq m\,\}$
when $\mf{g}=\mf{so}_{2m+1}$, and $\mf{g}=\mf{so}_{2m}$, respectively.
Let 
\begin{equation*}
\Delta^{\mf g}_t=
\dfrac{1}{\prod_{1\leq i<j\leq n}(1-tx_ix_j)}.
\end{equation*}
By using the Littlewood identity (when $t=1$), we have
\begin{equation}\label{eq:Littlewood identity}
\Delta^{\mf g}_t = 
\sum_{\la\in \cP_n^{(1,1)}}t^{|\la|/2} {\rm ch}V^\la_{{\rm GL}_n},
\end{equation}
where $|\la|=\sum_{i\geq 1}\la_i$ for $\la=(\la_i)_{i\geq 1}$.
Note that \eqref{eq:char of S(g)} can be obtained from $\Delta^{\mf g}_t$  by specializing it with respect to the torus of ${\rm SO}_n$ (see for example \cite[Section 2.2]{LL}).

For $\mu\in \mc{P}({\rm O}_n)$, put

\begin{equation*}\label{eq:double branching}
\left\ldbrack V^\la_{{\rm GL}_{n}}:V^\mu_{{\rm O}_{n}}\right\rdbrack = 
\begin{cases}
\left[V^\la_{{\rm GL}_{n}}:V^\mu_{{\rm O}_{n}}\right]+
\left[V^\la_{{\rm GL}_{n}}:V^{\ov{\mu}}_{{\rm O}_{n}}\right],& 
\text{if $\mu\neq \ov{\mu}$},\\
\left[V^\la_{{\rm GL}_{n}}:V^\mu_{{\rm O}_{n}}\right], & \text{if $\mu=\ov{\mu}$}.
\end{cases}
\end{equation*}

\begin{prop}\label{prop:formula for K}
For $\mu\in \cP_m$, we have 
\begin{equation*}
\begin{split}
\frac{K_{\mu 0}^{\mf{so}_{2m+1}}(t)}{\prod_{i=1}^m(1-t^{2i})} &=
\sum_{\la\in \cP^{(1,1)}_{2m+1}} 
\left\ldbrack V^\la_{{\rm GL}_{2m+1}}:V^\mu_{{\rm O}_{2m+1}}\right\rdbrack
t^{|\la|/2},\\
\frac{K_{\mu 0}^{\mf{so}_{2m}}(t)}{(1-t^m)\prod_{i=1}^{m-1}(1-t^{2i})} &=
\sum_{\la\in \cP^{(1,1)}_{2m}} 
\left\ldbrack V^\la_{{\rm GL}_{2m}}:V^\mu_{{\rm O}_{2m}}\right\rdbrack
t^{|\la|/2},
\end{split}
\end{equation*}
where we regard $\mu$ in $K^{\mf g}_{\mu 0}(t)$ as a dominant integral weight $\mu_1\epsilon_1+\dots + \mu_m\epsilon_m\in P_+$.
\end{prop}
\pf 
Suppose that $\mf{g}=\mf{so}_{2m+1}$. 
For $\mu\in \mc{P}({\rm O}_{2m+1})$ with $\ell(\mu)\leq m$, we have 
${\rm ch}V^\mu_{{\rm O}_{2m+1}}={\rm ch}V^{\ov{\mu}}_{{\rm O}_{2m+1}}={\rm ch}V^\mu_{\mf{so}_{2m+1}}$.
By taking restriction of \eqref{eq:Littlewood identity} with respect to ${\rm O}_{2m+1}$, we have
\begin{equation}\label{eq:S(g)-B}
{\rm ch}_t S(\mf{g})=
\sum_{\substack{\mu\in \mc{P}({\rm O}_{2m+1}) \\ \ell(\la)\leq m}}
\left(
\sum_{\la\in \cP^{(1,1)}_{2m+1}} t^{|\la|/2}
\left\ldbrack V^\la_{{\rm GL}_{2m+1}}:V^\mu_{{\rm O}_{2m+1}}\right\rdbrack 
\right)
{\rm ch}V^\mu_{\mf{so}_{2m+1}}.
\end{equation}

Next, suppose that $\mf{g}=\mf{so}_{2m}$.
For $\mu\in \mc{P}({\rm O}_{2m})$ with $\ell(\mu)< m$, we have 
${\rm ch}V^\mu_{{\rm O}_{2m}}={\rm ch}V^{\ov{\mu}}_{{\rm O}_{2m}}={\rm ch}V^\mu_{\mf{so}_{2m}}$. 
For $\mu\in \mc{P}({\rm O}_{2m})$ with $\ell(\mu)=m$, we have
${\rm ch}V^\mu_{{\rm O}_{2m}}=
{\rm ch}V^\mu_{\mf{so}_{2m}}+{\rm ch}V^{{\mu}^\sigma}_{\mf{so}_{2m}}$, where $\mu^\sigma = \mu_1\epsilon_1+\dots+\mu_{m-1}\epsilon_{m-1}-\mu_m\epsilon_m$.
Similarly, we have
\begin{equation}\label{eq:S(g)-D}
\begin{split}
{\rm ch}_t S(\mf{g})= &
\sum_{\substack{\mu\in \mc{P}({\rm O}_{2m}) \\ \ell(\mu)< m}}
\left(
\sum_{\la\in \cP^{(1,1)}_{2m}} t^{|\la|/2}
\left\ldbrack V^\la_{{\rm GL}_{2m}}:V^\mu_{{\rm O}_{2m}}\right\rdbrack 
\right)
{\rm ch}V^\mu_{\mf{so}_{2m}} \\
& +
\sum_{\substack{\mu\in \mc{P}({\rm O}_{2m}) \\ \ell(\mu)= m}}
\left(
\sum_{\nu\in \cP^{(1,1)}_{2m}} t^{|\nu|/2}
\left\ldbrack  V^\la_{{\rm GL}_{2m}}:V^\mu_{{\rm O}_{2m}}\right\rdbrack
\right)
\left({\rm ch}V^\mu_{\mf{so}_{2m}} + {\rm ch}V^{\mu^\sigma}_{\mf{so}_{2m}} \right).
\end{split}
\end{equation}

Now, combining \eqref{eq:S(g) and H(g)} and \eqref{eq:S(g)-B}, \eqref{eq:S(g)-D}, we obtained the identities.
\qed
\vskip 2mm

Suppose that $P=\bigoplus_{i=1}^n\Z \epsilon_i$ is the weight lattice of $\mf{gl}_n$.
For $1\leq i\leq n-1$, let $\varpi_i=\epsilon_1+\dots+\epsilon_i$ be the $i$th fundamental weight. 

Let $\mu \in \cP_n$ be given. We identify $\mu$ with $\mu_1\epsilon_1+\cdots+\mu_n\epsilon_{n}$.
Let $SST_n(\mu)$ (resp. $SST_n(\mu^\pi)$) be the  subset of $SST(\mu)$ (resp. $SST(\mu^\pi)$) consisting of $T$ with entries in $\{1,\dots,n\}$, which is a $\mf{gl}_n$-crystal with highest weight $\mu$. 
For $T\in SST_n(\mu)$ or $SST_n(\mu^\pi)$, put
\begin{equation*}
\varphi(T)= \sum_{i=1}^{n-1}\varphi_i(T)\varpi_i, \quad
\varepsilon(T)= \sum_{i=1}^{n-1}\varepsilon_i(T)\varpi_i.
\end{equation*}

\begin{df}\footnote{In this paper, we use the notation $\rho$ as a partition, not the half sum of positive roots.} \label{df:distinguished}
{\rm 
For $\rho \in \cP_n$, we say that $T$ is {\em $\rho$-distinguished} if 
\begin{equation*}
\varphi(T)=\la -\rho, \quad \varepsilon(T)= \delta -\rho,
\end{equation*}
for some $(\la,\delta)\in \cP^{(1,1)}_n\times \cP^{(2)}_n$.
}
\end{df}

We put
\begin{equation}
\begin{split}
D_n(\mu) & = \{\,T\in SST_n(\mu^{\pi})\,|\, \text{$T$ is $\rho$-distinguished for some $\rho\in\cP_n$}\,\},\\
\cP_T & = \{\,\rho\in \cP_n\,|\, \text{$T$ is $\rho$-distinguished}\,\}\quad\text{($T\in D_n(\mu)$)}.
\end{split}
\end{equation}

\begin{lem}\label{lem:rho_T}
For $T \in D_n(\mu)$, there exists a unique $\rho_T \in \cP$ such that $\cP_T = \rho_T + \cP^{(2,2)}_n$, where $\rho_T$ is determined by 
\begin{equation*}
\rho_T = \sum_{\substack{1\leq i\leq n-1\\ i\equiv 0\mod{2}}}(\varepsilon_{i}(T) \ \textrm{mod $2$})\varpi_{i}.
\end{equation*}
\end{lem}
\pf It follows from the same argument as in \cite[Lemma 4.4, Proposition 4.5]{LL} with Definition \ref{df:distinguished}. 
\qed
\vskip 2mm

\begin{df} \label{df:distinguished elements with flag conditions}
{\em 
We define $\underline{D}_n(\mu)$ to be the subset of $D_n(\mu)$ consisting of $T$ satisfying the condition \eqref{eq:flagged conditions}.
}
\end{df}
\vskip 2mm

\begin{prop}\label{prop:formula for K-2}
For $\mu\in \cP_m$, we have
\begin{equation*}
\sum_{\la\in \cP^{(1,1)}_{n}} 
\left\ldbrack V^\la_{{\rm GL}_{n}}:V^\mu_{{\rm O}_{n}}\right\rdbrack t^{|\la|/2}
=\frac{1}{\prod_{i=1}^m(1-t^{2i})}\sum_{T\in \mathbb{D}_n(\mu)}t^{|\varphi(T)+\rho_T|/2},
\end{equation*}
where
\begin{equation} \label{eq:distinguished set}
\mathbb{D}_n(\mu)=
\begin{cases}
\underline{D}_n(\mu)\sqcup \underline{D}_n(\ov{\mu}), & \text{if $\mu\neq \ov{\mu}$},\\
\underline{D}_n(\mu), & \text{if $\mu=\ov{\mu}$}.
\end{cases}
\end{equation}
\end{prop}
\pf 
Recall that we have bijections for $\mu\in \mc{P}({\rm O}_n)$
\begin{equation*}
\bigsqcup_{\la\in\cP_n^{(1,1)}}\texttt{LR}^{\mu}_{\lambda}(\mathfrak{d})
\quad\longrightarrow
\bigsqcup_{\la\in\cP_n^{(1,1)}}\bigsqcup_{\delta \in \cP_{n}^{(2)}} \ov{\texttt{LR}}^{\lambda'}_{\delta'\mu'}
\quad\longrightarrow
\bigsqcup_{\la\in\cP_n^{(1,1)}}\bigsqcup_{\delta \in \cP_{n}^{(2)}}\underline{\texttt {LR}}^{\la}_{\delta\mu},
\end{equation*}
where the first one is given in Theorem \ref{thm:main1} and the second one in Theorem \ref{thm:main result-flagged}.
By definition of $\underline{D}_n(\mu)$, we have a bijection
\begin{equation} \label{eq:bijection for ge-1}
\xymatrixcolsep{3pc}\xymatrixrowsep{0pc}\xymatrix{
\bigsqcup_{\la\in\cP_n^{(1,1)}}\bigsqcup_{\delta \in \cP_{n}^{(2)}}\underline{\texttt {LR}}^{\la}_{\delta\mu}
 \ar@{->}[r] & \bigsqcup_{T\in \underline{D}_n(\mu)}\{\,T\,\}\times\cP_T \\ 
T \ar@{|->}[r] & (T,\la-\varphi(T))=(T, \delta-\varepsilon(T))
}.
\end{equation}
By Lemma \ref{lem:rho_T}, we have a bijection
\begin{equation}\label{eq:bijection for ge-2}
\xymatrixcolsep{3pc}\xymatrixrowsep{0pc}\xymatrix{
\bigsqcup_{T\in \underline{D}_n(\mu)}\{\,T\,\}\times\cP_T  \ar@{->}[r] &  \bigsqcup_{T\in \underline{D}_n(\mu)}\{\,T \,\}\times\cP_n^{(2,2)}  \\ 
(T,\rho)  \ar@{|->}[r] & (T,\rho-\rho_T)
}.
\end{equation}
Therefore, we have from \eqref{eq:bijection for ge-1} and \eqref{eq:bijection for ge-2}

\begin{equation*}
\begin{split}
\sum_{\la\in \cP^{(1,1)}_{n}} 
\left\lbrack V^\la_{{\rm GL}_{n}}:V^\mu_{{\rm O}_{n}}\right\rbrack t^{|\la|/2}
&= \sum_{\la\in \cP^{(1,1)}_{n}}\sum_{\delta \in \cP^{(2)}_n}\underline{c}^{\lambda}_{\delta\mu} t^{|\la|/2} 
= \sum_{T\in \underline{D}_n(\mu)}\sum_{\rho\in \cP_T} t^{|\varphi(T)+\rho|/2} \\
&= \sum_{T\in \underline{D}_n(\mu)}t^{|\varphi(T)+\rho_T|/2}\sum_{\kappa\in \cP_n^{(2,2)}}t^{|\kappa|/2} \\
&= \sum_{T\in \underline{D}_n(\mu)}t^{|\varphi(T)+\rho_T|/2} \frac{1}{\prod_{i=1}^m(1-t^{2i})},
\end{split}
\end{equation*}
which implies the identity.
\qed\vskip 2mm

We have the following new combinatorial formulas for $K^{\mf{so}_{n}}_{\mu 0}(t)$.

\begin{thm}\label{thm:Kostka-Foulkes for BD}
For $\mu\in \cP_m$, we have
\begin{equation*}
\begin{split}
K^{\mf{so}_{2m+1}}_{\mu 0}(t)& = \sum_{T\in \mathbb{D}_{2m+1}(\mu)}t^{|\varphi(T)+\rho_T|/2}, \\
K^{\mf{so}_{2m}}_{\mu 0}(t)& = \frac{1}{1+t^m} \sum_{T\in \mathbb{D}_{2m}(\mu)}t^{|\varphi(T)+\rho_T|/2},
\end{split}
\end{equation*}
where $\mathbb{D}_n(\mu)$ is given in \eqref{eq:distinguished set}.
\end{thm}
\pf It follows from Propositions \ref{prop:formula for K} and \ref{prop:formula for K-2}.
\qed

\begin{rem}{\rm
Since $K^{\mf{so}_{2m}}_{\mu 0}(t)$ is a polynomial in $t$, the polynomial
$$
\sum_{T\in \mathbb{D}_{2m}(\mu)}t^{|\varphi(T)+\rho_T|/2}
$$
is divisible by $1+t^m$. 
From the positivity of Kostka-Foulkes polynomial $K^{\mf{so}_{2m}}_{\mu 0}(t)$, one may expect a decomposition of $\mathbb{D}_{2m}(\mu)= X_1\sqcup X_2$ together with a bijection $\tau : X_1 \longrightarrow X_2$ such that 
$$
|\varphi(\tau(T))+\rho_{\tau(T)}|= 2m +|\varphi(T)+\rho_T|.
$$
}
\end{rem}
\vskip 2mm

\begin{rem}
{\em 
In \cite{LL}, Lecouvey--Lenart provide a bijection between the distinguished tableaux for type $C_n$ and the symplectic King tableaux with weight $0$ (see \cite[Section 6.4]{LL} for more details).

We do not know whether there is an analogue of the above bijection which maps an orthogonal distinguished tableau (Definition \ref{df:distinguished elements with flag conditions}) to an orthogonal tableau (with weight 0), which is from already known models (for example, \cite{KS,KW, KT,O1,P2, Su2}) or a new one. In this case, the flag condition \ref{eq:flagged conditions} would be quite complicated under this correspondence if exists.
}
\end{rem}


\section{Proof of Main Theorem}\label{sec:proof of main}

\subsection{Outline}
The proof of Theorem \ref{thm:main1} is organized as follows. 

In subsection \ref{subsec:pf_main1}, we consider the case of $n-2\mu_1' \ge 0$, which is easier to deal with.
\vskip 2mm

\begin{itemize}
	\item[(1)] ({\em Well-definedness})
	First, we show that $\ov{\bf T}^{\texttt{tail}}\in \ov{\texttt{{LR}}}^{\lambda'}_{\delta'\mu'}$ (Corollary \ref{cor:well-definedness}).
	To do this, we study some properties of the sequences $(m_i)_{1 \le i \le p}$ and $(n_j)_{1 \le j \le q}$  associated to $\ov{\bf T}^{\texttt{tail}}$ with respect to sliding (Lemmas \ref{lem:sequence m_i} and \ref{lem:sequence n_i}), which implies that $\ov{\bf T}^{\texttt{tail}}$ satisfies \eqref{eq:condition_on_second_row}.
	\vskip 2mm
	
	\item[(2)] ({\em Injectivity})
	Second, we show that the map 
\begin{equation*} \label{eq:bijection in 6.1}
\xymatrixcolsep{3pc}\xymatrixrowsep{0pc}\xymatrix{
\texttt{{LR}}^{\mu}_{\lambda}(\mathfrak{d}) \ar@{->}[r] & \ \bigsqcup_{\delta \in \cP_{n}^{(2)}} \ov{\texttt{{LR}}}^{\lambda'}_{\delta'\mu'}  \\ 
\mathbf{T}   \ar@{|->}[r] & \ov{\bf T}^{\texttt{tail}} 
}.
\end{equation*}	
 is injective by using Proposition \ref{prop:body and tail} (Lemma \ref{lem:injective}).
 \vskip 2mm
 
 	\item[(3)] ({\em Surjectivity})
	Finally, we prove the above map is surjective, that is, for ${\bf W} \in \ov{\texttt{{LR}}}^{\lambda'}_{\delta'\mu'}$, there exists ${\bf T} \in \texttt{{LR}}^{\mu}_{\lambda}(\mathfrak{d})$ such that $\ov{\bf T}^{\texttt{tail}} = {\bf W}$.
	We use induction on $n$. The initial step when $n=4$ is proved in Lemma \ref{lem:fundamental case}. 
	Then based on this step, we construct ${\bf T} \in \texttt{{LR}}^{\mu}_{\lambda}(\mathfrak{d})$ in general in Lemma \ref{lem:surjective}
	
\end{itemize}
\vskip 2mm

In subsection \ref{subsec:negative case proof}, we consider the case of $n-2\mu_1' < 0$.
The proof is almost identical to the case of $n-2\mu_1' \ge 0$, but the major difficulty occurs when we consider the columns with odd height in $\overline{\mathbf{T}}(0)$ and $\mathbf{T}^{\textrm{sp}-}$ (cf.~Remark \ref{rem:admissibilty for spin column}(2)--(3)). 
To overcome this, we reduce the problems to the ones in the case of $n-2\mu_1' \ge 0$ so that we may apply the results (or the arguments in the proof) in subsection \ref{subsec:pf_main1}. 
\vskip 2mm

\subsection{Proof of Theorem \ref{thm:main1} when $n-2\mu_1' \ge 0$} \label{subsec:pf_main1}
Let $\mu\in {\mc P}({\rm O}_n)$ and $\la\in \cP_{n}$ be given. We assume that $n-2\mu'_1 \ge 0$. 
The proof for the case $n-2\mu'_1 < 0$ will be given in subsection \ref{subsec:negative case proof}. We keep the notations in Sections \ref{sec:separation} and \ref{sec:branching}. 

Suppose that $n=2l+r$, where $l\geq 1$ and $r=0,1$. Let ${\bf T} \in \texttt{LR}^{\mu}_{\lambda}(\mf d)$ be given 
with ${\bf T} = (T_l, \dots, T_1,T_0)$ as in \eqref{eq:T notation}. 
Let us assume that $r=0$ since the argument for $r=1$ is almost identical.
Write $\td{\bf T} = (\td{T}_{l-1}, \dots, \td{T}_1,\td{T}_0)$.
Let $\ov{\bf T}^{\texttt{body}}=H_{(\delta')^\pi}$ for some $\delta\in \cP^{(2)}$.
Let $s_1\leq \dots \leq s_p$ denote the entries in the first row, and $t_1\leq \dots \leq t_q$ the entries in the second row of $\ov{\bf T}^{\texttt{tail}}$.

\begin{lem} \label{lem:trivial invariant}
Suppose that $\ov{\bf T}=(U_{2l},\dots,U_1)\in {\bf E}^n$ under \eqref{eq:identification}.
If $T_{i+1}^{\texttt{\em R}}(1)<T_i^{\texttt{\em L}}(a_i)$, then
	\begin{equation*}
		U_{2i} = T_{i}^{\texttt{\em L}} \boxminus T_{i}^{\texttt{\em tail}}.
	\end{equation*}
In this case, $T_{i}^{\texttt{\em tail}}$ is the $(l-i+1)$-th column of $\ov{\bf{T}}^{\texttt{\em tail}}$ from the left.
\end{lem}
\pf 
If $T_{i+1}^{\texttt{R}}(1) < T_i^{\texttt{L}}(a_i)$, then by definition we have $\widetilde{T}_i$ has residue $0$. By Lemma \ref{lem:sliding one step}, $\widetilde{T}_i^{\texttt{R}}(1)<\widetilde{T}_{i-1}^{\texttt{L}}(a_{i-1})$. Inductively, we have $U_{2i} = T_{i}^{\texttt{L}} \boxminus T_{i}^{\texttt{tail}}$. By applying this argument together with Lemma \ref{lem:sliding one step}, we obtain the second statement.
\qed
\vskip 2mm

For simplicity, let us put ${\mathbb T}=\widetilde{\bf T}$. 
Let $\td{\mu}=(\mu_2, \mu_3, \dots)$ and $\zeta=(\delta_1, \dots, \delta_{n-1}) \in \cP_{n-1}^{(2)}$. 
By Lemmas \ref{lem:criterion_highest_weight_elt}, \ref{lem:sliding one step} and
Proposition \ref{prop:body and tail}, 
we have 
$\ov{\mathbb{T}}^{\texttt{body}}=H_{(\zeta')^\pi}$ and
$\ov{\mathbb T}^{\texttt{tail}} \in \texttt{LR}^{\xi'}_{\zeta'\td{\mu}'}$
where $\xi$ is given by $(\ov{\mathbb{T}}^{\texttt{tail}} \rightarrow H_{\zeta'}) = H_{\xi'}$.
Let $\widetilde{s}_1\leq \dots \leq \widetilde{s}_{p-1}$ be the entries in the first row of $\ov{\mathbb T}^{\texttt{tail}}$ and let $(\widetilde{m}_i)_{1 \le i \le p-1}$ (resp. $(m_i)_{1 \le i \le p}$) be the sequence associated to $\ov{\mathbb{T}}^{\texttt{tail}}$ (resp. $\ov{\bf T}^{\texttt{tail}}$) in Definition \ref{bounded orthogonal LR}. 
Note that $s_i = \widetilde{s}_{i-1}$ for $2\leq i\leq p$.
Put $\texttt{T}_i = T_{l-i+1}$ for $1\leq i\leq l$ and  
$\td{\texttt{T}}_j = \td{T}_{l-j}$ for $1\leq j\leq l-1$.
Assume that $\texttt{T}_i\in {\bf T}(\texttt{a}_i)$ for $1\leq i\leq l$.

\begin{lem} \label{lem:sequence m_i}
Under the above hypothesis, the sequences $(m_i)_{1 \le i \le p}$ and $(\widetilde{m}_i)_{1 \le i \le p-1}$ satisfy the relation
\begin{equation} \label{eq:inductive m_i}
m_i =  \widetilde{m}_{i-1}+\tau_i+1 \quad (2 \le i \le p),
\end{equation} 
where $\tau_i$ is given by
\begin{equation*} %
	\tau_i =  
	\begin{cases} 
			1, & \textrm{{\em if} $\texttt{\em T}_{i-1}^{\texttt{\em R}}(1)<\texttt{\em T}_i^{\texttt{\em L}}(\texttt{\em a}_i)$,} \\
			0, & \textrm{{\em if} $\texttt{\em T}_{i-1}^{\texttt{\em R}}(1)>\texttt{\em T}_i^{\texttt{\em L}}(\texttt{\em a}_i)$.} 
 	\end{cases}
\end{equation*}
\end{lem}
\pf 
Fix $i \ge 2$.
If $\texttt{T}_{i-1}^{\texttt{R}}(1)<\texttt{T}_i^{\texttt{L}}(\texttt{a}_i)$, then
by Lemma \ref{lem:description of S}(i) and \ref{lem:trivial invariant}
\begin{equation*}
	m_i = 2i-1, \quad \quad \widetilde{m}_{i-1} = 2i-3.
\end{equation*}
If $\texttt{T}_{i-1}^{\texttt{R}}(1)>\texttt{T}_i^{\texttt{L}}(\texttt{a}_i)$, then we have by Lemma \ref{lem:description of S}(ii) $m_i < 2i-1$. This implies $m_i = \widetilde{m}_{i-1}+1$. Hence we have \eqref{eq:inductive m_i}. 
\qed

\begin{lem} \label{lem:sequence n_i}
For $1\leq i \leq q$, we have
\begin{equation*} \label{eq:description of second condition}
		t_i = \texttt{\em T}_i^{\texttt{\em L}}(\texttt{\em a}_i-1) > \texttt{\em T}^{\texttt{\em R}}_{i}(1) = \delta^{\texttt{\em rev}}_{n_i}.
\end{equation*}
\end{lem}
\pf By Lemma \ref{lem:description of S}, it is easy to see that $t_i = \texttt{T}_i^{\texttt{L}}(\texttt{a}_i-1)$. 
Next we claim that $\texttt{T}^{\texttt{R}}_{i}(1) = \delta^{\texttt{rev}}_{n_i}$, which implies the inequality since ${\mf r}_{\texttt{T}_{i}}\leq 1$.
We use induction on $n$.
For each $i$, we define $\theta_i$ to be the number of $j$'s with $i+1\leq j\leq p$ such that $m_j < 2i+1$. 
Then we have
\begin{equation} \label{eq:description of n_i}
n_i = 2i + \theta_i.
\end{equation}

If $\theta_i = 0$, then $n_i = 2i$ and $m_{i+1}=2i+1$, which implies that
\begin{equation*}
	\texttt{T}_i^{\texttt{R}}(1) < \texttt{T}_{i+1}^{\texttt{L}}(\texttt{a}_{i+1}).	
\end{equation*} By applying Lemma \ref{lem:trivial invariant} on $\widetilde{\bf T}$, we have $\delta_{2i}^{\texttt{rev}} = \texttt{T}_i^{\texttt{R}}(1)$.

If $\theta_i > 0$, then we have by definition of $\theta_i$, 
\begin{equation} \label{eq:theta > 0}
\texttt{T}_{i-1}^{\texttt{R}}(1) > \texttt{T}_i^{\texttt{L}}(\texttt{a}_i).	
\end{equation}
Let $(\td{m}_i)_{1\leq i\leq p-1}$ and $(\td{n}_i)_{1\leq i\leq q-1}$ be the sequences in Definition \ref{bounded orthogonal LR} associated to $\widetilde{\bf T}$. 
Let $\widetilde{\theta}_i$ be defined in the same way with respect to $(\widetilde{m}_i)_{1\leq i\leq p-1}$. 
By definition of $\widetilde{\theta}_i$ and Lemma \ref{lem:sequence m_i}, we have {for $j \ge i+2$},
\begin{equation*}
	m_j < 2i+1 \ \Longrightarrow \ \widetilde{m}_{j-1} < 2i-\tau_j \le 2i-1. 	
\end{equation*} 
Thus we have $\widetilde{\theta}_i = \theta_i-1$.
By induction hypothesis, \eqref{eq:description of n_i}, and \eqref{eq:theta > 0}, we have
\begin{equation*}
	\texttt{T}^{\texttt{R}}_i(1) = \widetilde{\texttt{T}}_i^{\texttt{R}}(1) = \widetilde{\delta}_{\widetilde{n}_i}^{\texttt{rev}}=\widetilde{\delta}_{2i+\td{\theta}_i}^{\texttt{rev}} = \delta_{2i+\theta_i}^{\texttt{rev}} = \delta_{n_i}^{\texttt{rev}}.
\end{equation*}
\qed

\begin{cor}\label{cor:well-definedness}
Under the above hypothesis, we have $\ov{\bf T}^{\texttt{\em tail}}\in \ov{\texttt{\em LR}}^{\la'}_{\mu'\delta'}$.
\end{cor}
\pf 
It follows from Remark \ref{rem:existence m_i} and Lemma \ref{lem:sequence n_i}.
\qed

\begin{lem} \label{lem:injective}
The map ${\bf T} \longmapsto \ov{\bf T}^{\texttt{\em tail}}$ is injective on $\texttt{{\em LR}}^{\mu}_{\lambda}(\mathfrak{d})$.
\end{lem}
\pf Let ${\bf T}, {\bf S} \in \texttt{{LR}}^{\mu}_{\lambda}(\mathfrak{d})$ be given. 
Suppose that $\ov{\bf T}^{\texttt{tail}} = \ov{\bf S}^{\texttt{tail}}$.
We first claim that $\ov{\bf T}=\ov{\bf S}$. 
By Proposition \ref{prop:body and tail}(1), we have 
$\ov{\bf T}^{\texttt{body}}=H_{\delta'}$ and 
$\ov{\bf S}^{\texttt{body}}=H_{\chi'}$ for some $\delta, \chi\in \cP^{(2)}$.
Since $\ov{\bf T}^{\texttt{tail}} = \ov{\bf S}^{\texttt{tail}}$ and
$(\ov{\bf T}^{\texttt{tail}}\rightarrow H_{\delta'})=
(\ov{\bf S}^{\texttt{tail}}\rightarrow H_{\chi'})=H_{\la'}$,
we have $\delta=\chi$. Hence $\ov{\bf T}^{\texttt{body}}=\ov{\bf S}^{\texttt{body}}$, which implies $\ov{\bf T}=\ov{\bf S}$. Since the map ${\bf T}\mapsto \ov{\bf T}$ is reversible, we have ${\bf T}={\bf S}$. \qed
\vskip 3mm

Now, we will verify that the map in Theorem \ref{thm:main1}  is surjective.
Let ${\bf W} \in \ov{\texttt{LR}}^{\lambda'}_{\delta'\mu'}$ be given for some $\delta \in \cP_n^{(2)}$.
Let ${\bf V}= H_{(\delta')^{\pi}}$ and ${\bf X}$ be the tableaux of a skew shape $\eta$ as in \eqref{eq:shape after separation} with $n$ columns such that 
$$
{\bf X}^{\texttt{body}}={\bf V},\quad {\bf X}^{\texttt{tail}}={\bf W}.
$$
The semistandardness of ${\bf X}$ follows from Definition \ref{bounded orthogonal LR} and Remark \ref{rem:existence m_i}. %
Let $V_i$ and $W_i$ denote the $i$-th column of ${\bf V}$ and ${\bf W}$ from right, respectively.

Let us first consider the following, which is used in the proof of Lemma \ref{lem:surjective}.

\begin{lem} \label{lem:fundamental case}
Assume that $n=4$ and $\mu_1' = 2$. 
Then there exists ${\bf T} = (T_2, T_1) \in \texttt{\em LR}^{\mu}_{\lambda}(\mf d)$ such that 
$\ov{\bf T}={\bf X}$, that is, 
$\ov{\bf T}^{\texttt{\em body}} = {\bf V}$ and $\ov{\bf T}^{\texttt{\em tail}} = {\bf W}$. In fact, ${\bf T} = (T_2, T_1)$ is given as follows:
\begin{itemize}
	\item[(1)] If $m_2 = 3$, then
		\begin{equation*}
			(T_2^{\texttt{\em L}},T_2^{\texttt{\em R}}) = (V_4 \boxplus W_2, V_3), \ \ (T_1^{\texttt{\em L}},T_1^{\texttt{\em R}}) = (V_2 \boxplus W_1, V_1).
		\end{equation*}
		
	\item[(2)] If $m_2 = 2$, then
		\begin{equation*}
			(T_2^{\texttt{\em L}},T_2^{\texttt{\em R}}) = (V_4 \boxplus W_2,V_3^{\diamond}), \quad (T_1^{\texttt{\em L}},T_1^{\texttt{\em R}}) = (V_2^{\diamond} \boxplus W_1^{\diamond}, V_1),	
		\end{equation*} 
		where $V_3^{\diamond}$, $V_2^{\diamond}$ and $W_1^{\diamond}$ are given by
		\begin{equation*}
		\begin{split}
			& V_3^{\diamond} = \left(\dots, V_3(2), V_3(1), W_1(a_1), V_2(1)\right) \boxplus \emptyset, \\
			& V_2^{\diamond} = \left(\dots, V_2(4), V_2(3)\right), \quad W_1^{\diamond} = (W_2(2), W_1(a_1-1), \dots W_1(1)).
		\end{split}
		\end{equation*} 
\end{itemize}
\end{lem}
\pf By Remark \ref{rem:existence m_i} and definition of them, $T_1$ and $T_2$ are semistandard. Also, the residue ${\mf r}_i$ of $T_i$ is by Definition \ref{bounded orthogonal LR} at most $1$ for $i = 1, 2$.
It suffices to verify that $T_2\prec T_1$ since this implies 
$\ov{\bf T}^{\texttt{body}} = {\bf V}$ and 
$\ov{\bf T}^{\texttt{tail}} = {\bf W}$ by construction of ${\bf T}$.

Let $a_i$ be the height of $W_{i}$ for $i=1,2$. 
We have 
$V_i = (1, 2, \dots, \delta_i^{\texttt{rev}})$ for $1\leq i\leq 4$,
with 
$\delta_1^{\texttt{rev}} \le \delta_2^{\texttt{rev}} \le  \delta_3^{\texttt{rev}} \le \delta_4^{\texttt{rev}}$. 
Let $w({\bf X}) = w_1w_2\cdots w_m$. Put 
\begin{equation*}
	P_k = (w_k \rightarrow ( \cdots \rightarrow (w_3 \rightarrow (w_2 \rightarrow (w_1)))))
\end{equation*} 
for $k \le m$. Suppose that
\begin{equation*}
\begin{split}
	& w(V_1)w(V_2)w(V_3) = w_1w_2\cdots w_{s}, \\
	&  w(V_1)w(V_2)w(V_3)w(W_1)w(V_4)=w_1w_2\cdots w_{t},
\end{split}	
\end{equation*} 
for some $s \le t \le m$. 

{\bf\em Case 1. $m_2 = 3$.}
We first assume that ${\mf r}_1{\mf r}_2 = 0$.
	\begin{itemize}
		\item[(i)] ${\mf r}_1 = 0, \ {\mf r}_2=0$ : 
			It is obvious that $T_2 \prec T_1$.
		
		\item[(ii)]	${\mf r}_1 = 0, \ {\mf r}_2=1$ :
			  Definition \ref{def:admissibility}(1)-(i) follows from $\delta_2^{\texttt{rev}} \le \delta_3^{\texttt{rev}}$. Also, Definition \ref{def:admissibility}(1)-(ii) follows from $s_2 \ge \delta_3^{\texttt{rev}} \ge \delta_2^{\texttt{rev}}$. The semistandardness of ${\bf W}$ implies Definition \ref{def:admissibility}(1)-(iii). Thus we have  $T_2 \prec T_1$.
		
		\item[(iii)] ${\mf r}_1 = 1, \ {\mf r}_2=0$ : We may use the same argument as in (ii) to have $T_2 \prec T_1$.
	\end{itemize} 
	
Next, we assume that ${\mf r}_1{\mf r}_2 = 1$. 
Definition \ref{def:admissibility}(1)-(i) holds by definition of ${\bf T}$. 
Since ${\mf r}_1 = 1$ and ${\mf r}_2 = 1$, we have
\begin{equation} \label{eq:range of V_i(a_i)}
\delta_3^{\texttt{rev}} < W_1(a_1) \le \delta_4^{\texttt{rev}}, \quad  
\delta_1^{\texttt{rev}} < W_2(a_2) \le \delta_2^{\texttt{rev}}.
\end{equation}
By Lemma \ref{lem:criterion_highest_weight_elt},
\begin{equation} \label{eq:insertion V_i(a_i)}
\text{$(W_1(a_1) \rightarrow P_s)$ and $(W_2(a_2) \rightarrow P_t)$ 
are $\mf l$-highest weight elements}.
\end{equation} 
By \eqref{eq:range of V_i(a_i)} and \eqref{eq:insertion V_i(a_i)}, we have
\begin{equation} \label{eq:case 1 V_1(a_1) and V_2(a_2)}
W_1(a_1) = \delta_3^{\texttt{rev}}+1, \quad W_2(a_2) = \delta_1^{\texttt{rev}}+1.
\end{equation} 
This implies Definition \ref{def:admissibility}(1)-(ii). 
Definition \ref{def:admissibility}(1)-(iii) follows from 
the semistandardness of ${\bf W}$ and \eqref{eq:case 1 V_1(a_1) and V_2(a_2)}. 
Thus we have $T_2 \prec T_1$.
\vskip 3mm

{\bf\em Case 2. $m_2 = 2$.}
Since $m_2 = 2$, we have 
\begin{equation} \label{eq:case 2 delta_2 < delta_3}
\delta_2^{\texttt{rev}} < \delta_3^{\texttt{rev}}.	
\end{equation}
Otherwise, we have $W_1(a_1) > \delta_3^{\texttt{rev}}$, which contradicts to $m_2=2$. 
By definition of $m_2$, we have
\begin{equation} \label{eq:range of V_1(a_1)}
\delta_2^{\texttt{rev}} < W_1(a_1) < \delta_3^{\texttt{rev}}.
\end{equation} 
Note that if $W_1(a_1) = \delta_3^{\texttt{rev}}$, then by \eqref{eq:case 2 delta_2 < delta_3} the tableau $(W_1(a_1) \rightarrow P_s)$ cannot be an $\mf l$-highest weight element.
Since $(W_1(a_1) \rightarrow P_s)$ is an $\mf l$-highest weight element, 
\begin{equation} \label{eq:case 2 V_1(a_1)}
W_1(a_1) = \delta_2^{\texttt{rev}}+1.	
\end{equation} 
In particular, we have 
\begin{equation}\label{eq:delta_2+2 <= delta_3}
\delta_2^{\texttt{rev}}+2 \le \delta_3^{\texttt{rev}}.
\end{equation}
Since $W_2(a_2) \le W_1(a_1) < \delta_3^{\texttt{rev}}$ and 
$W_1(a_1) = \delta_2^{\texttt{rev}}+1$, we have ${\mf r}_2=1$. Also, since $\delta_3^{\texttt{rev}} \le \delta_4^{\texttt{rev}}$, it is clear that ${\mf r}_1=1$. 
Note that
\begin{equation*}
\delta_1^{\texttt{rev}} < W_2(a_2) \le W_1(a_1) = \delta_2^{\texttt{rev}}+1 < \delta_3^{\texttt{rev}}.
\end{equation*}
This implies that
\begin{equation}\label{eq:case 2 V_2(a_2)}
W_2(a_2) = \delta_1^{\texttt{rev}}+1,
\end{equation} 
since $(W_2(a_2) \rightarrow P_t)$ is an $\mf l$-highest weight element.

Now, Definition \ref{def:admissibility}(1)-(i) follows from \eqref{eq:delta_2+2 <= delta_3}, and
Definition \ref{def:admissibility}(1)-(ii) and (1)-(iii) follow from 
\eqref{eq:case 2 V_1(a_1)}, \eqref{eq:delta_2+2 <= delta_3}, \eqref{eq:case 2 V_2(a_2)} and the semistandardness of ${\bf W}$. Hence we have $T_2\prec T_1$.
\qed
\vskip 3mm

Let ${\mathbb X}$ be the tableau obtained from ${\bf X}$ by removing its leftmost column. Let $\td{\mu}=(\mu_2, \mu_3, \dots)$ and $\zeta=(\delta_1, \dots, \delta_{n-1}) \in \cP_{n-1}^{(2)}$. 
Since ${\mathbb X}$ is an $\mf l$-highest weight element
by Lemma \ref{lem:criterion_highest_weight_elt},
we have 
${\mathbb X}^{\texttt{body}}=H_{(\zeta')^\pi}$ and
${\mathbb X}^{\texttt{tail}} \in \texttt{LR}^{\xi'}_{\zeta'\td{\mu}'}$,
where $\xi$ is given by $({\mathbb X}^{\texttt{tail}} \rightarrow H_{\zeta'}) = H_{\xi'}$.

\begin{lem} \label{lem:inductive step}
We have
$
{\mathbb X}^{\texttt{\em tail}} \in \ov{\texttt{\em LR}}^{\xi'}_{\zeta'\td{\mu}'}.
$
\end{lem}
\pf
Let $(m_i)_{1\leq i\leq p}$ and $(n_i)_{1\leq i\leq q}$ be the sequences associated to ${\bf X}^{\texttt{tail}}={\bf W} \in \ov{\texttt{LR}}^{\lambda'}_{\delta'\mu'}$.
Let $\td{s}_1\leq \dots \leq \td{s}_{p-1}$ and $\td{t}_1\leq \dots \leq \td{t}_{q-1}$ be the entries in the first and second rows of ${\mathbb X}^{\texttt{tail}}$, respectively.

We define a sequence $1 \le \widetilde{m}_1 < \dots < \widetilde{m}_{p-1} \le n-1$ inductively as in Definition \ref{bounded orthogonal LR} with respect to $(\widetilde{s}_i)_{1 \le i \le p-1}$. 
Note that the sequence $(\widetilde{m}_i)_{1 \le i \le p-1}$ is well-defined by Remark \ref{rem:existence m_i}.
By Lemma \ref{lem:description of S}, we observe that
\begin{equation} \label{eq:m_i induction}
	\td{m}_i = \left\{
		\begin{array}{ll}
			1, & \textrm{if $i=1$}, \\
			m_{i+1}-1, & \textrm{if $i>1$ and $m_{i+1} < 2i+1$}, \\
			m_{i+1}-2, & \textrm{if $i>1$ and $m_{i+1} = 2i+1$}.
		\end{array}
		\right.
\end{equation} 
Let $(\td{n}_i)_{1\leq i\leq q-1}$ be the sequence with respect to $(\td{m}_i)_{1\leq i\leq p-1}$, that is,
\begin{equation*}
	\td{n}_i = \textrm{the $i$-th smallest integer in $\{ i+1, \dots, n-1 \} \setminus \{ \td{m}_{i+1}, \dots,\td{m}_{p-1}\}$}.	
\end{equation*} 
By \eqref{eq:m_i induction}, we obtain 
\begin{equation} \label{eq:inductive relation of n_i}
	\td{n}_i \le n_{i+1}-1,
\end{equation} 
and hence 
\begin{equation*} 
	\td{t}_i = t_{i+1} > \delta_{n_{i+1}}^{\texttt{rev}}=\widetilde{\delta}_{n_{i+1}-1}^{\texttt{rev}}	\ge \widetilde{\delta}_{\td{n}_i}^{\texttt{rev}}.
\end{equation*}
Therefore, we have ${\mathbb X}^{\texttt{tail}} \in \ov{\texttt{LR}}^{\xi'}_{\zeta'\td{\mu}'}$.
\qed
\vskip 3mm


%

\begin{lem} \label{lem:surjective}
There exists ${\bf T} \in \texttt{\em LR}^{\mu}_{\lambda}(\mf d)$ such that $\ov{\bf T}={\bf X}$, that is, 
\begin{equation*} 
\ov{\bf T}^{\texttt{\em body}} = {\bf V}, \quad \ov{\bf T}^{\texttt{\em tail}} = {\bf W}.	
\end{equation*}
\end{lem}
\pf
We use induction on $n\geq 2$. We may assume that $\mu_1' \neq 0$.

Let us first consider $n=3$. 
Note that ${\bf V} = (V_3, V_2, V_1)$, and ${\bf W}$ is a tableau of single-columned shape. 
Define $(T_2,T_1)$ by
\begin{equation*}
	(T_2^{\texttt{L}},T_2^{\texttt{R}}) = (V_3 \boxplus {\bf W},V_2), \quad T_1 = V_1.	
\end{equation*} 
Clearly $T_1$ and $T_2$ are semistandard.
By Definition \ref{bounded orthogonal LR}, the residue of $T_2$ is at most $1$. 
It is easy to check that $T_2 \prec T_1$. 
Therefore, ${\bf T} = (T_2, T_1) \in \texttt{LR}^{\mu}_{\lambda}(\mf d)$ and $\ov{\bf T}={\bf X}$.
Next, consider $n=4$. 
When $\mu_1' = 1$, apply the same argument as in the case of $n=3$. 
When $\mu_1' = 2$, we apply Lemma \ref{lem:fundamental case}.

Suppose that $n > 4$. Let us assume that $n=2l$ is even since the argument for $n$ odd is almost the same.
By Lemma \ref{lem:inductive step} and induction hypothesis, there exists ${\mathbb T}\in \texttt{LR}^{\td{\mu}}_{\xi}(\mf d)$ such that
\begin{equation*}
\ov{\mathbb T}^{\texttt{body}} = {\mathbb X}^{\texttt{body}}, \quad 
\ov{\mathbb T}^{\texttt{tail}} = {\mathbb X}^{\texttt{tail}},
\end{equation*}
where $\td{\mu}$ and ${\xi}$ are as in Lemma \ref{lem:inductive step}. 

Now, let us construct ${\bf T}=(T_l, \dots, T_1) \in \texttt{LR}^{\mu}_{\lambda}(\mf d)$ from ${\mathbb T}$, which satisfies $\ov{\bf T}={\bf X}$, 
by applying Lemma \ref{lem:fundamental case} repeatedly. 
 
Let ${\mathbb T} = (\td{T}_{l-1}, \dots, \td{T}_1,\td{T}_0)$
and let $a_i$ be the height of $\widetilde{T}_i^{\texttt{tail}}$ for $1\leq i\leq l-1$.
Put 
\begin{equation*}
{\mathbb U}=(\td{U}_{2l-1},\dots,\td{U}_2,\td{U}_1),
\end{equation*}
where
\begin{equation*}
\td{U}_{1} = \td{T}_0,\quad (\td{U}_{2i+1},\td{U}_{2i}) = (\td{T}_i^{\texttt{L}}, \td{T}_i^{\texttt{R}}) \quad (1\leq i\leq l-1).
\end{equation*} 
Let us define
\begin{equation*}
{\bf U}=(U_{2l},\dots,U_2,U_1).
\end{equation*}
First, let $U_1=\td{U}_1$ and let $U_{2l}$ be the leftmost column of ${\bf X}$. 
For $1\leq i\leq l-1$, let $(U_{2i+1},U_{2i})$ be defined in the following way.
Suppose that $a_i= 0$. Then we put
\begin{equation*}
U_{2i+1} = \td{U}_{2i+1},\quad U_{2i} = \td{U}_{2i}.
\end{equation*}
Suppose that $a_i\neq 0$. 
By Proposition \ref{prop:highest_weight_vectors}, we have $\widetilde{T}_i^{\texttt{L}}(a_i) \neq \widetilde{T}_i^{\texttt{R}}(1)$ for $1\leq i \leq l$.
If $\widetilde{T}_i^{\texttt{L}}(a_i) > \widetilde{T}_i^{\texttt{R}}(1)$, then
\begin{equation} \label{eq:shift the tails-1}
\begin{split}
U_{2i+1} =
\widetilde{T}_i^{\texttt{L}} \boxminus \widetilde{T}_i^{\texttt{tail}},  \quad
U_{2i}  =  
\widetilde{T}_i^{\texttt{R}} \boxplus \widetilde{T}_i^{\texttt{tail}}.
\end{split} 	
\end{equation}
If $\widetilde{T}_i^{\texttt{L}}(a_i) < \widetilde{T}_i^{\texttt{R}}(1)$, then
\begin{equation} \label{eq:shift the tails-2}
\begin{split}
U_{2i+1} & =
\left(\widetilde{T}_i^{\texttt{L}} \cup \left\{\widetilde{T}_i^{\texttt{R}}(1) , \widetilde{T}_i^{\texttt{L}}(a_i)\right\}\right) \boxplus \emptyset, \\
U_{2i} &  =  
\left( \td{T}_i^{\texttt{R}} \ \setminus \,  \left\{ \td{T}_i^{\texttt{R}}(1), \td{T}_i^{\texttt{R}}(2) \right\} \right) \boxplus 
\left( \left( \td{T}_i^{\texttt{tail}} \ \setminus \, \left\{ \td{T}_i^{\texttt{L}}(a_i) \right\}  \right)\cup \left\{ \td{T}_i^{\texttt{R}}(2) \right\}\right),
\end{split} 	
\end{equation}
where we identify a semistandard tableau of single-columned shape with the set of its entries.

Set
\begin{equation} \label{eq:definition of T}
	{\bf T} = (T_l, T_{l-1}, \dots, T_1), \ \text{where $(T_{i}^{\texttt{L}},T_{i}^{\texttt{R}}) = (U_{2i}, U_{2i-1})$ for $1\leq i \leq l$.}
\end{equation} 
We can check without difficulty that $T_i$ is semistandard, and  
the residue ${\mf r}_i$ of $T_i$ is at most $1$ by Lemma \ref{lem:sequence n_i} and \eqref{eq:inductive relation of n_i}. 

Next we show that $T_{i+1}\prec T_i$ and $(T_{i+1},T_i)\in {\bf H}^\circ((\mu'_{l-i},\mu'_{l-i+1}),4)$ for $1\leq i\leq l-1$, which implies that ${\bf T}\in {\bf H}^\circ(\mu,n)$. The proof is similar to the case of $n=4$ in Lemma \ref{lem:fundamental case}.

Let us prove $T_{i+1} \prec T_i$ inductively on $i$. 
For $i=1$, it follows from Lemma \ref{lem:fundamental case}.
Suppose that $T_{i}\prec T_{i-1}\prec \dots \prec T_1$ holds for given $i\geq 2$.
Consider 
\begin{equation*} %
	{\bf X}_i = (U_{2i+2}, \td{U}_{2i+1}, \td{U}_{2i},  U_{2i-1}).
\end{equation*}
By the admissibility on $\mathbb T$, ${\bf X}_i^{\texttt{tail}}$ is semistandard. 
It follows from \eqref{eq:shift the tails-1}, \eqref{eq:shift the tails-2}, Definition \ref{def:pseudo_H}(H1) on $\mathbb T$ and the induction hypothesis that ${\bf X}_i^{\texttt{body}}$ is equal to
$H_{\rho'^{\pi}}$,
for some $\rho = ({2a}, {2b}, {2c}, {2d})$ with $a \ge b \ge c \ge d \ge 0$, except the entries in the southeast corner and the next one to the left.

We remark that the map 
\begin{equation*}
{\bf X}_i=(U_{2i+2}, \td{U}_{2i+1}, \td{U}_{2i},  U_{2i-1}) \longmapsto 
(U_{2i+2}, {U}_{2i+1}, {U}_{2i},  U_{2i-1})=(T_{i+1},T_{i})
\end{equation*}
is the same as the map ${\bf X} \mapsto {\bf T}$ in Lemma \ref{lem:fundamental case}.
\vskip 3mm

{\bf\em Case 1.} {\em $\td{U}_{2i+1}(a_i) < \td{U}_{2i}(1)$ and $\td{U}_{2i+3}(a_{i+1}) < \td{U}_{2i+2}(1)$.}
First, we show that $\mf{r}_{i}=\mf{r}_{i+1}=1$.
By \eqref{eq:shift the tails-2} and \cite[Lemma 3.4]{K18-3}, we have 
\begin{equation} \label{eq:case 1 residue}
	U_{2i+2}(a_{i+1}) = \td{U}_{2i+2}(2) < \td{U}_{2i+2}(1) \le \td{U}_{2i}(1) = U_{2i+1}(1).
\end{equation}
By  \eqref{eq:shift the tails-2}, \eqref{eq:case 1 residue} and Proposition \ref{prop:highest_weight_vectors}, we have $\mf{r}_{i+1}=1$. Also, we have $\mf{r}_i=1$ by similar way.

Next, we verify Definition \ref{def:admissibility} (1)-(i), (ii) and (iii) for $(T_{i+1}, T_i)$.
The condition (1)-(i) follows from \eqref{eq:shift the tails-2}.
In this case, $T^{\texttt{R*}}_{i+1}$ and ${}^{\texttt{L}} T_i$ are given by
\begin{equation*}
	T^{\texttt{R*}}_{i+1} = \left( \td{U}_{2i+1}^{\texttt{body}} \cup \left\{\, \td{U}_{2i+1}(a_i) \, \right\} \right) \boxplus \emptyset, 
	\quad {}^{\texttt{L}} T_i = \left( \td{U}_{2i} \, \setminus \, \left\{ \, \td{U}_{2i}(1) \, \right\} \right) \boxplus \emptyset .
\end{equation*}
By Proposition \ref{prop:highest_weight_vectors} and the admissibility on $\mathbb{T}$, we have $T^{\texttt{R*}}_{i+1}(k) \le {}^{\texttt{L}} T_i(k)$. So the condition (1)-(ii) holds.

Now, we consider the condition (1)-(iii). In this case, ${}^{\texttt{R}} T_{i+1}$ and $T^{\texttt{L*}}_i$ are given by 
\begin{equation*}
\begin{split}
	{}^{\texttt{R}} T_{i+1} & = \left( \td{U}^{\texttt{body}}_{2i+1} \cup \left\{ \, \td{U}_{2i+1}(a_i) \, \right\} \right) \boxplus \left( \left( \td{U}^{\texttt{tail}}_{2i+3} \setminus \left\{ \, \td{U}_{2i+3}(a_{i+1}) \, \right\} \right) \cup \left\{ \, \td{U}_{2i}(1) \, \right\} \right) \\
	T^{\texttt{L*}}_i & =  \left( \td{U}_{2i} \, \setminus \, \left\{ \, \td{U}_{2i}(1) \,\right\} \right) \boxplus \left( \left( \td{U}_{2i+1}^{\texttt{tail}} \, \setminus \, \left\{ \, \td{U}_{2i+1}(a_i) \, \right\} \right) \cup \left\{ \, T_i^{\texttt{R}}(1) \, \right\} \right),
\end{split}
\end{equation*}
where 
\begin{equation} \label{eq:T_i^R(1)}
	T_i^{\texttt{R}}(1) =
	\left\{
\begin{array}{ll}
	\td{U}_{2i-1}(a_{i-1}+1),	& \textrm{if} \ \td{U}_{2i-1}(a_{i-1}) > \td{U}_{2i-2}(1) , \\
	\td{U}_{2i-2}(1),	& \textrm{if} \  \td{U}_{2i-1}(a_{i-1}) < \td{U}_{2i-2}(1).
\end{array}
\right.
\end{equation}
Note that we have by \cite[Lemma 3.4]{K18-3} and the admissibility of $\mathbb{T}$
\begin{equation} \label{eq:case 1 admissible condition (1)-(iii)}
	{}^{\texttt{R}} T_{i+1}(a_{i+1}) = \td{U}_{2i}(1) \le T_i^{\texttt{R}}(1) = T^{\texttt{L*}}_i(a_i).
\end{equation}
Then the condition (1)-(iii) for $(T_{i+1}, T_i)$ follows from \eqref{eq:case 1 admissible condition (1)-(iii)}, Proposition \ref{prop:highest_weight_vectors} and the admissibility of $\mathbb{T}$.

Finally, we have $(T_{i+1},T_i)\in {\bf H}^\circ((\mu'_{l-i},\mu'_{l-i+1}),4)$ by \eqref{eq:shift the tails-2}, induction hypothesis and Proposition \ref{prop:highest_weight_vectors}.

{\bf\em Case 2.} {\em $\td{U}_{2i+1}(a_i) > \td{U}_{2i}(1)$ and $\td{U}_{2i+3}(a_{i+1}) < \td{U}_{2i+2}(1)$.}
Since $\td{U}_{2i+1}(a_i) > \td{U}_{2i}(1)$, we have by the admissibility of $\mathbb{T}$
\begin{equation*}
	U_{2i+2}(a_{i+1}) = \td{U}_{2i+2}(2) < \td{U}_{2i+2}(1) \le \td{U}_{2i+1}(a_i+1) = U_{2i+1}(1).
\end{equation*}
Thus the residue $\mf{r}_{i+1}$ is equal to $1$.
If the residue $\mf{r}_i = 0$, then the admissibility of $(T_{i+1}, T_i)$ follows  immediately from the one of $\mathbb{T}$, and we have $(T_{i+1},T_i)\in {\bf H}^\circ((\mu'_{l-i},\mu'_{l-i+1}),4)$ by \eqref{eq:shift the tails-1}, \eqref{eq:shift the tails-2}, induction hypothesis and Proposition \ref{prop:highest_weight_vectors}.

We assume $\mf{r}_i = 1$. Then ${}^{\texttt{L}} T_i$, $T^{\texttt{R*}}_{i+1}$, $T^{\texttt{L*}}_i$ and ${}^{\texttt{R}} T_{i+1}$ are given by

\begin{equation} \label{eq:case 2 pairs for admissibility}
\begin{split}
	T^{\texttt{R*}}_{i+1} = \left( \td{U}^{\texttt{body}}_{2i+1} \, \setminus \, \left\{ \, \td{U}_{2i+1}(a_i+1) \, \right\} \right) \boxplus \emptyset \, ,  \quad  \quad \ \ \ 
	{}^{\texttt{L}} T_i = \left( \td{U}_{2i} \cup \left\{ \, \td{U}_{2i+1}(a_i) \, \right\} \right) \boxplus \emptyset \, , \quad \quad  & \\
	{}^{\texttt{R}} T_{i+1}	 = \left( \td{U}_{2i+1}^{\texttt{body}} \, \setminus \left\{ \, \td{U}_{2i+1}(a_i+1) \, \right\} \right) \boxplus \left( \left( \td{U}_{2i+3}^{\texttt{tail}} \, \setminus \left\{ \, \td{U}_{2i+3}(a_{i+1}) \, \right\}  \right) \cup \left\{ \, \td{U}_{2i+1}(a_i+1) \, \right\} \right), & \\
	T^{\texttt{L*}}_i = \left( \td{U}_{2i} \cup \left\{ \, \td{U}_{2i+1}(a_i) \, \right\} \right) \boxplus \left(  \left( \td{U}^{\texttt{tail}}_{2i+1} \, \setminus \, \left\{ \, \td{U}_{2i+1}(a_i) \, \right\} \right) \cup \left\{ \, T_i^{\texttt{R}}(1) \,  \right\}  \right),\quad \quad \quad \quad \quad \quad \ \ &
\end{split}
\end{equation}
where $T^{\texttt{R}}_i(1)$ is given as in \eqref{eq:T_i^R(1)}. 
By applying \cite[Lemma 3.4]{K18-3} on $\mathbb{T}$, we have 
\begin{equation} \label{eq:case 2 inequality 1}
	{}^{\texttt{R}} T_{i+1}(a_{i+1}) = \td{U}_{2i+1}(a_i+1) \le T_i^R(1) = T^{\texttt{L*}}(a_i).
\end{equation}
Now we apply a similar argument with {\em Case 1} to \eqref{eq:case 2 pairs for admissibility} with \eqref{eq:case 2 inequality 1} to obtain the admissibility of $(T_{i+1}, T_i)$ and $(T_{i+1},T_i)\in {\bf H}^\circ((\mu'_{l-i},\mu'_{l-i+1}),4)$ in this case.

{\bf\em Case 3.} {\em $\td{U}_{2i+1}(a_i) < \td{U}_{2i}(1)$ and $\td{U}_{2i+3}(a_{i+1}) > \td{U}_{2i+2}(1)$.}
The proof of this case is almost identical with {\em Case 2}. We leave it to the reader.

{\bf\em Case 4.} {\em $\td{U}_{2i+1}(a_i) > \td{U}_{2i}(1)$ and $\td{U}_{2i+3}(a_{i+1}) > \td{U}_{2i+2}(1)$.}
In this case, the claim follows immediately from \eqref{eq:shift the tails-1}, and the admissibility of $\mathbb{T}$.

Therefore, we have ${\bf T}\in {\bf H}^\circ(\mu,n)$. 
By Lemma \ref{lem:description of S}, we have ${\bf T}\equiv_{\mf l}\td{\bf T}\otimes U_{2l}$ and $\td{\bf T}=\mathbb{T}$ since ${\bf T}\in {\bf H}^\circ(\mu,n)$. This implies
\begin{equation}
{\bf T}\equiv_{\mf l} {\mathbb T} \otimes U_{2l} \equiv_{\mf l} {\mathbb X}\otimes U_{2l} \equiv_{\mf l} {\bf X},
\end{equation}
and hence ${\bf T}\in \texttt{LR}^{\mu}_{\la}(\mf{d})$. 
Since 
$\ov{\mathbb T}={\mathbb X}$, it follows from the inductive definition of $\ov{\bf T}$ that  
$\ov{\bf T} = {\bf X}$.
\qed
\vskip 3mm

{\em Proof of Theorem \ref{thm:main1} when $n-2\mu_1' \ge 0$.}
The map
\begin{equation} \label{eq:LR map}
\xymatrixcolsep{3pc}\xymatrixrowsep{0pc}\xymatrix{
\texttt{LR}^{\mu}_{\lambda}(\mathfrak{d})  \ \ar@{->}[r] & \underset{{\delta \in \cP_{n}^{(2)}}}{\bigsqcup}\ov{\texttt{LR}}^{\lambda'}_{\delta'\mu'}\\
{\bf T} \ar@{|->}[r] & \ov{\bf T}^{\texttt{tail}} }
\end{equation} 
is well-defined by Corollary \ref{cor:well-definedness}. 
Finally it is bijective by Lemmas \ref{lem:injective} and \ref{lem:surjective}.
\qed

\vskip 3mm


\subsection{Proof of Theorem \ref{thm:main1} when $n-2\mu'_1 < 0$} \label{subsec:negative case proof}
Let $\mu\in {\mc P}({\rm O}_n)$ and $\la\in \cP_{n}$ be given. 
We assume that $n-2\mu_1' < 0$.

Let ${\bf T} \in \texttt{LR}^{\mu}_{\lambda}(\mf d)$ be given with ${\bf T} = (T_l, \dots,T_{m+1},T_m,\dots,T_1,T_0)$ as in \eqref{eq:T notation-2}. 
We also use the convention for ${\bf T}$ in subsection \ref{subsec:pf_main1}.
Let $\ov{\bf T}$ be the one defined in Section \ref{subsec:negative case}.
Then we have $\ov{\bf T}^{\texttt{tail}} \in \texttt{LR}^{\lambda'}_{\mu' \delta'}$ by Proposition \ref{prop:body and tail-2}(2).
Let $\texttt{L}=2\mu'_1-n$. 
Choose $\kappa=(\kappa_{\scalebox{0.5}{$1$}},\dots,\kappa_{\scalebox{0.5}{$\texttt{L}$}})\in \cP^{(2)}$ such that $\kappa_i$ is sufficiently large.

Let $\eta, \chi\in \cP$ be given by 
\begin{equation} \label{eq:eta and xi}
\begin{split}
\eta &= \kappa \cup \la = (\kappa_{\scalebox{0.5}{$1$}},\dots,\kappa_{\scalebox{0.5}{$\texttt{L}$}},\la_1,\la_2\dots),\\
\xi &= \kappa \cup \delta =(\kappa_{\scalebox{0.5}{$1$}},\dots,\kappa_{\scalebox{0.5}{$\texttt{L}$}},\delta_1,\delta_2\dots).
\end{split}
\end{equation}

\begin{lem} \label{lem:m_i and n_i 2}
We have $\ov{\bf T}^{\texttt{\em tail}} \in \ov{\texttt{\em LR}}^{\lambda'}_{\mu' \delta'}$.
\end{lem}
\pf
Put ${\bf T} = (U_{2l}, \dots, U_{2m+1}, U_{2m}, \dots, U_0)$ under \eqref{eq:identification}. 
Let
\begin{equation*} 
\begin{split}
	\mathbb{B} & = (U_{2m}^{\downarrow}, \dots, U_0^{\downarrow}, H_{(1^{\kappa_{\scalebox{0.3}{$\texttt{L}$}}})}, \dots, H_{(1^{\kappa_{\scalebox{0.3}{$1$}}})}),
\end{split}
\end{equation*} where $U_i^{\downarrow} = (\dots, U_i(3), U_i(2)) \boxplus (U_i(1))$ for $0 \le i \le 2m$.
By the choice of $\kappa$, \,$\mathbb{B}$ is an ${\mf l}$-highest weight element,
and we note that 
\begin{equation*}
	\mu_1' = l+m+1, \quad \quad (\texttt{L}+2m+1) - 2(2m+1) = 0,
\end{equation*}
where $\texttt{L}+2m+1$ is the number of columns of $\mathbb{B}$ and $2m+1$ is the length of the first row of $\mathbb{B}^{\texttt{tail}}$.
Hence by Lemma \ref{lem:surjective}, there exists
${\bf B} = \left(X_{2m}, \dots, X_0, Y_{\texttt{L}}, \dots, Y_1\right) \in \texttt{LR}^{\dot{\mu}}_{\dot{\eta}}(\mf d)$
such that $\ov{\bf B} = \mathbb{B}$, where $\dot{\mu}' = (2m+1)$ and $\dot{\eta}$ is determined by $\mathbb{B} \equiv_{\mf l} H_{\dot{\eta}'}$.

Put
${\bf A} := \left(U_{2l}, \dots, U_{2m+1}, X_{2m}, \dots, X_0, Y_{\texttt{L}}, \dots, Y_1\right)$.
By construction of ${\bf B}$ and Corollary \ref{cor:l-equivalence under S} (cf. Remark \ref{rem:admissibilty for spin column}), it is straightforward that
\begin{equation} \label{eq:lemma 6.9 eq1}
	{\bf A} \in \texttt{LR}^{\mu}_{\eta}(\mf d), \quad \quad \ov{\bf A}^{\texttt{tail}} = \ov{\bf T}^{\texttt{tail}}.
\end{equation}

Let $(m_i)_{1 \le i \le p}$ be the sequence associated to $\ov{\bf A}^{\texttt{tail}}$, which is given as in Definition \ref{bounded orthogonal LR}.
Since by the construction $m_i \le \texttt{L}$ for all $1 \le i \le p$, 
the sequence $(m_i)_{1 \le i \le p}$ can be viewed as the sequence associated with $\ov{\bf T}^{\texttt{tail}}$ in Definition \ref{bounded orthogonal LR}.
Put $(n_i)_{1 \le i \le q}$ to be the sequence defined in Definition \ref{bounded orthogonal LR} with respect to $(m_i)_{1 \le i \le p}$. 
By Lemma \ref{lem:sequence n_i} with \eqref{eq:lemma 6.9 eq1}, the sequence $(n_i)_{1\le i \le q}$ satisfies \eqref{eq:condition_on_second_row} with respect to $\ov{\bf T}^{\texttt{tail}}$. Hence we have $\ov{\bf T}^{\texttt{tail}} \in \ov{\texttt{LR}}^{\lambda'}_{\mu' \delta'}$.
\qed

\vskip 2mm

Hence the map \eqref{eq:LR map} is well-defined by Proposition \ref{prop:body and tail-2}(2) and Lemma \ref{lem:m_i and n_i 2}. 
It is also injective since Lemma \ref{lem:injective} still holds in this case. 
So it remains to verify that the map is surjective.

Let ${\bf W} \in \ov{\texttt{LR}}^{\lambda'}_{\mu' \delta'}$ be given for some $\delta \in \cP_n^{(2)}$.
Let ${\bf V}= H_{(\delta')^{\pi}}$ and ${\bf X}$ be the tableau of a skew shape $\eta$ as in \eqref{eq:shape after separation} with $n$ columns such that 
${\bf X}^{\texttt{body}}={\bf V}$ and ${\bf X}^{\texttt{tail}}={\bf W}$.
As in the case of $n-2\mu'_1\geq 0$, ${\bf X}$ is semistandard.

Put
\begin{equation} \label{eq:construction of Y and Z}
\begin{split}
{\bf Y} & =(Y_{\texttt{L}}, \dots, Y_1),\\
{\bf Z} & = (X_n,\dots,X_1,Y_{\texttt{L}}, \dots, Y_1),
\end{split}
\end{equation} 
where $Y_i = H_{(1^{\kappa_i})}$ for $1\leq i \leq \texttt{L}$.

\begin{lem}\label{lem:Z LR}
We have ${\bf Z}^{\texttt{\em tail}}\in \ov{\texttt{\em LR}}^{\eta'}_{\mu' \xi'}$.
\end{lem}
\pf By construction of ${\bf Z}$, we have ${\bf Z}\equiv_{\mf l} H_{\eta'}$.
Let $(m_i)_{1\leq i\leq p}$ and $(n_i)_{1\leq i\leq q}$ be the sequences associated to ${\bf W}\in \ov{\texttt{LR}}^{\lambda'}_{\mu' \delta'}$.
Since $\kappa_i$ is sufficiently large, we have
${\bf Z}^{\texttt{tail}}\in \ov{\texttt{LR}}^{\eta'}_{\mu' \xi'}$ with respect to the same sequences $(m_i)_{1\leq i\leq p}$ and $(n_i)_{1\leq i\leq q}$.
\qed
\vskip 2mm

Note that ${\bf Z}\in {\bf E}^{M}$ where $M=n+\texttt{L}=2\mu'_1$.
By Lemma \ref{lem:Z LR}, we may apply Theorem \ref{thm:main1} for $M-2\mu'_1=0$ to conclude that there exists a unique ${\bf R}\in {\bf T}(\mu,M)$ such that
\begin{equation}\label{eq:+inductive step 1}
\ov{\bf R} = {\bf Z}.
\end{equation}
Suppose that ${\bf R}=(R_{M},\dots,R_1)\in {\bf E}^{M}$ under \eqref{eq:identification}.
Put
${\bf S} = (R_{2\texttt{L}}, \dots, R_1)$. 
Note that $2\texttt{L}<M$ with $M-2\texttt{L}=n$, and ${\bf S}\in {\bf T}((1^\texttt{L}),2\texttt{L})$.
If $\ov{\bf S} = (\ov{S}_{2\texttt{L}}\dots, \ov{S}_1)\in {\bf E}^{2\texttt{L}}$, 
then we have by Corollary \ref{cor:l-equivalence under S} and \eqref{eq:+inductive step 1} 
\begin{equation*} \label{eq:U}
(\ov{S}_{\texttt{L}}\dots, \ov{S}_1) = {\bf Y}. 
\end{equation*}

Now, we put
\begin{equation*}
{\bf T} = (R_M,\dots,R_{2\texttt{L}+1}, \ov{S}_{2\texttt{L}},\dots,\ov{S}_{\texttt{L}+1}) \in {\bf E}^n,
\end{equation*}
under \eqref{eq:identification}.
Then it is straightforward to check that ${\bf T} \in {\bf T}(\mu,n)$.
Since ${\bf Z}\in \texttt{LR}^{\mu}_{\eta}(\mf d)$, we have ${\bf T} \in \texttt{LR}^{\mu}_{\lambda}(\mf d)$ by construction of ${\bf T}$ and Lemma \ref{lem:criterion_highest_weight_elt}. Finally, by \eqref{eq:+inductive step 1} and Corollary \ref{cor:l-equivalence under S}, we have 
\begin{equation*}
\ov{\bf T}^{\texttt{body}} = {\bf X}^{\texttt{body}}, \quad 
\ov{\bf T}^{\texttt{tail}} = {\bf X}^{\texttt{tail}}.
\end{equation*} Hence, the map \eqref{eq:LR map} is surjective.
\qed

\appendix
\section{Index of notation}
\begin{itemize}
	\item[\ref{subsec:crystals} :] 
	$\equiv_{\mf g}$, $\equiv$
	\item[\ref{subsec:notations} :] 
	$\cP$, $\cP_{\ell}$, $\cP^{(2)}$, $\cP^{(1,1)}$, $\cP^{(2,2)}$, $\cP^{(2)}_{\ell}$, $\cP^{(1,1)}_{\ell}$, $\cP^{(2,2)}_{\ell}$, $\ell(\lambda)$, $SST(\lambda/\mu)$, $w(T)$, $a \rightarrow T$, $w \rightarrow T$, $S \rightarrow T$, $\la^{\pi}$, $H_{\la}$, $H_{\la^\pi}$, $\texttt{LR}^{\lambda}_{\mu \nu}$, $\texttt{LR}^{\lambda}_{\mu \nu^\pi}$, $c^{\lambda}_{\mu \nu}$, $\psi : \texttt{LR}^{\lambda'}_{\mu' \nu'} \rightarrow \texttt{LR}^{\lambda}_{\mu \nu^\pi}$, $S^i$, $H^i$, $Q(S \rightarrow H_{\mu'})$, $U_i$, $H_i$, $Q(U \rightarrow H_{\mu})$
	\item[\ref{subsec:spinor} :]
	$\mathfrak{l}$, ${\rm ht}(T)$, $T^{\texttt{body}}$, $T^{\texttt{tail}}$, $U(i)$, $U[i]$, $\boxplus$, $\boxminus$, $\lambda(a, b, c)$, ${\mf r}_T$, $\mathcal{E}$, $\mathcal{F}$, ${\bf T}(a)$, $\ov{\bf T}(0)$, ${\bf T}^{\texttt{sp}}$, ${\bf T}^{\texttt{sp}+}$, ${\bf T}^{\texttt{sp}-}$, $T^{\texttt{L}*}$, $T^{\texttt{R}*}$, ${}^{\texttt{L}} T$, ${}^{\texttt{R}} T$, $T \prec S$, $\mathcal{P}({\rm O}_n)$, $\Lambda(\mu)$, $q_{\pm}$, $r_{\pm}$, $\ov{\mu}$, $\widehat{\bf T}(\mu, n)$, ${\bf T}(\mu, n)$
	\item[\ref{subsec:l-highest weight element} :]
	${\bf H}(\mu, n)$, ${\bf H}^{\circ}(\mu, n)$
	\item[\ref{subsec:sliding} :] ${\bf E}^n$, $\mathcal{E}_j$, $\mathcal{F}_j$, $\mathcal{S}_j$	%
	\item[\ref{subsec:separation} :]
	$\widetilde{\bf T}$, $\ov{\bf T}$, $\ov{\bf T}^{\texttt{body}}$, $\ov{\bf T}^{\texttt{tail}}$
	\item[\ref{subsec:combi_branching} :]
	$\texttt{LR}^{\mu}_{\la}({\mf d})$, $c^{\mu}_{\la}({\mf d})$, $\delta^{\texttt{rev}}$, 
	$m_i$, $n_j$,
	$\ov{\texttt{LR}}^{\lambda'}_{\delta' \mu'}$, $\ov{c}^{\la}_{\delta \mu}$, 
	$\underbar{\texttt{LR}}^{\lambda}_{\delta \mu}$, $\underline{c}^{\la}_{\delta \mu}$,
	\item[\ref{subsec:generalized exponent} :]
	$V^{\mu}_{\mf g}$, $K^{\mf g}_{\mu 0}(t)$
	\item[\ref{subsec:combi model of generalized expoentns} :]
	$\varphi(T)$, $\varepsilon(T)$, $D_n(\mu)$, $\cP_T$, $\rho_T$, $\underbar{$D$}_n(\mu)$, $\mathbb{D}_n(\mu)$, 
\end{itemize}


\end{document}